\documentclass[12pt, pdftex]{amsart}
\usepackage[margin=1in]{geometry}
\usepackage[utf8]{inputenc}
\usepackage{tikz}
\usetikzlibrary{arrows.meta,decorations.pathreplacing,decorations.markings,shapes,calc}
\usetikzlibrary{matrix,arrows,decorations.pathmorphing,backgrounds,decorations.markings,positioning}
\usepackage{graphics}
\newcommand{\rotatesim}{\rotatebox[origin=c]{90}{$\sim$}}
\usepackage[all]{xy}
\usepackage{comment}
\usepackage{amsmath}
\usepackage{amssymb}
\usepackage{amsthm}
\usepackage{here}
\usepackage{amscd} 
\usepackage{tikz-cd}
\usepackage{mathrsfs}
\usepackage{mathtools}
\usepackage{mathabx}
\usepackage{scalefnt}
\usepackage{url}
\theoremstyle{plain}
\newtheorem{thm}{Theorem}[section]
\newtheorem{lem}[thm]{Lemma}
\newtheorem{cor}[thm]{Corollary}
\newtheorem{prop}[thm]{Proposition}
\newtheorem{conj}[thm]{Conjecture}
\theoremstyle{definition}

\newtheorem{rem}[thm]{Remark}

\newtheorem{que}[thm]{Question}

\newtheorem{ex}[thm]{Example}

\numberwithin{equation}{section}
\def\A{{\mathbb A}}

\def\Q{{\mathbb Q}}
\def\R{{\mathbb R}}
\def\Z{{\mathbb Z}}
\def\C{{\mathbb C}}
\def\P{{\mathbb P}}
\def\B{{\mathbb B}}

\def\Aut{\mathop{\mathrm{Aut}}\nolimits}

\def\Ker{\mathop{\mathrm{Ker}}\nolimits}
\def\Pic{\mathop{\mathrm{Pic}}\nolimits}

\def\SL{\mathop{\mathrm{SL}}\nolimits}

\def\det{\mathop{\mathrm{det}}\nolimits}
\def\dim{\mathop{\mathrm{dim}}\nolimits}

\def\Stab{\mathop{\mathrm{Stab}}\nolimits}

\def\ss{\mathop{\mathrm{ss}}}
\def\c{\mathop{\mathrm{c}}}

\def\ord{\mathrm{ord}}
\def\GIT{\mathrm{GIT}}
\def\K{\mathrm{K}}
\def\tor{\mathrm{tor}}
\def\BB{\mathrm{BB}}

\def\L{\mathscr{L}}
\def\H{\mathscr{H}}
\def\g{\mathfrak{g}}

\def\N{\mathscr{N}}
\def\M{\mathcal{M}}

\def\OO{\mathscr{O}}

\def\A{\mathcal{A}}
\def\D{\mathscr{D}}

\def\a{\alpha}
\def\b{\beta}
\def\g{\gamma}
\def\d{\delta}
\def\e{\epsilon}
\def\l{\langle}
\def\r{\rangle}

\def\SO{\mathop{\mathrm{SO}}\nolimits}

\def\Sym{\mathop{\mathrm{Sym}}\nolimits}

\def\U{\mathrm{U}}
\def\O{\mathrm{O}}

\allowdisplaybreaks[4]

\newcommand{\defeq}{\vcentcolon=}

\begin{document}

\title[Compactifications of the Eisenstein ancestral Deligne-Mostow variety]
{Compactifications of the ancestral Eisenstein Deligne-Mostow variety}

\author{Klaus Hulek}
\address{K.H: Institut f\"ur Algebraische Geometrie, Leibniz University Hannover, Welfengarten 1, 30060 Hannover, Germany}
\email{hulek@math.uni-hannover.de}
\author{Shigeyuki Kond\={o}}
\address{S.K: Graduate School of Mathematics, Nagoya University, Nagoya, 464-8602, Japan}
\email{kondo@math.nagoya-u.ac.jp}
\author{Yota Maeda}
\address{Y.M: Advanced Research Laboratory, Technology Infrastructure Center, Technology Platform, Sony Group Corporation, 1-7-1 Konan, Minato-ku, Tokyo, 108-0075, Japan}
\email{y.maeda.math@gmail.com}

\date{\today}

\maketitle
\begin{abstract}
All arithmetic non-compact ball quotients by Deligne-Mostow's unitary monodromy group arise as sub-ball quotients of either of two spaces called ancestral cases, corresponding to Gaussian or Eisenstein Hermitian forms respectively.
In a previous paper, we investigated the compactifications of the Gaussian Deligne-Mostow variety.
Here we work on the remaining case, namely the ring of Eisenstein integers. This variety is related to the moduli space of unordered 12 points on $\P^1$.
In particular, we show that Kirwan's partial resolution of the moduli space is not a semi-toroidal compactification and Deligne-Mostow's period map does not lift to the unique toroidal compactification. 
We give two interpretations of these phenomena in terms of the log minimal model program and automorphic forms.
As an application, we prove that the above two compactifications are not (stacky) derived equivalent, as the $DK$-conjecture predicts.
Furthermore, we construct an automorphic form on the moduli space of non-hyperelliptic curves of genus 4, which is isogenous to the Eisenstein Deligne-Mostow variety,  giving another intrinsic proof, independent of lattice embeddings, of a result by Casalaina-Martin, Jensen and Laza.
\end{abstract}
\tableofcontents

\section{Introduction}
\subsection{Compactifications of $\M_{g,n}$}
The coarse moduli space $\M_{g,n}$ of smooth ordered distinct $n$-pointed curves of genus $g$ is one of the most studied objects in algebraic geometry.
This is a quasi-projective variety, which, however, is not projective as both conditions {\em smooth} and {\em distinct} are open conditions that are not preserved under (natural) degenerations.  
This poses the obvious problem to try and understand the limits of $n$-pointed curves in families.
This in turn leads to the notion of a compactification ${\M}^{\c}_{g,n}$ for the moduli spaces $\M_{g,n}$.
Obviously, such compactifications are not uniquely determined and one wants to single out those which are well-behaved from the point of view of the moduli problem. 
In particular, one asks that a compactification ${\M}^{\c}_{g,n}$ should
\begin{enumerate}
    \item[(A)] be a moduli space itself, extending the moduli problem $\M_{g,n}$ naturally,
    \item[(B)] have mild singularities, and
    \item[(C)] have a normal crossing boundary ${\M}^{\c}_{g,n}\setminus\M_{g,n}$.
\end{enumerate}
In the case of $\M_g$, the most famous and best-investigated solution is the Deligne-Mumford compactification $\overline{\M}_{g}$.
The situation is more complicated when one considers $n$-pointed curves.
In general, the Deligne-Mumford compactification $\overline{\M}_{g,n}\supset\M_{g,n}$ always exists. Here the underlying limit curves are nodal curves such that the $n$-pointed curves have 
finite automorphism groups.
This satisfies the above conditions (A), (B) and (C).
However, studying possible limits of pointed curves $(C,p_1,\cdots,p_n)\in\M_{g,n}$, one is also led to other natural compactifications, depending on the point of view one takes. Understanding this in some depth is already an interesting problem in genus $0$ with many questions unresolved. This naturally leads to the moduli space of $n$ distinct ordered points on $\P^1$.
A very natural approach is to take the GIT quotient $\M_{\ord, n}^{\GIT}\supset \M_{0,n}$, namely
\[\M_{\ord, n}^{\GIT}\defeq (\P^1)^n /\!/_{\OO(1,\cdots,1)}\SL_2(\C).\]
Roughly speaking, this parametrizes the tuples $(p_1,\cdots,p_n)$ where at most $\lfloor n/2 \rfloor$ points coincide. This satisfies property (A), but not necessarily (B) and (C). It is an obvious question to ask what the difference is between 
$\overline{\M}_{0,n}$ and $\M_{\ord, n}^{\GIT}$. Many mathematicians have worked on this, often using the log minimal model (LMMP) as a guiding principle. This includes:
\begin{enumerate}
    \item[(A)] Keel, Kapranov and Losev-Manin studied the birational contraction morphism $\overline{\M}_{0,n} \to \M_{\ord,n}^{\GIT}$, which can be described by combinatorial data.
    \item[(B)] Hassett \cite[Section 2]{Has03} introduced the moduli space of weighted pointed stable rational curves $\overline{\M}_{0,\A}$ for a weight datum $\A\in\Q^n$.
    The compactifications $\overline{\M}_{0,\A}$ generalize $\overline{\M}_{0,n}$ and $\M_{\ord,n}^{\GIT}$, providing a whole hierarchy of compactifications of $\M_{0,n}$.
    \item[(C)] Kiem-Moon \cite[Theorem 1.1]{KM11} studied the factorization
\begin{align}
\label{morphism:blow-up_sequence_ordered_npoints}
    \overline{\M}_{0,n} \stackrel{\varphi_{\lfloor n/2 \rfloor-2}}{\to} \overline{\M}_{0,\A_{\lfloor n/2 \rfloor-3}} \stackrel{\varphi_{\lfloor n/2 \rfloor-3}}{\to} \cdots \stackrel{\varphi_2}{\to}\overline{\M}_{0,\A_1} \stackrel{\varphi_1}{\to} \M^{\GIT}_{\ord,n}
\end{align}
of the above contraction for symmetric linearizations.
Here, the weight data $\A_i$ is given by the symmetric $n$-tuple with value $1/(\lfloor n/2\rfloor+1-i)+\epsilon$ 
for any small $\epsilon>0$.
Below, we simply write this as $\A_i=n(1/(\lfloor n/2\rfloor+1-i)+\epsilon)$.
\end{enumerate}
Kiem and Moon used this to determine the log canonical models of $\overline{\M}_{0,n}$, and along the way, also obtained the following result.
\begin{thm}[{\cite[Theorem 1.1]{KM11}}]
\label{previous_thm:kirwan_ordered}
When $n$ is even, the  blow-up $\overline{\M}_{0,\A_1} \to \M_{\ord,n}^{\GIT}$ coincides with the Kirwan desingularization  $\M^{\K}_{\ord,n}$ of 
$\M_{\ord,n}^{\GIT}$ at the $\frac{1}{2}\binom{n}{n/2}$ polystable points:
\begin{equation}\label{equ:KirwanequalHassett}
\M^{\K}_{\ord,n}\cong\overline{\M}_{0,\A_1}.
\end{equation}
\end{thm}

It is also natural to consider the case of unordered $n$ points on $\P^1$. The GIT moduli space is then given by 
\[\M_n^{\GIT}\defeq\P^n/\!/_{\OO(1)}\SL_2(\C)\]
where the action of $\SL_2(\C)$ on $\P^n$ is descended from the diagonal action of $\SL_2(\C)$ on $(\P^1)^n$ via the isomorphism $ (\P^1)^n/S_n \cong \P^n$.  
By taking the $S_n$-quotient of (\ref{morphism:blow-up_sequence_ordered_npoints}),  we obtain a blow-up sequence for the moduli space of unordered $n$ points on $\P^1$, namely
\[
    \widetilde{\M}_{0,n} \to \widetilde{\M}_{0,\A_{\lfloor n/2 \rfloor-3}} \to \cdots \to\widetilde{\M}_{0,\A_1} \to \M_n^{\GIT},\]
where $\widetilde{\M}_{0,n}\defeq\overline{\M}_{0,n}/S_n$ and $\widetilde{\M}_{0,\A_i}\defeq\overline{\M}_{0,\A_i}/S_n$.
Here the spaces $\widetilde{\M}_{0,n}$ and $\widetilde{\M}_{0,\A_i}$ are the moduli space of unordered distinct $n$ points on $\P^1$ and unordered weighted pointed stable rational curves respectively, which are compactifications of $\M_{0,n}/S_n$.
The motivation for the work of Kiem-Moon was to prove a conjecture by Hassett \cite[Acknowledgment]{KM11}, which says  that, for $n$ even, $\widetilde{\M}_{0,\A_1}$ is the (weighted) blow-up of $\M^{\GIT}_{n}$
at the unique singular (polystable) point. Indeed, Theorem \ref{previous_thm:kirwan_ordered} of Kiem and Moon confirms this for the case of ordered points. This is one of the starting points of our discussion. 
The above discussion provides one solution that characterizes the blow-up, but from the perspective of GIT quotients, it is also natural to consider the following question, which can be seen as an analog of Theorem \ref{previous_thm:kirwan_ordered}.
\begin{que}
\label{que:kirwan_blow-up}
    When $n$ is even, does the blow-up $\widetilde{\M}_{0,\A_1}\to\M^{\GIT}_n$ coincide with the Kirwan partial desingularization of $\M^{\GIT}_n$ at the unique polystable point?
\end{que}

\subsection{Main results}
\label{subsec:main_results}
There is also a completely different angle from which one can approach Question \ref{que:kirwan_blow-up}:
for certain values of $n$, the space $\M^{\GIT}_n$ also has a realization as a ball quotient due to the Deligne-Mostow theory \cite{DM86, Mos86}.
This means that there are various period maps 
\[\Phi_{\mathbf{w}}:\M_{0,n}/S_n\hookrightarrow \B^{n-3}/\Gamma_{\mathbf{w}}\]
for an arithmetic non-cocompact unitary group $\Gamma_{\mathbf{w}}$ stabilizing certain Hermitian forms, according to a weight $\mathbf{w}\in\Q^n$ listed in \cite[Table 2, 3]{GKS21}.
This applies in particular to $n=8$ and $n=12$. These two cases are called the {\em ancestral} cases. The reason for this is that all other Deligne-Mostow varieties, 
where the monodromy group of Appell-Lauricella hypergeometric functions is arithmetic, and the corresponding GIT quotient has strictly semistable points, 
can be realized as sub-ball quotients of these two varieties. 
The case $n=8$ is the maximal ball quotient where the Hermitian form is defined over the Gaussian integers, whereas the case $n=12$ is maximal for the Eisenstein integers. We refer the reader to
\cite[Section 2]{GKS21} for a more detailed discussion of this topic. 
For the case $n=8$ we showed in \cite{HM25} that Question \ref{que:kirwan_blow-up} has indeed a negative answer. 

The current paper is dedicated to the remaining ancestral case, namely the moduli space of unordered 12 points on $\P^1$, which is 
related to the ring of Eisenstein integers. 
For simplicity, we will from now on omit $n=12$ in our notation, that is, $\M_{\ord}^{\GIT}$ means $\M_{\ord, 12}^{\GIT}$, etc.
We will describe the relevant arithmetic group $\Gamma$ acting on $\B^9$ in detail in Subsection \ref{subsection:DMiso}. 
We also recall that any arithmetic ball quotient $\B^n/\Gamma$ has two natural algebraic compactifications, namely the Baily-Borel compactification 
$\overline{\B^n/\Gamma}^{\BB}$ and the toroidal compactification
$\overline{\B^n/\Gamma}^{\tor}$. The latter is unique for ball quotients as the relevant fan lives in a 1-dimensional space and this allows no choices. 
As the above discussion makes it clear, the ball quotient $\B^9/\Gamma$, which contains the moduli space of 12 unordered points, has several 
interesting compactifications. In the present paper we will mainly concentrate on 
the Kirwan blow-up $\M^{\K}$, the moduli space of unordered weighted pointed stable rational curves $\widetilde{\M}_{0,\A_1}$ and the toroidal compactification $\overline{\B^9/\Gamma}^{\tor}$. 
Obviously, it is natural to ask how these relations are related to each other. To simplify the terminology we say that two compactifications of (an open subset of) a ball quotient are 
{\em naturally isomorphic} if the identity extends to an isomorphism of the compactifications.
We first recall that Gallardo, Kerr and Schaffler clarified the relationship between the spaces $\widetilde{\M}_{0,\A_1}$ and the toroidal compactification $\overline{\B^9/\Gamma}^{\tor}$.
\begin{thm}[{\cite[Theorem 1.1]{GKS21}}]
\label{teo:GKSisomunordered}
The two compactifications  $\widetilde{\M}_{0,\A_1}$ and  $\overline{\B^9/\Gamma}^{\tor}$  are naturally  isomorphic to each other.
\end{thm}
Our first main theorem relates the Kirwan blow-up $\M^{\K}$  of  $\M_{\ord}^{\GIT}$ and the toroidal compactification $\overline{\B^9/\Gamma}^{\tor}$ of the Deligne-Mostow variety.
\begin{thm}[{Theorems \ref{thm:lift}, \ref{thm:not_k_equiv}, \ref{thm:not_semi_toric}}, Remark \ref{rem:another_proof_semi_toroidal}]
    \label{mainthm:semi-toroidal}
    Regarding the geometry of $\M^{\K}$ and $\overline{\B^9/\Gamma}^{\tor}$, the following holds.
    \begin{enumerate}
        \item The Kirwan blow-up $\M^{\K}$ is not a semi-toroidal compactification.
        \item Neither the Deligne-Mostow isomorphism $\phi$ nor its inverse $\phi^{-1}$ lift to a morphism between the Kirwan blow-up $\M^{\K}$ and the unique toroidal compactification $\overline{\B^9/\Gamma}^{\tor}$.
\item The varieties $\M^{\K}$  and $\overline{\B^9/\Gamma}^{\tor}$ are not $K$-equivalent and hence, in particular, not isomorphic as abstract varieties.
    \end{enumerate}
\end{thm}
Combining Theorems \ref{teo:GKSisomunordered} and \ref{mainthm:semi-toroidal} (1) gives a negative answer to 
Question \ref{que:kirwan_blow-up} in the unordered case.
\begin{cor}
    \label{maincor:not_isom}
     The two compactifications $\widetilde{\M}_{0,\A_1}$ and $\M^{\K}$ of the ball quotient $\B^9/\Gamma$ are not isomorphic to each other. 
\end{cor}
Moreover, in the process of showing Theorem \ref{mainthm:semi-toroidal} (3), we also obtain the following more intuitive result from the LMMP point of view.

\begin{thm}[{Theorem \ref{thm:not_semi_toric}}]
\label{mainthm:not_log_crepant}
The two compactifications $\M^{\K}$ and $\overline{\B^9/\Gamma}^{\tor}$ are not log $K$-equivalent when taking the boundaries as the sum of the strict transform of the discriminant divisors with its standard coefficients and the exceptional divisors with coefficient 1.
\end{thm}

Theorem \ref{mainthm:semi-toroidal} 
can be approached from (at least) two angles.
We will give a first proof of Theorem \ref{mainthm:semi-toroidal}  in Theorems  \ref{thm:lift} and \ref{thm:not_k_equiv}. This is based on Luna slice computations and automorphic forms, see Subsection \ref{subsec:auxiliary_calculations}.
This proof is motivated by the aim of understanding more detailed geometric phenomena by examining the moduli spaces $\overline{\M}_{0,\A}$, as introduced by Hassett \cite{Has03}, and also uses 
automorphic forms \cite{All00}.
In Subsection \ref{subsection:semi-toroidal}, we will further give a proof of Theorem \ref{mainthm:semi-toroidal} (1)
in terms of LMMP. This uses the work of Alexeev-Engel and of Odaka \cite{AE23, AEH24, Oda22} who 
characterize semi-toroidal compactifications by the property that they lie between the (in this case unique) 
toroidal compactification and the Baily-Borel compactification. Their work uses concepts from LMMP and recent progress in $K$-stability.

We also discuss the relationship between the different compactifications $\M^{\K}$ and $\overline{\B^9/\Gamma}^{\tor}$ of the ball quotient $\B^9/\Gamma$ from the perspective of derived geometry.  
Since $K$-equivalence is expected to be reflected by derived categories via the $DK$-conjecture, it is natural to regard the categorical comparison of $\M^K$ and $\B^9/\Gamma^{\tor}$ as a natural endpoint of our study of these geometric compactifications. Indeed, one of the most fundamental questions in derived algebraic geometry, historically proposed by Bondal-Orlov \cite{BO02} and Kawamata \cite{Kaw02}, is the \textit{$DK$-conjecture}, see \cite[Conjecture 1.2]{Kaw02}.
This asserts that for two smooth projective varieties, $D$-equivalence is equivalent to $K$-equivalence.
This was generalized to projective varieties with only finite quotient singularities by Kawamata \cite[Conjecture 2.2]{Kaw05}, \cite[Conjecture 4.1]{Kaw18}.
By putting the boundaries as 0 there, one can conjecture the following.
\begin{conj}[{\cite[Conjecture 2.2]{Kaw05}, \cite[Conjecture 4.1]{Kaw18}}]
\label{mainconj:generalized DK conjecture}
    Let $X$ and $Y$ be projective varieties with only finite quotient singularities,
and let $\mathcal{X}(X)$ and $\mathcal{X}(Y)$ be the associated stacks (see below for the precise definition and \cite{Kaw04, Kaw05}).
Then, the varieties $X$ and $Y$ are $K$-equivalent if and only if there exists an equivalence as triangulated categories between the bounded derived categories of coherent sheaves on $\mathcal{X}(X)$ and $\mathcal{X}(Y)$:
\[D^{\mathrm{b}}(\mathrm{{Coh}}(\mathcal{X}(X)))\cong D^{\mathrm{b}}(\mathrm{{Coh}}(\mathcal{X}(Y))).\]
\end{conj}
For the precise notion, see \cite{Kaw02, Kaw04, Kaw05,Kaw18} and Subsection \ref{subsec:concepts}.
It is known that the converse direction (from $D$ to $K$) is not true in general \cite[Remark 2.5]{Kaw18}.
However, as an application of Theorem \ref{mainthm:semi-toroidal}, we, in fact, prove the following, which is the converse in Conjecture \ref{mainconj:generalized DK conjecture} in our case.
Below, for a projective variety $X$ with only finite quotient singularities, we denote by $D(\mathcal{X}(X))$ the derived category of the coherent sheaves on the associated stack $\mathcal{X}(X)$; see Subsection \ref{subsec:concepts} for more details.
\begin{thm}[{Theorem \ref{thm:derived_equivalent}}]
\label{mainthm:derived_equivalent}
The bounded derived categories of coherent sheaves on the stacks associated with the two varieties $\M^{\K}$ and $\overline{\B^9/\Gamma}^{\tor}$ are not equivalent as triangulated categories:
\[D(\mathcal{X}(\M^{\K}))\not\cong D(\mathcal{X}(\overline{\B^9/\Gamma}^{\tor})).\]   
In other words, the two compactifications $\M^{\K}$ and $\overline{\B^9/\Gamma}^{\tor}$ are not stacky $D$-equivalent.
\end{thm}
This result seems interesting because there is very little research on the derived algebraic geometry of modular varieties.
Moreover, Conjecture \ref{mainconj:generalized DK conjecture} is, to the best of the authors' knowledge, partially known only for toric or toroidal varieties, first treated in \cite[Theorem 4.2]{Kaw05}. 
This is due to a series of papers by Kawamata  \cite{Kaw06, Kaw13, Kaw16}.
To avoid confusion we remark that in Kawamata's definition \cite[Definition 4.1]{Kaw05}, all points of a toroidal variety must have toric neighborhoods, a property which is typically not satisfied by 
toroidal compactifications of locally symmetric spaces.
To prove the above theorem, we use recently developed LMMP methods and classical results on automorphic forms by Mumford.
In the case we know that the Baily-Borel compactification has finite quotient singularities, such as the moduli spaces of cubic surfaces, we can moreover prove a similar assertion for these spaces.
\begin{thm}[{Theorem \ref{thm:derived_cubic}}]
\label{mainthm:derived_cubic}
    For the moduli spaces of cubic surfaces, no two of the three categories $D(\mathcal{X}(\M_{\mathrm{cub}}^{\K}))$, $D(\mathcal{X}(\overline{\B^9/\Gamma_{\mathrm{cub}}}^{\tor}))$ and $D(\mathcal{X}(\overline{\B^9/\Gamma_{\mathrm{cub}}}^{\BB}))$ are equivalent to each other as triangulated categories.    
In other words, the three compactifications $\M_{\mathrm{cub}}^{\K}$, $\overline{\B^9/\Gamma_{\mathrm{cub}}}^{\tor}$ and $\overline{\B^9/\Gamma_{\mathrm{cub}}}^{\BB}$ are mutually not stacky $D$-equivalent.
\end{thm}

In this paper, we also compute the cohomology of the various spaces which appear in our discussions. This was indeed the starting point of some of this work in recent years.
The first named author, together with Casalaina-Martin, Grushevsky and Laza, noticed that for both, moduli  of cubic surfaces and cubic threefolds,
the Betti numbers of the Kirwan blow-up 
and the toroidal compactifications of the corresponding ball quotients coincide \cite{CMGHL23a}. 
This naturally led to the question of how precisely these spaces are related. 
In the case of 8 points on $\P^1$ the analogous question was treated in \cite{HM25}. 
In our case, the result is:
\begin{thm}[{Theorems \ref{thm:coh_previous_work}, \ref{thm:coh_M^K}, \ref{thm:coh_tor}}]
All the odd degree cohomology of the following projective varieties vanishes.
In even degrees, their Betti numbers are given by:
\begin{align*}
\renewcommand*{\arraystretch}{1.2}
\begin{array}{l|cccccccccc}
\hskip2cm j&0&2&4&6&8&10&12&14&16&18\\\hline
\dim H^j(\M_{\ord}^{\K})&1&474&991&1618&2410&2410&1618&991&474&1\\
\dim IH^j(\M_{\ord}^{\GIT})&1&12&67&232&562&562&232&67&12&1\\
\dim H^j(\M^{\K})&1&2&3&4&5&5&4&3&2&1\\
\dim IH^j(\M^{\GIT})&1&1&2&2&3&3&2&2&1&1\\
\dim H^j(\overline{\B^9/\Gamma}^{\tor})&1&2&3&4&5&5&4&3&2&1\\
\dim IH^j(\overline{\B^9/\Gamma}^{\BB})&1&1&2&2&3&3&2&2&1&1
\end{array}
\end{align*}
In particular, all  Betti numbers of $\M^{\K}$ and $\overline{\B^9/\Gamma}^{\tor}$ coincide.
\end{thm}
We recall that $\M^{\GIT} \cong \overline{\B^9/\Gamma}^{\BB}$ and hence these spaces must obviously have the same Betti numbers.

In Section \ref{section:Automorphic forms} we shall slightly change the perspective and complete the picture by investigating two related automorphic forms, 
one of which was found by Allcock in \cite{All00}. 
In particular, we construct an automorphic form on another important 9-dimensional ball quotient $\B^9/\Gamma_{\mathrm{nh}}$.
To put this into perspective we first recall some facts about the moduli space of cubic threefolds. 
Allcock, Carlson and Toledo \cite{ACT11}, and Looijenga and Swiestra \cite{LS07}, realized this moduli space as a ball quotient.
More precisely, the moduli space of smooth cubic threefolds is isomorphic to an open subset of a 10-dimensional ball quotient $\B^{10}/\Gamma'$, where the complement consists of two irreducible Heegner divisors, namely
 $H_{\mathrm{n}}'$  (the locus of singular cubics, also referred to as the nodal or discriminant locus) and $H_{\mathrm{h}}'$ (the locus of chordal cubics, also referred to as the hyperelliptic locus).
 Both of these divisors are closely related to 9-dimensional ball quotients, namely $\B^9/\Gamma_{\mathrm{nh}}$ and the ball quotient $\B^9/\Gamma$ which we have been considering so far and which contains the 
 moduli space of 12 unordered points or, equivalently, the moduli space of hyperelliptic curves of genus 5.
The ball quotient $\B^9/\Gamma_{\mathrm{nh}}$ (see below and Subsection \ref{subsec:A complex ball uniformization} for more details) maps to the discriminant locus $H_{\mathrm{n}}'$ which is a 
 normalization onto the image by \cite[Proposition 5.3]{CMJL12}. It contains the moduli space $\mathcal M_4$ of curves of genus 4 as an open subset.
When the discriminant locus $H_{\mathrm{n}}'$ is restricted to $\B^9/\Gamma_{\mathrm{nh}}$, it decomposes into two divisors $H_{\mathrm{n}}=\H_{\mathrm{n}}/\Gamma_{\mathrm{nh}}$ and 
$H_{\mathrm{vt}}=\H_{\mathrm{vt}}/\Gamma_{\mathrm{nh}}$ which in turn are ball quotients, as we shall discuss in Subsection \ref{subsec:A complex ball uniformization}. 
We remark here that the birational isomorphism between $\B^9/\Gamma_{\mathrm{nh}}$ and $\M_4$ was given by using the theory of periods of $K3$ surfaces.

Accordingly, in Section \ref{section:Automorphic forms} we construct related automorphic forms and show a non-trivial relationship between Heegner divisors 
on $\B^9/\Gamma_{\mathrm{nh}}$. 
Here $\Gamma_{\mathrm{nh}}\defeq \U(\Lambda_{\mathrm{nh}})$ where $\Lambda_{\mathrm{nh}}$  is the Hermitian form over $\Z[\omega]$ associated to the lattice $L\defeq U\oplus U(3) \oplus E_8^{\oplus 2}$.
In this paper, we give two proofs and an interpretation in terms of vector-valued modular forms; see below and Subsection \ref{subsec:second_proof}. In particular, we prove the following theorem, where $\widetilde{\O}^+(L)$
is the stable orthogonal group of $L$ of real spinor norm 1 (or, equivalently, fixing a chosen connected component of the homogenous domain associated to $L$). We refer to Subsection  \ref{AutoForm} for more details.     
\begin{thm}[{Theorem \ref{AutoMain}, Corollary \ref{AutoCor}}]
\label{mainthm:modular_form}
There exists a holomorphic automorphic form $\Psi_{\mathrm{B}}$ on $\B^9$ of weight $51$ and with character $\det$ (with respect to 
$\widetilde{\Gamma}_{\mathrm{nh}}\defeq \Gamma_{\mathrm{nh}}\cap \widetilde{\O}^+(L)$) whose zero divisor is given by
\[3(\H_{\mathrm{n}} + 28\H_{\mathrm{h}} + 3\H_{\mathrm{vt}}).\]
Moreover, we can write $\Psi_{\mathrm{B}}=\Psi\vert_{\B^9}$ for an automorphic form $\Psi$ on $\mathbb{D}_L$ (type IV domain) where $\Psi$ is the Borcherds lift of a weakly holomorphic modular form.
\end{thm}

We remark that the claim $\Psi_{\mathrm{B}}=\Psi\vert_{\B^9}$ for the Borcherds lift $\Psi$ of a weakly holomorphic form on $\B^9$, is different from constructing the form $\Psi_{\mathrm{B}}$ as a quasi-pullback. Indeed, this particular claim follows from 
Bruinier's converse theorem \cite{Bru02, Bru14, Ma19}. 
Moreover, we compute the vector-valued modular form whose Borcherds lift is $\Psi$ explicitly in Subsection \ref{subsec:Ma's result} in terms of theta functions. 

The first proof uses quasi-pullbacks of automorphic forms similar to \cite{CMJL12}.
Though they take the quasi-pullback of $\Phi_{12}$ to $\B^9$ directly, based on the communication with the second named author \cite[Acknowledgements]{CMJL12}, in this paper, we pass through $\D_L$ (a type IV domain).
The method of the second proof is due to \cite{AF02, Fre03}.
This approach has the advantage that it can be applied regardless of an embedding $L\hookrightarrow II_{2,26}$.
In this sense, our second proof gives a new intrinsic proof of \cite[Theorem 5.11]{CMJL12}.
Moreover, not taking quasi-pullbacks, we show that the automorphic form is obtained from the Borcherds product.

The Heegner divisors $\H_{\mathrm{n}}, \H_{\mathrm{h}},  \H_{\mathrm{vt}}$ correspond, in terms of moduli of genus 4 curves, to the nodal locus, to hyperelliptic genus $4$ curves together with a pair of points conjugate under the hyperelliptic involution, and the locus with a vanishing theta constant respectively. For details see (\ref{def:divisors_kondo}) and the discussion therein, in particular Proposition \ref{prop:moduli_of_nonhyperelliptic_genus_4}.
Note that the moduli space of non-hyperelliptic curves $\M_4^{\mathrm{nh}}(g_3^1)$ of genus 4 endowed with a $g_3^1$ is an open subset of $\B^9/\widetilde{\Gamma}_{\mathrm{nh}}$. 
We give two rational maps $\pi_{\mathrm{arith}}, \pi_{\mathrm{geom}}$ 
from $\M_4^{\mathrm{nh}}(g_3^1)$ to $\B^9/\Gamma$, one of them defined lattice-theoretically and
the other geometrically (Remarks \ref{rem:DM_relation}, \ref{rem:prob_modular_forms}).  It would be interesting to investigate their relationship (Question \ref{que:comparisonmaps}).
For a more detailed discussion of the relationship between the above moduli spaces from the point of view of ball quotients, we refer to Remark \ref{rem:DM_relation}.
Finally, we would like to mention a recent paper by Looijenga \cite{Loo23} in which he gave a different ball quotient structure on $\M_4^{\mathrm{nh}}(g_3^1)$ without resorting to the theory of periods of $K3$ surfaces.

\subsection{Relation to previous work}
\label{subsec: previous work}
Obviously, some of the results and the techniques used in this paper are related to other cases, notably cubic surfaces \cite{CMGHL23a}, and 8 points which we treated in \cite{HM25}.
There are, however, also a number of differences and new aspects which have not been treated before. 
First of all, we would like to point out that the situation is more complicated here, compared to the case of 8 points. The reason is that the latter case  satisfies the INT condition, while in our case it corresponds to the 
$\Sigma$INT condition in the notion of the Deligne-Mostow theory \cite{DM86, Mos86, KLW87}. 
In the case of 8 points the relationship between the unordered and the ordered cases is very straightforward. Both are ball quotients and the unordered ball quotient model is an $S_8$ quotient of the ordered
model. This extends to the Baily-Borel and the toroidal compactifications. In the case of 12 points, we do not know whether the moduli space of ordered points has a ball quotient model with respect to a suitable period map.
In fact, we conjecture that this is not the case (Remark \ref{rem:several_compactifications2} and Question \ref{que:ball_quot_model}). As a consequence, some of the arguments used in the 8 points case are no longer available here. We replace them by using the results of Keel and MacKernan on the birational 
geometry of moduli spaces of weighted sets of points on $\P^1$.      

A further new aspect of this paper is that we also study aspects of derived algebraic geometry. This was not treated in either \cite{CMGHL23a} or in \cite{HM25}. In fact, to the best of our knowledge, very little systematic 
research has been done on the topic of derived algebraic geometry of ball quotients or, in fact more generally, modular varieties.
The current paper further contains a detailed discussion of the construction of some automorphic forms that arise as Borcherds lifts. 

Finally, a very good reason to study the case of 12 points in some detail is that it is the second ancestral
Deligne-Mostow variety. Having thus completed both ancestral cases, for the Eisenstein as well as for the Gaussian integers, our results allow us to study the other Deligne-Mostow varieties more systematically.
We are planning to return to this in a future paper.

\subsection{Basic concepts used in the present paper }
\label{subsec:concepts}
At this point, we would like to comment on some important concepts that play a role in this paper and the proof of our theorems. First of all, we refer the reader to Section  \ref{section:preliminaries} for a more precise 
description of the ball quotient $\B^9/\Gamma$ and the GIT quotient $\M^{\GIT}$. 
The Kirwan blow-up was introduced in the seminal work of Kirwan \cite{Kir85}. 
For a discussion of the Kirwan blow-up in the related cases of cubic surfaces, cubic threefolds and of $8$ points on $\P^1$, we refer to \cite{CMGHL23a, CMGHL23b, HM25}. 

Second, in view of Theorems \ref{mainthm:semi-toroidal}, \ref{mainthm:not_log_crepant}, \ref{mainthm:derived_equivalent}, \ref{mainthm:derived_cubic} and Conjecture \ref{mainconj:generalized DK conjecture}, we would like to recall the 
fundamental notions of $K$-equivalence and $D$-equivalence.
As we mention in Subsection \ref{subsec:Applications of the main result}, these concepts are expected to coincide under suitable conditions through the Serre functor (Remark \ref{rem:Serre_functor}), but this has not yet been fully understood.
In this paper, a pair $(X,A)$ shall mean that $X$ is a normal projective $\Q$-Gorenstein variety, $A$ is a $\Q$-divisor on $X$ and $K_X+A$ is $\Q$-Cartier.
We recall that two normal projective $\Q$-Gorenstein varieties $X$ and $Y$ (resp. pairs $(X,B)$ and $(Y,C)$) are called \textit{$K$-equivalent} (resp. \textit{log $K$-equivalent}) 
if there is a smooth projective (resp. normal projective) variety $Z$ dominating $X$ and $Y$ birationally  
\[\xymatrix{
  & Z \ar[ld]_{f_X} \ar[rd]^{f_Y} & \\
  X \ar@{<-->}[rr] && Y
}\]
such that $f_X^*K_X \sim_{\Q} f_Y^*K_Y$ (resp. $f_X^*(K_X+B) \sim_{\Q} f_Y^*(K_Y+C)$). 
If the latter condition is satisfied for $Z=X$ and a fixed birational morphism $f_Y:X\to Y$, we use the terminology  $f_Y$ \textit{log-crepant}.
We note that for $K$-equivalent varieties, the top intersection numbers of the canonical divisors are equal: $K_X^n=K_Y^n$, where $n$ is the dimension of $X$ and $Y$.

For a smooth projective variety $X$, we can define the associated derived category $D(X)\defeq D^{\mathrm{b}}(\mathrm{Coh}(X))$ of bounded complexes of coherent sheaves on $X$.
Of course, we can also consider $D(X)$ for a singular variety $X$, but this cannot be expected to contain much information about $K_X$ or satisfy the original $DK$-conjecture; see Remark \ref{rem:Serre_functor}.
This leads us to the notion of derived categories of associated stacks.
For a projective variety $X$ with only finite quotient singularities, one can define an associated stack $\mathcal{X}(X)$, taking local covers, and its derived category $D(\mathcal{X}(X))\defeq D^{\mathrm{b}}(\mathrm{Coh}(\mathcal{X}(X)))$ following Kawamata \cite{Kaw04}.
We call two smooth projective varieties $X$ and $Y$ (resp. projective varieties $X$ and $Y$ with only finite quotient singularities) \textit{$D$-equivalent} (resp. \textit{stacky $D$-equivalent})
if there exists an equivalence $D(X)\cong D(Y)$ (resp. $D(\mathcal{X}(X))\cong D(\mathcal{X}(Y))$) as triangulated categories.
Note that if $X$ and $Y$ are smooth, then the two notions coincide.
The categories $D(X)$ and $D(\mathcal{X}(X))$ are known to encode important information on the canonical divisor $K_X$.

\begin{rem}
    \label{rem:Serre_functor}
    We briefly recall why we work with smooth stacks, rather than the original singular variety:
    \begin{enumerate}
        \item Let $X$ be a smooth projective variety.
    Then it is known that $D(X)$ has the (unique) Serre functor $\mathcal{S}_X(-)\defeq -\otimes\OO(K_X)[\dim X]$.
    This functor encodes information on the canonical divisor $K_X$ and is useful for studying derived categories.
    For example, a fully faithful functor between derived categories, commuting with $\mathcal{S}_X$, gives a categorical equivalence of triangulated categories.
    Combining this with Orlov's representability theorem \cite[Theorem 2.2]{Orl97}, one can obtain that, if $K_X$ or $-K_X$ is big, $D$-equivalence implies $K$-equivalence \cite[Theorem 2.3]{Kaw02}.
    However, if $X$ is not smooth, then $D(X)$ does not necessarily admit a Serre functor, 
    which is one of the reasons why one considers the associated smooth stack $\mathcal{X}(X)$.
    For a projective variety $X$ with only finite quotient singularities, there exists the Serre functor $S_{\mathcal{X}(X)}(-)\defeq -\otimes\omega_{\mathcal{X}(X)}[\dim X]$  on the derived category $D(\mathcal{X}(X))$, 
    which enables us to analyze the connection between canonical bundles and derived categories. 
    Analogous to the case of smooth projective varieties, it is known that for projective varieties with only finite quotient singularities, having big (anti-)canonical bundles, stacky $D$-equivalence implies $K$-equivalence, a fact which we will use  in the proof of 
    Theorem \ref{thm:derived_equivalent}.
    \item The definition of the associated stacks originally arose from the study of the Francia flop \cite[Example 4.3]{Kaw18}.
    By generalizing this, the log version of Conjecture \ref{mainconj:generalized DK conjecture} (from $K$ to $D$) holds for toroidal varieties in relation to the derived McKay correspondence; see \cite[Theorem 1.1]{Kaw16}  \cite[Sections 8, 9]{Kaw18}.
    \end{enumerate}
\end{rem}

We shall also use semi-toroidal compactifications, a concept first introduced by Looijenga in a series of papers \cite{Loo85, Loo86, Loo03a, Loo03b}. His main motivation was to define a more flexible notion of 
compactifications of locally symmetric varieties that appear naturally in the context of certain moduli problems. 
In recent years Alexeev, Engel and Odaka, in particular, have renewed the study of these compactifications from a different modular point of view.
Roughly speaking, semi-toroidal compactifications use semi-fans rather than fans (in the case of the type IV domains). They lie 
between the Baily-Borel and toroidal compactifications (both of which are special cases of semi-toroidal compactifications), and can in fact be characterized by this property; see also Subsection \ref{thm:not_semi_toric}.
Semi-toroidal compactifications are closely related to hyperplane arrangements and can be given as the projective spectrum of a suitable ring of automorphic forms with poles along the hyperplane arrangements.

In the light of these discussions, we shall propose a number of questions (Questions \ref{que:ball_quot_model}, \ref{que:quotient_singularities},  \ref{que:derived_categories}, \ref{que:comparisonmaps}), 
concerning such diverse areas as automorphic forms, singularities, derived algebraic geometry and moduli theory.

\subsection{Organization of the paper}
We will now outline the structure of the paper. We first recall the ball-quotient model and the Deligne-Mostow isomorphism in Section \ref{section:preliminaries}. 
There we give a proof of Theorem \ref{mainthm:semi-toroidal} (2) by using local computation via the Luna slice theorem.
In Section \ref{sec:HassettKeelappl} we recapitulate part of the Hassett-Keel program and use this to compute the canonical bundle of the Kirwan blow-up of the moduli space of ordered 12 points.
We also study the geometry of the exceptional divisors of the ball-quotient model in the ordered case.  Section \ref{sec:discnonK} contains one of the main results of the paper, namely Theorem  \ref{mainthm:semi-toroidal} (3) which states that the Kirwan blow-up $\M^{K}$ and the toroidal compactification $\overline{\B^9/\Gamma}^{\tor}$ are not $K$-equivalent and hence not isomorphic as abstract varieties. 
We then discuss 
connections with the LMMP and semi-toroidal compactifications in Subsection \ref{subsection:semi-toroidal}, where we show in Theorem \ref{mainthm:semi-toroidal} (1) that 
$\M^{\K}$ is not a semi-toroidal compactification.
As applications, in Section \ref{section:two_applications}, we give another proof of Theorem \ref{mainthm:semi-toroidal} (2) and prove that these spaces are not stacky $D$-equivalent (Theorem \ref{mainthm:derived_equivalent}).
The cohomology of the varieties concerned is computed in Section \ref{sec:cohomology}. 
In the rest of the paper, we consider automorphic forms related to our varieties.
Section \ref{section:Automorphic forms} is devoted to a discussion of the moduli space of non-hyperelliptic curves of genus 4 and the construction of an associated automorphic form.
Finally, we will review the automorphic form by Allcock, observe that it can be obtained by quasi-pullback and compare it with our automorphic form.

\subsection{Notation and Diagram}
Throughout the paper, we will work over the complex numbers $\C$.
To improve the readability of the paper and to fix the notation we provide a list of the moduli spaces which play a role in this paper:
\begin{itemize}
    \item The space $\M_{\ord}^{\GIT}\defeq(\P^1)^{12}/\!/\SL_2(\C)$ is the moduli space of ordered 12 points on $\P^1$. 
    \item The spaces  $\M^{\GIT}\defeq\M^{\GIT}/S_{12}\cong\P^{12}/\!/\SL_2(\C)\cong\overline{\B^9/\Gamma}^{\BB}$ are the moduli space of unordered 12 points on $\P^1$ and the corresponding realization as a ball quotient. 
    \item The space $\overline{\B^9/\Gamma}^{\tor}$ is the unique toroidal compactification of $\B^9/\Gamma$.
    \item The space $\M^{\K}_{\ord}$ (resp. $\M^{\K}$) is the Kirwan blow-up at $\frac{1}{2}\binom{12}{6}$ (resp. 1) polystable point(s) introduced in \cite{Kir85}.
    \item The compactified space $\overline{\M}_{0,12}$ is the Deligne-Mumford compactification of the configuration space $\M_{0,12}$.
    \item The compactified space ${\overline{\M}_{0,12(w)}}$ is the Hassett space of 12 points on $\P^1$ with symmetric weights $(w, \ldots , w)\in\Q^{12}$. Note that (\ref{morphism:blow-up_sequence_ordered_npoints}) is now written as 
    \begin{align*}
    \overline{\M}_{0,12} \stackrel{\varphi_4}{\to} \overline{\M}_{0,12\left(\frac{1}{4}+\epsilon\right)}\stackrel{\varphi_3}{\to}\overline{\M}_{0,12\left(\frac{1}{5}+\epsilon\right)} \stackrel{\varphi_2}{\to}\overline{\M}_{0,12\left(\frac{1}{6}+\epsilon\right)} \stackrel{\varphi_1}{\to} \M^{\GIT}_{\ord}.
\end{align*}
Each morphism is a blow-up whose center is the locus where certain points coincide and the exceptional divisors are vital divisors; see Section  \ref{sec:HassettKeelappl}.
\end{itemize}
The relationship between the spaces appearing in this paper is described in Figure \ref{fig:compactifications}.
Before discussing the properties of the morphisms and varieties appearing in this diagram we will explain the arrows appearing in this figure.
The morphisms $\varphi_i$ are contractions defined in the blow-up sequence (\ref{morphism:blow-up_sequence_ordered_npoints}).
The morphism $\overline{\Phi}_{12\left(\frac{1}{6}\right)}$ was introduced in \cite[Theorem 1.1]{GKS21}, compactifying the Deligne-Mostow period map $\Phi_{12\left(\frac{1}{6}\right)}:\M_{0,12}/S_{12}\hookrightarrow\B^9/\Gamma$.
The morphisms $\varphi_1:\M_{\ord}^{\K}\to\M_{\ord}^{\GIT}$ and $f:\M^{\K}\to\M^{\GIT}$ are the Kirwan blow-ups introduced in \cite{Kir85}.
The Kirwan blow-up $f$ and the toroidal blow-up $\pi$ are blow-ups whose centers are the same point (the cusp).
The birational map $g$ is the lift of the Deligne-Mostow isomorphism $\phi$ and the birational map which we are particularly interested in.

\begin{figure}[h]
\centering
\label{fig:compactifications}
\[
  \begin{tikzpicture}[
    labelsize/.style={font=\scriptsize},
    isolabelsize/.style={font=\normalsize},
    hom/.style={->,auto,labelsize},
    scale=1,
  ]
  
  \node(K) at (0,4) {$\M_{\ord}^{\K}$};
  \node(Qu) at (0,6) {$\overline{\M}_{0,12\left(\frac{1}{6}+\epsilon\right)}$};
  \node(S) at (0,8) {$\overline{\M}_{0,12}$};
  \node(rtor) at (12,4) {$\overline{\B^9/\Gamma}^{\tor}$}; 
  \node(git) at (0,0) {$\M^{\GIT}_{\ord}$};
  \node(cbb) at (12,0) {$\overline{\B^9/\Gamma}^{\BB}$};
  \node(gitl) at (6,0) {$\M^{\GIT}$}; 
  \node(k) at (6,4) {$\M^{\K}$};
  \node(s) at (12,8) {$\M_{\ord}^{\K}/S_{12}=\widetilde{\M}_{0,12\left(\frac{1}{6}+\epsilon\right)}$};
  
  \draw[hom,swap] (Qu) to
    node[align=center,yshift=-2pt]{}
    node[swap]{\normalsize $\rotatesim$}
  (K);

  \draw[hom] (gitl) to
    node[align=center]{\scriptsize $\phi$}
    node[swap]{\normalsize$\sim$}
  (cbb);

  \draw[hom] (s) to
    node[align=center]{\scriptsize $\overline{\Phi}_{12\left(\frac{1}{6}\right)}$}
    node[swap,xshift=2pt]{\normalsize $\rotatesim$}
  (rtor);
  

  \draw[hom] (git) to node{$/S_{12}$} (gitl);
  \draw[hom,swap] (K) to node{$\varphi_1$} (git);

  \draw[hom,swap] (k) to node{$f$} (gitl);
  \draw[hom,dashed] (k) to node{$g$} (rtor);

  \draw[hom] (rtor) to node{$\pi$} (cbb);
  \draw[hom](K) to node{$/S_{12}$}(s);
  \draw[hom](S) to node[swap]{$\varphi_2\circ\varphi_3\circ\varphi_4$}(Qu);
  
  \end{tikzpicture}
\]
\caption{Relationship between the compactifications of moduli spaces}
\label{com_diag}
\end{figure}
For the arrows and varieties in the figure, the following properties hold:
\begin{enumerate}
    \item The morphisms $\phi$ and $\overline{\Phi}_{12\left(\frac{1}{6}\right)}$ are  isomorphisms given by \cite{Mos86} and \cite[Theo\-rem 1.1]{GKS21}.
    \item The varieties $\overline{\M}_{0,12\left(\frac{1}{6}+\epsilon\right)}$ and $\M^{\K}_{\ord}$ are isomorphic by the construction of the blow-up sequence (\ref{morphism:blow-up_sequence_ordered_npoints}); see \cite[Theorem 1.1]{KM11} and Theorem \ref{previous_thm:kirwan_ordered}.
    \item The Kirwan blow-up $\M^{\K}_{\ord}$ is nonsingular \cite[Section 4]{KM11}.
\end{enumerate}
For the discussion in Subsection \ref{subsec:auxiliary_calculations} and Section \ref{section:semi-toroidal_compactifications}, we introduce the notation of discriminant divisors and vital divisors.

\begin{itemize}
    \item The divisor $\H$ is defined in (\ref{eq:discriminant_divisor_on_ball_quotient}). It is a Heegner divisor in $\B^9$, whose $\Gamma$-quotient is denoted by $H$. We write its closure in the Baily-Borel and toroidal compactifications by $\overline{H}^{\BB}$ and $\overline{H}^{\tor}$ respectively.
    \item The divisor $\D$
    is the discriminant divisor in $\M^{\GIT}$ and is isomorphic to $\overline{H}^{\BB}$ through the Deligne-Mostow isomorphism $\phi$.
    \item We denote the strict transform of $\D$ with respect to the Kirwan blow-up $f$ by $\widetilde{\D}$.
    \item In the case of ordered points, we use a similar notation as above, adding the subscript ``ord": $\D_{ord}$, $\widetilde{\D_{\ord}}$ and $\Delta_{\ord}$.
    \item The divisors $D_k^{(\ell)}$ lie in the space $\overline{\M}_{0,12(\epsilon+\frac{1}{7-\ell})}$, which are called vital divisors studied in \cite{Has03, KM11, KM13}.
    Roughly speaking, these are the strict transform of the locus where at least $k$ points coincide on $\M_{\ord}^{\GIT}$ through some composition of $\varphi_i$'s.
\end{itemize}

\subsection*{Acknowledgements}
The authors wish to express their thanks to Sebastian Casalaina-Martin for helpful discussions.
The first named author is grateful to Brendan Hassett for enlightening discussions about moduli spaces of pointed curves. He is also grateful to Caucher Birkar for inviting him to give an online seminar at Tsinghua University. It was during this seminar that the question about the $DK$-equivalence was raised. 
The third named author would also like to thank Shouhei Ma for his insightful suggestion on automorphic forms, Makoto Enokizono, Teppei Takamatsu for the discussion on derived categories and LMMP, and Takuya Yamauchi for letting us know the papers about linear algebraic groups.
We are also grateful to Daniel Allcock for feedback on a preliminary draft of this paper.
This research was supported through the program ``Oberwolfach Research Fellows" by the Mathematisches Forschungsinstitut Oberwolfach in 2023.  The second named author is partially supported by JSPS Grant-in-Aid for Scientific Research (A) No.20H00112.
\section{Preliminaries}
\label{section:preliminaries}
In this section, we recall some details about the Deligne-Mostow isomorphism and perform some local computations.

\subsection{The Deligne-Mostow isomorphism}
\label{subsection:DMiso}
We denote the ring of Eisenstein integers by $\Z[\omega]$, where $\omega$ is a primitive cube root of unity. 
The Hermitian forms that we will consider have values in $\frac{1}{\sqrt{-3}}\Z[\omega]$. More precisely, we consider the Hermitian lattice $\Lambda$ of signature $(1,9)$ whose underlying 
integral lattice is $U^{\oplus 2}\oplus E_8^{\oplus 2}$. 
In this paper, we follow the convention that root lattices, such as $E_8$, are negative definite lattices, a convention which we also use in Section \ref{section:Automorphic forms}.
More concretely, this is given by the Gram matrix     
\[\frac{1}{\sqrt{-3}}\begin{pmatrix}
0 & 1 \\
-1 & 0 \\
\end{pmatrix}
\oplus
\frac{-1}{\sqrt{-3}}\begin{pmatrix}
\sqrt{-3} & 0 &-1 & -1 \\
0 & \sqrt{-3} & -1 & 1 \\
1 & 1 & \sqrt{-3} & 0 \\
1 & -1 & 0 & \sqrt{-3}
\end{pmatrix}^{\oplus 2};
\]
see also \cite[Subsection 7.1.1]{CMGHL23a}.
We denote by $(-,-)$ the Hermitian form.
Let $\Gamma\defeq\U(\Lambda)$ be the integral unitary group.
This group acts on the 9-dimensional ball
\[\B^9\defeq\{v\in\P(\Lambda\otimes_{\Z[\omega]}\C)\mid (v,v)>0\},\]
which is a type I Hermitian symmetric domain.
By the work of Allcock \cite[Theorem 5.1]{All00} 
\begin{align}
\label{eq:DM}
    (\M_{\ord}^{\GIT}/S_{12}\cong) \M^{\GIT}\cong \overline{\B^9/\Gamma}^{\BB}.
\end{align}
Under the isomorphism  (\ref{eq:DM}) the unique polystable point $c_{6,6}$, corresponding to two different points with multiplicity 6 each, is sent to the unique Baily-Borel cusp.
Now, let 
\begin{align}
\label{eq:discriminant_divisor_on_ball_quotient}
    \H\defeq\bigcup_r\{[v]\in\B^9\mid (v,r) = 0\}
\end{align}
where $r$ runs $(-1)$-vectors in $\Lambda$ (it should be borne in mind that we use the $\frac{1}{\sqrt{-3}}\Z[\omega]$-valued Hermitian form here). 
Then, the discriminant divisor in $\M^{\GIT}$, parametrizing stable but not distinct $12$-tuples, corresponds to  $H\defeq\H/\Gamma$, see  \cite[Subsection 7.3]{Loo03a}, and this is irreducible.
In other words, the isomorphism (\ref{eq:DM}) extends the isomorphism between open varieties:
\[\M_{0,12}/S_{12}\cong (\B^9\setminus \H)/\Gamma.\]
By Deligne-Mostow theory this is essentially the period map for the abelian varieties constructed from the quotient of the Jacobian of a $\mu_6$-cover of $\P^1$ branched along unordered 12 points.

\begin{rem}
\label{rem:several_compactifications1}
For future use, we remark that there is a connection with semi-toroidal compactifications. By  \cite[Definition 4.2]{Loo03a} the Baily-Borel compactification of any ball quotient or locally symmetric space of type IV
is a semi-toroidal compactification, namely the one obtained by taking the empty hyperplane arrangement, see \cite[Definition 4.2]{Loo03a}. Hence, in our case, starting with the isomorphism 
    \[(\P^{12})^{\mathrm{s}}/\!/\SL_2(\C) \cong \B^9/\Gamma\]
and applying \cite[Theorem 7.1]{Loo03a} to the empty hyperplane arrangement, we recover (\ref{eq:DM}), corresponding to the Baily-Borel system \cite[Definition 4.2]{Loo03a}.
In this framework, moreover, it is known that the Deligne-Mumford compactification $\overline{M}_{0,12}$ is known to be the minimal normal crossing blow-up of $\M^{\GIT}$ with respect to the hyperplane arrangement $\D$ \cite[Subsection 7.3]{Loo03a} (see also \cite{KM11}).
In other words, $\D$ can be stratified by loci describing how points collide, corresponding to lower dimensional ball quotients in $\overline{\B^9/\Gamma}^{\BB}$, and $\overline{M}_{0,12}$ is the minimal compactification in the sense that all strata become normal crossing divisors.
Note that the toroidal compactification $\overline{\B^9/\Gamma}^{\tor}$ is the first step of this sequence of blow-ups, see \cite[Theorem 1.1]{GKS21}.
\end{rem}

\begin{rem}
\label{rem:several_compactifications2}
It should be noted that the case of 12  \textit{ordered} points behaves differently from the other ancestral case, namely 8 points. In the latter case, both the moduli of unordered and ordered points are 
Deligne-Mostow varieties. Indeed, these two moduli spaces have  5-dimensional ball quotient models $\B^5/\Gamma_{\ord}$ and $\B^5/\Gamma$ respectively.
The arithmetic groups $\Gamma$  and $\Gamma_{\ord}$ are the unitary group $\U(U\oplus U(2)\oplus D_4^{\oplus 2})$ 
and its {\em stable}  subgroup. The latter is defined as the subgroup acting trivially on the discriminant
$(U \oplus U(2)\oplus D_4^{\oplus 2})^*/U\oplus U(2)\oplus D_4^{\oplus 2}\cong (\Z/2\Z)^6$.
Here and below, we denote by $M^*$ the dual lattice of a given lattice $M$.
Note that 
$\Gamma/\Gamma_{\ord}\cong S_8$ acts on the set of non-isotropic vectors $(\Z/2\Z)^3\subset(\Z/2\Z)^6$ transitively, and that 
$S_8\cong \O^+(6, \mathbb{F}_2)$ is the even orthogonal group over $\mathbb{F}_2$ in dimension 8, see also \cite[p. 22]{Con85}. The action of the group $S_8$ on $\B^5/\Gamma_{\ord}$ corresponds to the 
permutation on the set of 8 points.

The situation is also similar for 5 and 6 points where both the moduli spaces of ordered and unordered points are ball quotients with the covering map given by a suitable quotient of arithmetic groups.
For 6 points this  is $\O^-(4,\mathbb{F}_3)\cong \Z/2\Z\times S_6$ (cf. \cite[p. 4]{Con85}, \cite{Kon13}) 
where $\O^-(4,\mathbb{F}_3)$ is the odd orthogonal group. Finally,  for 5 points we have  $\O(3,\mathbb{F}_5)\cong \Z/2\Z\times S_5$ (cf. \cite[p. 2]{Con85}, \cite{Kon07}).
In both cases, the factor $\Z/2\Z$ goes to the identity under the projection to the projective unitary group and hence acts trivially on the complex ball. 
Note that in all of these cases ($n=5,6,8)$ the condition INT of \cite{DM86} is satisfied, which implies that the moduli space of ordered points is also a ball quotient. 

However, in the case of 12 points, the INT condition is not satisfied, hence this moduli problem does not appear in \cite[(14.4)]{DM86}, but it does appear in \cite[Section 5]{Mos86} (see \cite[Appendix]{Thu98} for a complete list of INT and $\Sigma$INT). 
This fits with the fact that the symmetric group $S_{12}$ does not have an exceptional isomorphism as above; see \cite[p. 91]{Con85}. It is, however, conceivable that the moduli space of 12 ordered 
points still has a ball quotient model, e.g. via a different period map. We do not know whether this is the case. 
We state this as a formal problem below.

\begin{que}
\label{que:ball_quot_model}
    Does the moduli space $\M_{\ord}^{\GIT}$ of ordered 12 points on $\P^1$ have a ball-quotient model?
    In other words, is there an arithmetic subgroup $\Gamma_{\ord}$ so that 
    \[\M_{\ord}^{\GIT}\cong \overline{\B^9/\Gamma_{\ord}}^{\BB}?\]
    A stronger question is to ask, whether such a $\Gamma_{\ord}$ exists, with the additional property that  $\Gamma/\Gamma_{\ord}\cong S_{12}$ and the above isomorphism is compatible with the natural 
    action of $S_{12}$ on     
    \[\M_{0,12}\cong (\B^9\setminus \H)/\Gamma_{\ord}.\]
    Note that the two ball models may come from a different period map.
\end{que}
\end{rem}

\subsection{Computations involving the Luna slice}
\label{subsec:auxiliary_calculations}
In this subsection, we collect the result regarding the local computations and conclude Theorem \ref{mainthm:semi-toroidal} (2).
We will keep this discussion short since many of the arguments are similar to analogous considerations in \cite[Section 3]{HM25}.    
We denote by $x_0, x_1$ the homogeneous coordinates on $\P^1$.
Let $p_{6,6}$ be the point corresponding to $x_0^6x_1^6$, which is mapped to the unique polystable orbit $c_{6,6}$ in $\M^{\GIT}$.
First, we recall the local description in terms of the Luna slice theorem.

\begin{lem}
\label{lem:Luna_slice}
\begin{enumerate}
    \item The stabilizer $R$ of $p_{6,6}$ and its connected component of the identity are described as
\begin{align*} 
    R&\defeq\Stab(p_{6,6})=\Bigl\{\begin{pmatrix}
\lambda & 0 \\
0 & \lambda^{-1} \\
\end{pmatrix}\in\SL_2(\C)\Bigr\}\bigcup \Bigl\{\begin{pmatrix}
0 & \lambda \\
-\lambda^{-1} & 0 \\
\end{pmatrix}\in\SL_2(\C)\Bigr\}\cong \C^{\times}\rtimes S_2,\\
    R^{\circ}&\defeq\Stab(p_{6,6})^{\circ}\cong\C^{\times}.
    \end{align*}  
    \item A Luna slice for $p_{6,6}$, normal to the orbit $\SL_2(\C)\cdot \{p_{6,6}\}\subset\P^{12}$, 
is isomorphic to $\C^{10}$, spanned by the $10$ monomials 
\[x_0^{12}, \quad x_1^{12},\quad x_0^{11}x_1,\quad x_0x_1^{11},\quad x_0^{10}x_1^2,\quad x_0^2x_1^{10},\quad x_0^9x_1^3,\quad x_0^3x_1^9,\quad x_0^8x_1^4, \quad x_0^4x_1^8\]
in the tangent space $H^0(\P^1, \OO_{\P^1}(12))$.
Projectively, 
\begin{align*}
    \P^{10}&=\{\a_0x_0^{12}+\a_1x_1^{12}+\b_0x_0^{11}x_1+\b_1 x_0x_1^{11}+\g_0x_0^{10}x_1^2+\g_1x_0^2x_1^{10} \\
    &+\d_0 x_0^9x_1^3+\d_1 x_0^3x_1^9 + \e_0 x_0^8x_1^4+ \e_1 x_0^4x_1^8+kx_0^4x_1^4\}\\
    &\subset \P H^0(\P^1,\OO_{\P^1}(12)) =\P^{12}.
\end{align*}
\end{enumerate}
  
\end{lem}
\begin{proof}
Since these follow from straightforward calculations which are similar to those in \cite[Subsection 4.3.1]{CMGHL23a}, \cite[Lemma 3.3, 3.4]{CMGHL23b} and \cite[Lemma 3.1, 3.2]{HM25}, we omit the details.
\end{proof}

The following Theorem turns out to be a crucial observation. 

\begin{thm}
\label{thm:nonord_nontransversal}
The strict transform $\widetilde{\D}$ of the discriminant divisor $\D\subset\M^{\GIT}$ and the exceptional divisor $\Delta$ of the Kirwan blow-up $\M^{\K}\to\M^{\GIT}$
do not meet generically transversally in $\M^{\K}$.
\end{thm}
\begin{proof}
This is a local calculation using the Luna slice from Lemma \ref{lem:Luna_slice}. The structure of the proof is the same as in the proof of \cite[Theorem 3.4]{HM25}. Recall that the exceptional Kirwan
divisor lies over the point $c_{6,6}$  which is the orbit of the point $p_{6,6}$ given by $x_0^6x_1^6$. This brings the versal deformation of $x^6$ into play, which is given by the sextic 
polynomial $x^6+\a_0 x^4+ \b_0 x^3 + \g_0 x^2 + \d_0 x +\e_0$.  
A computer-based calculation shows that the discriminant $d(\a_0, \b_0, \g_0, \d_0, \e_0)$  is of the form 
\begin{align}
    \label{eq:discriminant}
    - 46656 \e_0^5 + (\mathrm{terms\ of\ degree}\geq 6).
\end{align}
The rest of the argument is then a lengthy, but straightforward, calculation as in the proof of \cite[Theorem 3.4]{HM25}.  
\end{proof}

\begin{rem}\label{rem:nontransversal}
We find it noteworthy that this situation can be observed for moduli of cubic surfaces, 8 points on $\P^1$ and again for 12 points on $\P^1$, but that the situation changes each time one 
goes to the natural level covers of marked cubic cases and ordered sets of points respectively.     
\end{rem}

Next, we give a first proof of Theorem \ref{mainthm:semi-toroidal} (2) in the spirit of  \cite[Section 4.3]{CMGHL23b}.

\begin{thm}
\label{thm:lift}
Neither the Deligne-Mostow isomorphism $\phi$ nor its inverse $\phi^{-1}$ lift to a morphism between the Kirwan blow-up $\M^{\K}$ and the toroidal compactification $\overline{\B^9/\Gamma}^{\tor}$.
\end{thm}
\begin{proof}
We first note that the birational lift $g$ of $\phi$ is not an isomorphism. This follows, since by Theorem \ref{thm:nonord_nontransversal} the Kirwan exceptional divisor $\Delta$ and the discriminant $\D$ do not intersect 
transversally in $\M^{\K}$, whereas the closure $\overline{H}^{\tor}$ of the discriminant  $H$ and the toric boundary $T$  intersect generically transversally in $\overline{\B^9/\Gamma}^{\tor}$ by Remark \ref{rem:transversally}. Now assume that 
$g: \M^{\K} \to \overline{\B^9/\Gamma}^{\tor}$ is a morphism. Since we know that the restriction of $g$ induces an isomorphism 
$g|_{\M^{\K} \setminus \Delta} : \M^{\K} \setminus \Delta \cong \overline{\B^9/\Gamma}^{\tor} \setminus H$, it follows that the exceptional Kirwan divisor $\Delta$ must be mapped to the toric
boundary $T$. Since $\M^{\K}$ and $\overline{\B^9/\Gamma}^{\tor}$ are normal varieties, $g$ cannot be a bijection and must hence be a small contraction.
However, this contradicts the fact 
that $\overline{\B^9/\Gamma}^{\tor}$ is $\Q$-factorial. The same argument also applies to $g^{-1}$.   
\end{proof}
\begin{rem}
    \label{rem:another_proof_semi_toroidal}
    \begin{enumerate}
        \item   We remark that in our previous paper \cite[Theorem 1.1, Remark 3.6]{HM25}, we used the coincidence of Betti numbers of two spaces to prove the nonexistence of a morphism.
    We note that this point can be avoided by discussing the geometry as in the above proof.
    \item     Using recent developments in $K$-stability and the LMMP, another proof  for the part of Theorem \ref{thm:lift}, saying that $\phi^{-1}$ does not lift to a morphism can be derived from 
    Theorem \ref{thm:not_semi_toric} and \cite{AE23, AEH24}; see \cite[Remark 4.8]{HM25}.
    \end{enumerate}

\end{rem}

\section{The Hassett-Keel program and an application}
\label{sec:HassettKeelappl}
The Hassett-Keel program aims at giving a modular interpretation of the configuration spaces and their log canonical models.
We want to give two applications for this program: the geometry of the toroidal compactification and the description of the canonical bundles.
For the precise definition of the spaces $\overline{\M}_{0,\A}$ and the vital divisors $D_{\ell}^{(k)}$, see \cite[Sections 2, 7]{Has03} and \cite{KM11}.
In this context, \cite[Theorem 1.1]{KM11} gives the following blow-up sequence:
\begin{equation}
\label{equ:blowupsequenceHK}
\overline{\M}_{0,12}\cong\overline{\M}_{0,12\left(\frac{1}{3}+\epsilon\right)}\xrightarrow{\varphi_4} \overline{\M}_{0,12\left(\frac{1}{4}+\epsilon\right)}\xrightarrow{\varphi_3}\overline{\M}_{0,12\left(\frac{1}{5}+\epsilon\right)}\xrightarrow{\varphi_2}\overline{\M}_{0,12\left(\frac{1}{6}+\epsilon\right)}\cong\M^{\K}_{\ord}\xrightarrow{\varphi_1}\M_{\ord}^{\GIT}.
\end{equation}
Here the right-hand isomorphism is (\ref{equ:KirwanequalHassett}) in Theorem \ref{previous_thm:kirwan_ordered}.
As mentioned in the introduction, the Deligne-Mumford compactification $\overline{\M}_{0,12}$ is a normal crossing compactification of $(\P^{12})^{\mathrm{s}}/\!/\SL_2(\C)$. 
Note that the number of the strictly semi-stable orbits of $\M^{\GIT}_{\ord}$  is $\frac{1}{2}\binom{12}{6}$. These points correspond to the cusps (which are singular points) of $\M^{\GIT}_{\ord}$ and are blown up in the Kirwan blow-up.

By \cite[Proposition 5.4, Lemma 5.5]{KM11} the canonical bundles of all spaces in the blow-up sequence \ref{equ:blowupsequenceHK} are known: 
\begin{align*}
    K_{\M^{\GIT}_{\ord}}&=-\frac{2}{11}D_2^{(0)}\\
K_{\M^{\K}_{\ord}}&=-\frac{2}{11}D_2^{(1)}+\frac{14}{11}D_6^{(1)}\\
K_{\overline{\M}_{0,12\left(\frac{1}{5}+\epsilon\right)}}&=-\frac{2}{11}D_2^{(2)}+\frac{13}{11}D_5^{(2)}+\frac{14}{11}D_6^{(2)}\\
K_{\overline{\M}_{0,12\left(\frac{1}{4}+\epsilon\right)}}&=-\frac{2}{11}D_2^{(3)}+\frac{10}{11}D_4^{(3)}+\frac{13}{11}D_5^{(3)}+\frac{14}{11}D_6^{(3)}\\
K_{\overline{\M}_{0,12}}&=-\frac{2}{11}D_2^{(3)}+\frac{10}{11}D_4^{(4)}+\frac{5}{11}D_4^{(4)}+\frac{13}{11}D_5^{(4)}+\frac{14}{11}D_6^{(4)}.\\
\end{align*}
Here the number 11 comes from $n-1=11$.
Note that $D_2^{(1)}/S_{12}$ (resp. $D_6^{(1)}/S_{12}$) is the unique discriminant divisor $\overline{H}^{\tor}$ (resp. the exceptional divisor $T$) of the toroidal compactification $\pi:\overline{\B^9/\Gamma}^{\tor} (\cong \M^{\K}_{\ord} )\to\overline{\B^9/\Gamma}^{\BB}$.
\begin{rem}
\label{rem:transversally}

    Comparing the stability conditions in the sense of Deligne-Mumford in \cite{AL02}, the quotient  $\overline{\M}_{0,12}/S_{12}$ is exactly the Deligne-Mumford compactification of $\M_{0,12}/S_{12}$, and hence a 
    normal crossing compactification of $\M_{0,12}/S_{12}$.
    Hence, the divisors $D_2^{(1)}/S_{12}$ and $D_6^{(1)}/S_{12}$  intersect generically transversally.
    This implies that the boundary and discriminant divisors generically intersect transversally on the toroidal compactification $\overline{\B^9/\Gamma}^{\tor}$, which is an analog to \cite[Theorem 3.14]{HM25}.
\end{rem}

Additionally, from \cite[Lemma 5.3]{KM11}, the pullback of the blow-up $\varphi_1$ can be described as 
\begin{align}
\label{eq:discrepancy_ordered_boundary}
    \varphi_1^*(D_2^{(0)}) = D_2^{(1)} + 15 D_6^{(1)}.
\end{align}
Hence, substituting this in the above, it follows that 
\begin{align}
\label{eq:discrepancy_ordered_canonical_bundles}
K_{\M^{\K}_{\ord}}=-\frac{2}{11}\varphi_1^*(D_2^{(0)})+4D_6^{(1)}
\end{align}
where we have used the identification $\M^{\K}_{\ord}\cong \overline{\M}_{0,12\left(\frac{1}{6}+\epsilon\right)}$. 
Let $\D_{\ord}$ be the union of the discriminant divisors in $\M_{\ord}^{\GIT}$ and denote their strict transforms in $\M_{\ord}^{\K}$ by $\widetilde{\D_{\ord}}$.
We further denote the exceptional divisor in $\M_{\ord}^{\K}$ by $\Delta_{\ord}$.
Then (\ref{eq:discrepancy_ordered_boundary}) and (\ref{eq:discrepancy_ordered_canonical_bundles}) can be rewritten as
\begin{align}
\label{eq:canonical_bundle_kirwan_ord}
     \varphi_1^*(\D_{\ord}) &= \widetilde{\D_{\ord}} + 15\Delta_{\ord}\notag\\
     K_{\M^{\K}_{\ord}}&=-\frac{2}{11}\varphi_1^*(\D_{\ord})+4\Delta_{\ord}. 
\end{align}

The boundary of $\M^{\K}_{\ord}$ consists of $\frac{1}{2}\binom{12}{6}$ components $\Delta_{\ord,i} \cong \P^4\times\P^4$. 
The latter isomorphism follows from the moduli interpretations of $\overline{\M}_{0,\A}$ \cite[Remark 4.6]{Has03}.
Applying a similar computation as in \cite[Proposition 4.3]{HM25}, we can describe the normal bundle of the boundary in the Kirwan blow-up.
\begin{prop}
\label{prop:normalbundleKirwan}
    The normal bundle of the irreducible components of the boundary divisors $\Delta_{\ord,i}\cong \P^4\times\P^4$ in $\M^{\K}_{\ord}$ is 
    \[\N_{\Delta_{\ord,i}/\M^{\K}_{\ord}} \cong \OO_{\P^4\times\P^4}(-1,-1).\]
\end{prop}
\begin{proof}
From (\ref{eq:canonical_bundle_kirwan_ord}), we have
     \[(K_{\M^{\K}_{\ord}}+\Delta_{\ord,i})\vert_{\Delta_{\ord,i}}=(-\frac{2}{11}\varphi_1^*(\D_{\ord})+4\Delta_{\ord}+\Delta_{\ord,i})\vert_{\Delta_{\ord,i}}.\]
Since $\M^{\K}_{\ord}$ is smooth, the left-hand side is 
\begin{align*}
(K_{\M^{\K}_{\ord}}+\Delta_{\ord,i})\vert_{\Delta_{\ord,i}}&=K_{\P^4\times\P^4}\\
&=\OO(-5,-5)
\end{align*}
by the adjunction formula.
The right-hand side is 
\[5\Delta_{\ord,i}\vert_{\Delta_{\ord,i}},\]
which proves the claim.
\end{proof}

\section{Discrepancies and the failure of (log) $K$-equivalence. }
\label{sec:discnonK}
We are now ready to prove one of the main results of our paper, namely the claim that $\M^{\K}$ and  $\overline{\B^9/\Gamma}^{\tor}$ are not $K$-equivalent. For this purpose, we first compute the following.   
\begin{lem}
\label{lem:intersetionT}
The top intersection number of the toroidal boundary $T$ of  $\overline{\B^9/\Gamma}^{\tor}$ is given by
\begin{align}
\label{eq:top_self_intersection_toroidal}
T^{9}=\frac{\binom{8}{4}\cdot\frac{1}{2}\binom{12}{6}}{12!}=\frac{7}{144\cdot 6!}.    
\end{align}
\end{lem}
\begin{proof}
We first recall that there is an isomorphism $\overline{\B^9/\Gamma}^{\tor}\cong\M_{\ord}^{\K}/S_{12}$. 
This follows from Diagram \ref{fig:compactifications}, more precisely the fact that $\overline{\Phi}_{12\left(\frac{1}{6}\right)}$ is an isomorphism, as shown in \cite{Mos86} and \cite[Theo\-rem 1.1]{GKS21}. 
The rest is a straightforward calculation using Proposition \ref{prop:normalbundleKirwan},
together with $(\OO_{\P^4\times\P^4}(-1,-1))^8=\binom{8}{4}$ and the fact that the number of the polystable points in $\M_{\ord}^{\GIT}$ is $\frac{1}{2}\binom{12}{6}$. 
\end{proof}

We shall now compute the canonical bundles of these spaces in terms of discrepancies.
As we saw in the proof of Theorem \ref{thm:not_semi_toric} we have, using the notation of \cite[Lemma 6.4]{CMGHL23b},  $c=10$ and $|G_X|=|G_F|=2$; there is no divisorial locus having a strictly bigger stabilizer than $G_X$.
Note that the GIT quotient is $\Q$-Gorenstein by Proposition \ref{prop:singularity_ball_quotient}.
Hence, on the one hand, it follows that
\begin{align}
\label{eq:canonical_bundles_Kirwan}
   K_{\M^{\K}} = f^*K_{\M^{\GIT}} + 9\Delta. 
\end{align}
On the other hand, the proof of \cite[Proposition 5.8]{CMGHL23b} leads to 
\begin{align}
\label{eq:canonical_bundle_ball_quotients_pre}
    K_{\overline{\B^9/\Gamma}^{\tor}} = \pi^*K_{\overline{\B^9/\Gamma}^{BB}} + \frac{2+a(\Delta,\M^{\GIT}_{\ord})-r(T)}{r(T)} T,
\end{align}
where $r(T)$ is the branch index of the finite quotient map $\pi:\M_{\ord}^{\K} \to \M_{\ord}^{\K}/S_{12}\cong\overline{\B^9/\Gamma}^{\tor}$ along $T$ and $a(\Delta,\M^{\GIT}_{\ord})$ is the discrepancy of $\Delta$ defined by
\[f^*(\D_{\ord})=\widetilde{\D_{\ord}}+a(\Delta,\M^{\GIT}_{\ord})\Delta_{\ord}.\]
By (\ref{eq:discrepancy_ordered_boundary}), we have $a(\Delta,\M^{\GIT}_{\ord})=15$.
In the case of 8 points, the ramification order $r(T)$ was computed in \cite{HM25} using the ball quotient model for the ordered case. This argument is not available here. 
Instead, we use the diagram in Figure \ref{figure:branch_index}.
\begin{figure}[h]
  \[\begin{CD}
   \overline{\M}_{0,12}   @>{\widetilde{\Pi}}>>\widetilde{\M}_{0,12}=\overline{\M}_{0,12}/S_{12}  \\
  @VVV    @VVV \\
     \M_{\ord}^{\K}  @>{\Pi}>>  \M_{\ord}^{\K}/S_{12}\cong\overline{\B^9/\Gamma}^{\tor}
  \end{CD}\]
  \caption{Blow-up sequence and symmetries}
\label{figure:branch_index}
\end{figure}

In Figure \ref{figure:branch_index}, the two morphisms $\Pi$ and $\widetilde{\Pi}$ are $S_{12}$-quotients and the vertical morphisms are blow-ups along the locus where points coincide, as described in Section \ref{sec:HassettKeelappl}.
By \cite[Lemma 3.4]{KM13}, the quotient map $\widetilde{\Pi}$ is unramified in codimension 1 outside the discriminant divisor $D_2^{(3)}$.
This implies that $\Pi$ is unbranched along the toroidal boundary $T$, thus we have $r(T)=1$.
Hence Equation (\ref{eq:canonical_bundle_ball_quotients_pre}) can be rewritten as 
\begin{align}
\label{eq:canonical_bundle_ball_quotients}
    K_{\overline{\B^9/\Gamma}^{\tor}} = \pi^*K_{\overline{\B^9/\Gamma}^{\BB}} + 16 T.
\end{align}

We can now prove Theorem \ref{mainthm:semi-toroidal} (3).
\begin{thm}
\label{thm:not_k_equiv}
    The varieties $\M^{\K}$ and $\overline{\B^9/\Gamma}^{\tor}$ are not $K$-equivalent.
\end{thm}
\begin{proof}
    If these two varieties were $K$-equivalent then, as we recalled in the introduction, the top intersection numbers of their canonical bundles would be the same.  By (\ref{eq:canonical_bundles_Kirwan}) and (\ref{eq:canonical_bundle_ball_quotients}) this is equivalent to
    \[(9\Delta)^9 = \left(16 T\right)^9.\]
    Combined with (\ref{eq:top_self_intersection_toroidal}), this implies
    \begin{align}
    \label{eq:self_intersection_of_Delta}
        \Delta^9 = \frac{16^9\cdot7}{9^9\cdot 144\cdot 6!}.
    \end{align}
    Our calculations in Subsection \ref{subsec:auxiliary_calculations} for the Luna slice show that the order of the stabilizer of a
     point in $\Delta$ is not divisible by 3. Hence, arguing as in \cite[Proposition 3.7]{HM25}, it follows that
     $\Delta^9\in\frac{1}{e}\Z$ with $3\nmid e$, and this is a contradiction.
\end{proof}

Using the classical results of Mumford, we shall now discuss the relation with automorphic forms.
\begin{rem}
    \label{rem:lattice_theoretic_branch}
    The uniformization map $\B^9\to\B^9/\Gamma$ is ramified in codimension 1 along the discriminant divisor $\H$ with branch index 6, see (\ref{eq:discriminant_divisor_on_ball_quotient}).
    This can be proven by a lattice theoretic computation of which we give a sketch here. 
    For a $(-1)$-vector $\ell\in \Lambda$ and unit element $\xi\in\Z[\omega]^{\times}\setminus\{1\}$, we define a unitary reflection by
    \[\sigma_{\ell,\xi}(r)\defeq r+(1-\xi)\frac{(\ell,r)}{2}r\in \Lambda\otimes\Q(\omega), \, r \in \Lambda.\]
    Note that all branch divisors of the morphism $\B^9\to \B^9/\Gamma$ arise as the fixed divisors by such a reflection; see \cite[Cororally 3]{Beh12} for the classification of unitary reflections.
    One easily checks that the element $\sigma_{\ell,\xi}$ is not contained in $\Gamma$, except for $\xi=-\omega$.
    If $\xi=-\omega$, this reflection $\sigma_{\ell,-\omega}$ is an element of order 6 in $\U(\Lambda\otimes\Q(\omega))$ and hence called a \textit{hexaflection}. Its fixed point set is a Heegner divisor $\H_0$.
    By a straightforward computation one can prove that this hexaflection is contained in $\Gamma$, more precisely $\sigma_{\ell,-\omega}\in\Gamma$. 
    The ramification divisor $\H$ is the union of all $\Gamma$-translates of $\H_0$. We denote its image in $\B^9/\Gamma$ by $H$.   
    The map $\B^9\to\B^9/\Gamma$ is branched exactly along $H$ with ramification order 6. 
\end{rem}
By the Hirzebruch proportionality principle, 
the canonical bundle of our ball quotient can be described as
\begin{align}
    K_{\overline{\B^9/\Gamma}^{\BB}} &= 10\L -\frac{5}{6}\overline{H}^{\BB}\label{eq:canonical_bundle_bb}\\
    K_{\overline{\B^9/\Gamma}^{\tor}} &= 10\L -\frac{5}{6}\overline{H}^{\tor} - T \label{eq:canonical_bundle_toroidal}
\end{align}
where $\L$ is the automorphic $\Q$-line bundle of (arithmetic) weight 1 on the Baily-Borel compactification $\overline{\B^9/\Gamma}^{\BB}$.
By abuse of notation, we use the same notation $\L$ as the $\Q$-line bundle on the Baily-Borel compactification and its pullback to the toroidal compactification.
The coefficients $10$ in the above presentation come from the fact that the compact dual of $\B^9$ is $\P^9$, whose canonical bundle is $\OO(-10)$.
Here, we need to be careful with the notion of \textit{weight}.
On the one hand, a section of $k\L$ for a $k\in\Z$ defines an automorphic form of weight $k$, which means that the function satisfies the modular symmetry with respect to the automorphic factor to the power of $k$, and this quantity $k$ is called the \textit{arithmetic weight}.
On the other hand, Mumford et al. \cite[Chapter IV]{AMRT10} studied the extendability of pluricanonical forms to toroidal boundaries.
The weight they used there is now called \textit{geometric weight}, which emphasizes the power of the canonical bundle of which a given pluricanonical form is a section.
Explicitly, the geometric weight $k$ corresponds to the arithmetic weight $10k$ in our case.

Combining (\ref{eq:canonical_bundle_ball_quotients}), (\ref{eq:canonical_bundle_bb}) and (\ref{eq:canonical_bundle_toroidal}), we now obtain
\begin{align}
    \label{eq:discrepancy_pi_H}
    \pi^*\overline{H}^{\BB} = \overline{H}^{\tor} -18T.
\end{align}

\begin{rem}
\label{rem:compnormalbundle}
Allcock \cite[Theorem 7.1]{All00} constructed an automorphic form $\Psi_{\mathrm{A}}$, a third root of $\Psi_1$ (Example \ref{ExampleBorcherdsProd} (2)), of weight 44 with respect to the group $\Gamma$, vanishing exactly on the ramification 
divisor $\H \subset \B^9$ with multiplicity 1 (see also Subection \ref{subsection:allcock}).
For the Baily-Borel compactification, this implies that
\[44\L=\frac{1}{6}\overline{H}^{\BB},\]
Combining this with (\ref{eq:canonical_bundle_bb}), we obtain  
\begin{align}
    \label{eq:canonical_bundle_bb in terms of L}
    K_{\overline{\B^9/\Gamma}^{\BB}} = -210\L.
\end{align}
Now, we can deduce
\[44\L=\frac{1}{6}\overline{H}^{\tor} - 3T\]
on the toroidal compactification by (\ref{eq:discrepancy_pi_H}).
Substituting this to (\ref{eq:canonical_bundle_toroidal}), we have
\begin{align}
    \label{eq:canonical_bundle_toroidal in terms of L}
    K_{\overline{\B^9/\Gamma}^{\tor}} = -210\L-16T.
\end{align}
\end{rem}

\section{Log minimal model program and semi-toroidal compactifications}
\label{section:semi-toroidal_compactifications}
\subsection{Review of the log minimal model program}
In this section, we will briefly recall the relevant notions of LMMP which are necessary for the discussion of semi-toroidal compactifications later on.
For the basic definitions see \cite[Section 4]{Fuj17} and also \cite[Subsection 4.3]{HM25}, where we use this in an analogous situation. 

Before stating our claims and giving the proofs, we will review the various notions of singularities.
For a pair $(X,\Delta_X)$, we use the notion of \textit{kawamata log terminal}, \textit{divisorial log terminal}, \textit{purely log terminal} and \textit{log canonical} as defined in \cite[Definitions 2.34, 2.37]{KM98}.
In short, we often call these \textit{klt}, \textit{dlt}, \textit{plt} and \textit{lc}.
If $\Delta_X=0$, then klt, dlt and plt are called \textit{log terminal}.
A pair $(X,\Delta_X)$ is said to be \textit{quasi-divisorial log terminal} if there is a finite surjective morphism $f:Y\to X$ for a dlt pair $(Y,\Delta_Y)$ such that $K_Y+\Delta_Y=f^*(K_X+\Delta_X)$.  
    The reason why we introduce this notion is that in the above situation, given $f$ as stated, $(X,\Delta_X)$ is klt (resp. plt, lc) if and only if $(Y,\Delta_Y)$ is klt (resp. plt, lc),  but a similar argument does not hold for dlt; see \cite[Corollary 2.43]{Kol13}.

Now, we start with the singularities of the compactifications of the ball quotient $\B^9/\Gamma$.
    Let $T$ be the exceptional divisor of the blow-up $\pi:\overline{\B^9/\Gamma}^{\tor}\to\overline{\B^9/\Gamma}^{\BB}$.
    A classical result of Mumford implies that the blow-up 
\begin{align}
\label{morphism:bq_log_crepant}
    \pi:\left(\overline{\B^9/\Gamma}^{\tor}, \frac{5}{6}\overline{H}^{\tor}+T\right)\to\left(\overline{\B^9/\Gamma}^{\BB}, \frac{5}{6}\overline{H}^{\BB}\right)
\end{align}
is log-crepant.
\begin{prop}
    \label{prop:singularity_ball_quotient}
    \begin{enumerate}
        \item  The varieties $\overline{\B^9/\Gamma}^{\BB}$ and $\overline{\B^9/\Gamma}^{\tor}$ are $\Q$-Gorenstein.
    Moreover, $\overline{\B^9/\Gamma}^{\tor}$ is locally $\Q$-factorial.
    In particular, $K_{\overline{\B^9/\Gamma}^{\BB}} + \frac{5}{6}\overline{H}^{\BB}$ and  $K_{\overline{\B^9/\Gamma}^{\tor}}+\frac{5}{6}\overline{H}^{\tor}+T$ are $\Q$-Cartier.
    \item $(\overline{\B^9/\Gamma}^{\BB}, \frac{5}{6}\overline{H}^{\BB})$ has log canonical singularities.
    \item $(\overline{\B^9/\Gamma}^{\tor}, \frac{5}{6}\overline{H}^{\tor}+T)$ has quasi-divisorial log terminal singularities.
    \item $\overline{\B^9/\Gamma}^{\BB}$ has log terminal singularities.
    Moreover, $\B^9/\Gamma$ is rationally connected and $(\B^9/\Gamma)\setminus H$ is affine.
    \end{enumerate}
\end{prop}
\begin{proof}
     We take a neat normal subgroup $\Gamma'\lhd \Gamma$ so that $\overline{\B^9/\Gamma'}$ is smooth and there is a finite surjective Galois cover $f:\overline{\B^9/\Gamma'}^{\tor}\to\overline{\B^9/\Gamma}^{\tor}$
     with covering group $\Gamma/\Gamma'$.
     By $T'$ we denote the toroidal boundary of the neat cover $\overline{\B^9/\Gamma'}^{\tor}$, which is simple normal crossing. 
   
    (1) Since $\overline{\B^9/\Gamma'}^{\tor}$ is smooth, 
    it follows from \cite[Lemma 5.16]{KM98} that $\overline{\B^9/\Gamma}^{\tor}$ is $\Q$-factorial.
    For the Baily-Borel compactification, (\ref{eq:canonical_bundle_bb in terms of L}) implies that $K_{\overline{\B^9/\Gamma}^{\BB}}$ is $\Q$-Cartier because $\L$ is a $\Q$-line bundle.

    (2) Since there is a finite surjective morphism $\overline{\B^9/\Gamma'}\to\overline{\B^9/\Gamma}$,  \cite[Corollary 2.43]{Kol13} implies that the pair $(\overline{\B^9/\Gamma'}^{\BB},0)$ 
    is log canonical if and only if $(\overline{\B^9/\Gamma}^{\BB},\frac{5}{6}\overline{H}^{\BB})$ is.
    Hence, it suffices to prove that $(\overline{\B^9/\Gamma'}^{\BB},0)$ is log canonical.
    This follows from the existence of the log-crepant resolution $(\overline{\B^9/\Gamma}^{\tor},T')$ in a similar way as (\ref{morphism:bq_log_crepant}).
    
    (3) By the choice of a suitable finite cover $\overline{\B^9/\Gamma'}^{\tor}$ the pair $(\overline{\B^9/\Gamma'}^{\tor}, T')$ is dlt. To prove the claim it is enough to show that 
    \begin{equation}
    \label{equ:comparingKnormalcover}
   K_{\overline{\B^9/\Gamma'}^{\tor}} + T'=   f^*\left(K_{\overline{\B^9/\Gamma}^{\tor}} + \frac{5}{6}\overline{H}^{\tor} + T\right).
   \end{equation}
    For this we denote by $\overline{H'}^{\tor}\subset\overline{\B^9/\Gamma'}^{\tor}$ the closure of the ramification divisor (with index 6) mapping to $\overline{H}^{\tor}\subset\overline{\B^9/\Gamma}^{\tor}$.
    We then have $f^*\overline{H}^{\tor}=6\overline{H'}^{\tor}$ and $f^*T=nT'$ where $n$ is the ramification index of the toroidal divisors (which are all the same since we have chosen a normal subgroup).
    The Riemann-Hurwitz formula then  tells us that 
    \[K_{\overline{\B^9/\Gamma'}^{\tor}}=f^*K_{\overline{\B^9/\Gamma}^{\tor}} + (n-1)T' + 5\overline{H'}^{\tor}=f^*\left(K_{\overline{\B^9/\Gamma}^{\tor}} + \frac{5}{6}\overline{H}^{\tor} + T\right) -T',\]
    thus showing (\ref{equ:comparingKnormalcover}).

(4) In our case, by the moduli interpretation, the closure of the discriminant divisor $\overline{H}^{\BB}$ contains the unique Baily-Borel cusp.
Hence, $\overline{\B^9/\Gamma}^{\BB}$ contains no naked cusp in the sense of \cite[Definition 2.7]{MO23}.
Note that the automorphic form given by Allcock as in Remark \ref{rem:compnormalbundle} satisfies the assumption in \cite[Theorem 2.4 (1)]{MO23}.
Hence, \cite[Corollaries 2.8, 2.10]{MO23} implies that $\overline{\B^9/\Gamma}^{\BB}$ is log terminal and rationally connected, and $(\B^9/\Gamma)\setminus H$ is affine.
\end{proof}
We note that in general the Baily-Borel compactifications have bad singularities along the cusps.
However, in some special cases, they can be explicitly computed; the moduli space of cubic surfaces is isomorphic to $\P(1,2,3,4,5)$ and the cusp is even a smooth point, see the proof of \cite[Lemma 3.1]{CMGHL23b}.
On this basis, we propose the following problem.
\begin{que}
    \label{que:quotient_singularities}
    Does $\overline{\B^9/\Gamma}^{\BB}$ have finite quotient singularities?
\end{que}
Not only is this interesting for the geometry of the ball quotient, but there is also an application concerning derived algebraic geometry; see Subsection \ref{subsection:derived_category}.
Now, we can restate the results of Baily-Borel and Mumford, specializing them to our case.

\begin{thm}
\label{thm:bb_tor_mmp}
The Baily-Borel compactification and the toroidal compactification of $\B^9/\Gamma$ have the following properties in LMMP:
    \begin{enumerate}
        \item $(\overline{\B^9/\Gamma}^{\BB}, \frac{5}{6}\overline{H}^{\BB})$ is the log canonical model and a minimal model of itself.
        \item $(\overline{\B^9/\Gamma}^{\tor}, \frac{5}{6}\overline{H}^{\tor}+T)$ is a log minimal model of itself.
    \end{enumerate}
\end{thm}
\begin{proof}
These statements follow from the classical results of Baily-Borel and Mumford; see \cite[Corollaries 3.3, 3.5]{Ale96}.
In other words, the first (resp. second) statement follows from the projective construction of the Baily-Borel compactification, which implies the ampleness of $\L$ (resp.  the fact that $\pi^*\L$ coincides with $K_{\overline{\B^9/\Gamma}^{\tor}}+\frac{5}{6}\overline{H}^{\tor}+T$ on the toroidal compactification, which can be obtained as a corollary of Hirzebruch's proportionality principle by Mumford).
The coefficient $5/6$ comes from the lattice-theoretic computation discussed in Remark \ref{rem:lattice_theoretic_branch}, see also (\ref{eq:canonical_bundle_bb}), (\ref{eq:canonical_bundle_toroidal}), and \cite[Theorems 3.3, 3.5]{Ale96}.
We already know that the two pairs are log canonical pairs by Proposition \ref{prop:singularity_ball_quotient}.
\end{proof}

Since we have a dominant rational map from $\Sym^{12}(\P^1) \cong \P^{12}$ to the moduli space of unordered $12$ points on $\P^1$ the ball quotient $\overline{\B^9/\Gamma}^{\BB}$ is clearly unirational 
and thus has
Kodaira dimension $-\infty$. 
In fact, it is even known to be rational \cite{Kat84}.
In the following remark, we give a different argument for the statement concerning the Kodaira dimension by using Allcock's automorphic form, see also \cite[Theorem 4.4]{GHS08} where a similar method was employed.     
\begin{rem} One can use Allcock's automorphic form to show  $\kappa(\overline{\B^9/\Gamma}^{\BB})= \kappa(\overline{\B^9/\Gamma}^{\tor})= - \infty$. We prove this by contradiction. Let $X$ be a smooth projective model 
of $\overline{\B^9/\Gamma}^{\BB}$ which admits a non-trivial  $k$-fold pluricanonical form. By the usual correspondence between automorphic forms and pluricanonical forms on locally symmetric varieties, such a form comes from a form $F_{10k}(dZ)^k$ on $\B^9$ where $dZ$ is the standard volume element and $F_{10k}$ is a cusp form of weight $10k$ vanishing of sufficiently high order at the cusp. 
Since we have ramification of order $6$ along $H$, it follows that $F_{10k}$ must vanish of order $5k$ along $H$, see also the expression (\ref{eq:canonical_bundle_bb}) for the canonical bundle.
On the other hand, we know from Remark \ref{rem:compnormalbundle} that Allcock's automorphic form $\Psi_{\mathrm{A}}$, 
which has weight 44, vanishes exactly along $H$ with vanishing order 1. Now dividing $F_{10k}$ by $\Psi_{\mathrm{A}}^{5k}$ 
we obtain an automorphic form of negative weight, a contradiction.  
\end{rem}

Next, we shall move on to the Kirwan blow-up.
Let $\Delta$ be the exceptional divisor of the Kirwan blow-up $f:\M^{\K}\to\M^{\GIT}$.
Further, let $\D$ be the discriminant divisor on $\M^{\GIT}$ and let $\widetilde{\D}$ be its strict transform.
We first prove an analog of Proposition \ref{prop:singularity_ball_quotient} (1).
      \begin{prop}
      \label{prop:qgorenstein_kirwan}
          The Kirwan blow-up $\M^{\K}$ is locally $\Q$-factorial and log terminal.    In particular, $K_{\M^{\K}}+\frac{5}{6}\widetilde{\D}+\Delta$ is $\Q$-Cartier.
      \end{prop}   
      \begin{proof}
          By the construction of $\M^{\K}$, all points in $\M^{\K}$ are stable.
    Hence, the stabilizer of any point is finite and thus,  locally, $\M^{\K}$ has finite quotient singularities.
    This implies that $\M^{\K}$ is locally $\Q$-factorial by \cite[Proposition 5.15]{KM98}.
    Also, it is known that a finite quotient singularity is log terminal.
      \end{proof}
Below, for a pair $(Y,\Delta_Y)$ and a birational map $f:X\dashrightarrow Y$, we denote by $a(E_i, Y,\Delta_Y)$ the \textit{discrepancy} as given by 
\[K_X + (f^{-1})_*\Delta_Y = f^*(K_Y+\Delta_Y) + \sum_i a(E_i, Y,\Delta_Y)E_i \]
where the $E_i$ are the exceptional divisors of $f$.

Below, we shall compute discrepancies explicitly.
We denote by $a(\Delta, \M^{\GIT}, \frac{5}{6}\D)\in\Q$ the discrepancy of the Kirwan blow-up $f:(\M^{\K}, \frac{5}{6}\widetilde{\D} + \Delta)\to (\M^{\GIT}, \frac{5}{6}\D)$:
\begin{align}
\label{eq:def_log_discrepancy}
    K_{\M^{\K}} + \frac{5}{6}\widetilde{\D} = f^*\left(K_{\M^{\GIT}} + \frac{5}{6}\D\right) + a\left(\Delta, \M^{\GIT}, \frac{5}{6}\D\right)\Delta.
\end{align}
One can characterize semi-toroidal compactifications by using discrepancies.
\begin{prop}
\label{prop:discrepancy_semi_toric}
    If the discrepancy $a\left(\Delta, \M^{\GIT}, \frac{5}{6}\D\right)>-1$, then $(\M^{\K}, \frac{5}{6}\widetilde{\D} + \Delta)$ is not a log minimal model (of itself).
    Hence, $\M^{\K}$ is not a semi-toroidal compactification.
\end{prop}

\begin{proof}
By a similar computation, as in the proof of \cite[Proposition 4.7]{HM25}, we can show that if the discrepancy is greater than $-1$, then $(\M^{\K}, \frac{5}{6}\widetilde{\D} + \Delta)$ is not 
the log minimal model of itself.
This assumption exactly corresponds to $a\left(\Delta, \M^{\GIT}, \frac{5}{6}\D\right)>-1$.
Hence we conclude that $\M^{\K}$ is not a semi-toroidal compactification by \cite[Theorem 3.1]{Oda22}.
\end{proof}

We conclude this subsection with a discussion concerning the relation between discrepancies and log $K$-equivalence.
It is known that any two log minimal models of a given variety are (log) $K$-equivalent by the negativity lemma \cite[Lemma 2.3.26]{Fuj17}; see also \cite[Lemma 4.3.2]{Fuj17}.
Here, we shall discuss the converse.

\begin{prop}
\label{prop:log_K_equiv_criterion}
Let $(X,A)$ and $(Y,B)$ be two pairs.
We assume that there is a birational map 
$p:Y\dashrightarrow X$ with $p_*K_Y=K_X$, $p_*B=A$ and $K_Y+B=p^*(K_X+A)+E$ for a divisor $E\neq 0$ in $\Pic(Y)\otimes\Q$.
Then the two pairs are not log $K$-equivalent.
\end{prop}
\begin{proof}
    For any birational morphisms $f_X:Z\to X$ and $f_Y:Z\to Y$ from a projective variety $Z$, we have 
    \[f_Y^*(K_Y+B)=f_X^*(K_X+A) + f_Y^*E.\]
    Since $f_Y$ is a proper birational morphism and $Y$ is normal, it follows that $f_{Y,*}\OO_Z=\OO_Y$.
    This implies that the induced morphism between the Picard groups $f_Y^*:\Pic(Y)\to\Pic(Z)$ is injective by the assumption on the normality of $Y$.
    By our assumption on $E$, it follows that $f_Y^*E\neq 0$.
    Therefore, the pairs $(X,A)$ and $(Y,B)$ are never log $K$-equivalent,
\end{proof}

\begin{cor}
\label{cor:log_K_equiv_minimal}
 Let  $(X, A)$ and $(Y,B)$ be pairs.
 Assume the following:
 \begin{itemize}
     \item[(A)] $(X,A)$ is a log minimal model of itself.
     \item[(B)]     There is a birational map $f:Y\dashrightarrow X$ so that $p_*K_Y=K_X$,  $p_*B=A$, and $p$ is not log-crepant.
 \end{itemize}
Then, these two pairs are not log $K$-equivalent.
\end{cor}
\begin{proof}
    We can write 
    \[K_Y+B=p^*(K_X+A)+ \sum m_{\alpha}E_{\alpha}\]
    where the $E_{\alpha}$ are exceptional divisors.
    By an application of the negativity lemma \cite[Lemma 4.3.2]{Fuj17}, since $(X,A)$ is a log minimal model (A), we see that $E=\sum m_{\alpha}E_{\alpha}$ is an effective divisor, i.e. all coefficients $m_{\alpha} \geq 0$.
  By (B), we also know that $E$ cannot be the trivial divisor, i.e. at least one $m_{\alpha} > 0$,  and hence $E \neq 0 \in \Pic(Y)\otimes\Q$. 
    We can now conclude the claim from Proposition \ref{prop:log_K_equiv_criterion}.
\end{proof}

\subsection{Proof of the main result}
\label{subsection:semi-toroidal}
By the results of Alexeev-Engel \cite{AE23} and Odaka \cite{Oda22}, 
semi-toroidal compactifications are precisely the compactifications lying between the toroidal compactification and the Baily--Borel compactification.
We use this characterization to show that $\M^K$ is not semi-toroidal.
For this, we shall first compute the discrepancy of the Kirwan blow-up.
Let $f':\widetilde{(\P^{12})^{\mathrm{ss}}} \to (\P^{12})^{\mathrm{ss}}$ be the Kirwan blow-up (before taking quotients). Its center is $\SL_2(\C)\cdot\{p_{6,6}\}$.
By Kirwan's theory, the action of $\SL_2(\C)$ lifts to $\widetilde{(\P^{12})^{\mathrm{ss}}}$ and all points are stable.
Let $\D'$ (resp. $\Delta'$) be the discriminant (resp. exceptional) divisor on $(\P^{12})^{\mathrm{ss}}$ (resp. $\widetilde{(\P^{12})^{\mathrm{ss}}}$). 
Taking GIT quotients, we denote by $\D$ (resp. $\Delta$) the corresponding discriminant (resp. exceptional) divisor on $\M^{\GIT}$ (resp. $\M^{\K}$). 
\begin{thm}
\label{thm:not_semi_toric}
The following holds:
\begin{itemize} 
\item[(1)]  The discrepancy $a(\Delta, \M^{\GIT}, \frac{5}{6}\D) = \frac{2}{3}>-1$ and hence  $\M^{\K}$ is not a semi-toroidal compactification.
\item[(2)]  The two pairs $(\M^{\K}, \frac{5}{6}\widetilde{\D} + \Delta)$ and $(\overline{\B^9/\Gamma}^{\tor}, \frac{5}{6}\overline{H}^{\tor} + T)$ are not log $K$-equivalent.
 \end{itemize}   
\end{thm}
\begin{proof}
From Lemma \ref{lem:Luna_slice} (2), we find that $c=10$ in the notation of  \cite[Lemma 6.4]{CMGHL23b} and hence we obtain
\begin{align*}
    K_{\widetilde{(\P^{12})^{\mathrm{ss}}}} 
    & = f^{'*}K_{(\P^{12})^{\mathrm{ss}}} + 9\Delta'\\
    & = f^{'*}\left(K_{(\P^{12})^{\mathrm{ss}}} + \frac{5}{6}\D'\right) + \frac{2}{3}\Delta' -\frac{5}{6}\widetilde{\D'}.
\end{align*}
In other words, 
\[K_{\widetilde{(\P^{12})^{\mathrm{ss}}}} + \frac{5}{6}\widetilde{\D'}= f^{'*}\left(K_{(\P^{12})^{\mathrm{ss}}} + \frac{5}{6}\D'\right) + \frac{2}{3}\Delta' .\]
Combined with \cite[Remark 6.7]{CMGHL23b} and a similar discussion as in the proof of \cite[Proposition 4.7]{HM25}, this calculation descends to the GIT quotient and implies 
that the discrepancy $a(\Delta, \M^{\GIT}, \frac{5}{6}\D) = \frac{2}{3}$ and hence $\M^{\K}$ is not a semi-toroidal compactification by Proposition \ref{prop:discrepancy_semi_toric}. This proves item (1).

Item (2) essentially follows from item (1), Theorem \ref{thm:bb_tor_mmp} (1) and Corollary \ref{cor:log_K_equiv_minimal}.
Indeed, item (1) implies that (B) in Corollary \ref{cor:log_K_equiv_minimal} are satisfied for two pairs $(\M^{\K}, \frac{5}{6}\widetilde{\D} + \Delta)$ and $(\M^{\GIT}, \frac{5}{6}\D)$.
Theorem \ref{thm:bb_tor_mmp} (1) implies (A) for the GIT pair.
Now, we can apply Corollary \ref{cor:log_K_equiv_minimal} to $(\M^{\K},\frac{5}{6}\widetilde{\D}+\Delta)$ and $(\M^{\GIT}, \frac{5}{6}\D)$ with the exceptional divisor $E=\Delta$.
It follows that these two pairs are not log $K$-equivalent.
Here we note that by the transitivity of $K$-equivalence, the three pairs $(\M^{\GIT}, \frac{5}{6}\D), (\overline{\B^9/\Gamma}^{\BB}, \frac{5}{6}\overline{H}^{\BB})$, $(\overline{\B^9/\Gamma}^{\tor}, \frac{5}{6}\overline{H}^{\tor}+T)$ are log $K$-equivalent.
Therefore, we conclude that $(\M^{\K}, \frac{5}{6}\widetilde{\D} + \Delta)$ and $(\overline{\B^9/\Gamma}^{\tor}, \frac{5}{6}\overline{H}^{\tor}+T)$ are not log $K$-equivalent.
\end{proof}

As a corollary of the computation of the discrepancies we can now clarify the singularities; compare with Proposition \ref{prop:singularity_ball_quotient}.

\begin{cor}
    \label{cor:singularity_kirwan}
    \begin{enumerate}
        \item $(\M^{\K}, \frac{5}{6}\widetilde{\D}-\frac{2}{3}\Delta)$ is log canonical, but not purely log terminal.
    \item $(\M^{\K}, \frac{5}{6}\widetilde{\D}+\Delta)$ is not log canonical.
    \end{enumerate}
\end{cor}
\begin{proof}    
   (1) Our assertion follows from the following statements:
    \begin{itemize}
        \item[(A)]  $(\M^{\GIT}, \frac{5}{6}\D)$ is log canonical. (Proposition \ref{prop:singularity_ball_quotient} (2)).
        \item[(B)]  The discrepancy $a(\Delta, \M^{\GIT}, \frac{5}{6}\D)$ is equal to $\frac{2}{3}$. (Theorem \ref{thm:not_semi_toric} (1)).
    \end{itemize}
    We shall show that $(\M^{\K}, \frac{5}{6}\widetilde{\D}-\frac{2}{3}\Delta)$ is log canonical, but not purely log terminal.
    Here, we adopt the definition in \cite[Definition 2.34]{KM98}, and hence negative coefficients of boundaries are allowed.
To prove that $(\M^{\K}, \frac{5}{6}\widetilde{\D}-\frac{2}{3}\Delta)$ is log canonical but not purely log terminal,
let us first take any resolution of singularities $g:X\to\M^{\K}$ and denote its exceptional divisors by $\{E_i\}$ so that 
\begin{align}
\label{eq:desingularization_of_kirwan}
    K_X + (g^{-1})_*\left(\frac{5}{6}\widetilde{\D}-\frac{2}{3}\Delta\right)=g^*\left(K_{\M^{\K}} + \frac{5}{6}\widetilde{\D}-\frac{2}{3}\Delta\right) + \sum_ia_iE_i
\end{align}
where $a_i=a(E_i, \M^{\K}, \frac{5}{6}\widetilde{\D}-\frac{2}{3}\Delta)$.
Below, we shall prove that $a_i \geq -1$ for any $i$ and $a_j=-1$ for at least one $j$.

By considering the composition of $g:X\to\M^{\K}$ and $f:\M^{\K}\to\M^{\GIT}$, item (A) leads us to the following representation 
\begin{align}
    \label{eq:B}
    K_X + (f\circ g)^{-1}_*\frac{5}{6}\D &=(f\circ g)^*\left(K_{\M^{\GIT}}+\frac{5}{6}\D\right)+\sum_ib_iE_i + b(g^{-1})_*\Delta\\
    & = g^*\left(K_{\M^{\K}}+\frac{5}{6}\widetilde{\D}-\frac{2}{3}\Delta\right) + \sum_ib_iE_i + b(g^{-1})_*\Delta\notag
\end{align}
where $b_i=a(E_i, \M^{\GIT}, \frac{5}{6}\D)\geq -1$ and $b=a((g^{-1})_*\Delta, \M^{\GIT}, \frac{5}{6}\D)\geq -1$ by the claim that $(\M^{\GIT}, \frac{5}{6}\D)$ is log canonical. 
The last equation follows from
\begin{align}
    \label{eq:log_discrepanci_D}
    K_{\M^{\K}} + \frac{5}{6}\widetilde{\D} = f^*\left(K_{\M^{\GIT}}+\frac{5}{6}\D\right) + \frac{2}{3}\Delta,
\end{align}
which is deduced from (\ref{eq:def_log_discrepancy}) and Theorem \ref{thm:not_semi_toric} (1).

We write $g^*\widetilde{\D}=\widetilde{\D}'+\sum_i d_iE_i$ and $g^*\Delta=\widetilde{\Delta}+ \sum_i e_iE_i$  where $\widetilde{\D}'$ and $\widetilde{\Delta}$ are the strict transform of  $\widetilde{\D}$ and $\Delta$ via $g$.
Comparing the coefficients of $E_i$ in (\ref{eq:desingularization_of_kirwan}) and (\ref{eq:B}), it follows that $a_i=b_i$, which implies that $(\M^{\K}, \frac{5}{6}\widetilde{\D}-\frac{2}{3}\Delta)$ is log canonical because $b_i\geq -1$.
Now, since $\Delta$ is singular by the local description, there is some $j$ such that $g(E_j)\subset\Delta$.
In other words, for such an index $j$, we have that $f\circ g(E_j)$ is the unique Baily-Borel cusp in $\overline{\B^9/\Gamma}^{\BB}$.
We also know from \cite[Lemma 2.9 (1)]{MO23} that the log canonical center of $(\overline{\B^9/\Gamma}^{\BB}, \frac{5}{6}\D)$ is the unique Baily-Borel cusp.
This implies that $a_j=b_j=-1$ and hence the pair $(\M^{\K}, \frac{5}{6}\widetilde{\D}-\frac{2}{3}\Delta)$ is not purely log terminal.

(2) By \cite[Lemma 2.27]{KM98} we now obtain $a(E_i, \M^{\K}, \frac{5}{6}\widetilde{\D}-\frac{2}{3}\Delta) \geq a(E_i, \M^{\K}, \frac{5}{6}\widetilde{\D}+\Delta)$ in general.
Setting $i=j$ in this inequality, \cite[Lemma 2.27]{KM98} also implies that the inequality is strict. 
However, we have already seen that $a(E_j, \M^{\K}, \frac{5}{6}\widetilde{\D}-\frac{2}{3}\Delta)=-1$.
Therefore, it follows that $a(E_i, \M^{\K}, \frac{5}{6}\widetilde{\D}+\Delta)< -1$ and thus $(\M^{\K}, \frac{5}{6}\widetilde{\D}+\Delta)$ is not log canonical.

\end{proof}
\begin{rem}
\label{rem:singularities}
Here we summarize some observations surrounding Propositions  \ref{prop:singularity_ball_quotient}, \ref{prop:qgorenstein_kirwan} and Corollary \ref{cor:singularity_kirwan}, concerning the singularities of the compactifications of modular varieties.
\begin{enumerate}
    \item We remark that toroidal compactification $\overline{\B^n/\Gamma}^{\tor}$ of ball quotients are always locally $\Q$-factorial since they have finite quotient singularities. This follows from the same proof as given in 
    Proposition \ref{prop:singularity_ball_quotient} (1).   
    \item Similarly, Kirwan blow-ups $\M^{\K}$ of GIT-quotients $\M^{\GIT}$  are always locally $\Q$-factorial. The reason is that the  Kirwan blow-up replaces 
    the  singularities of the GIT quotient with finite quotient singularities, since, in the course of the 
    consecutive partial blow-ups, all polystable points are replaced by stable points
    \item The last part of Propositions \ref{prop:singularity_ball_quotient} (1) and \ref{prop:qgorenstein_kirwan} can be also deduced from (\ref{eq:canonical_bundles_Kirwan}) and (\ref{eq:canonical_bundle_ball_quotients}) respectively. This gives $\Q$-Gorensteiness directly and does not require $\Q$-factoriality.
    \item The last part of Proposition \ref{prop:singularity_ball_quotient} (4) is an analog of the quasi-affiness of the moduli space of $K3$ surfaces \cite{BKPS98} and cubic surfaces \cite[Section 2.7]{DvGK05}.
    \item Propositions \ref{prop:qgorenstein_kirwan} and  Corollary \ref{cor:singularity_kirwan} above shed some light on the properties of the Kirwan blow-up in connection with LMMP. 
    By construction, the singularities of $\M^{\K}$ are typically better than those of $\M^{\GIT}$.
   However, the pair $(\M^{\K}, \frac{5}{6}\widetilde{\D} + \Delta)$ behaves worse than the log canonical pair $(\M^{\GIT}, \frac{5}{6}\D)$, which means that the pair $(\M^{\K}, \frac{5}{6}\widetilde{\D} + \Delta)$ is 
   not well behaved from the point of view of LMMP. This may be due to the fact that the factor $5/6$ arises from the ball quotient picture, but has no natural explanation from a GIT point of view. 
   Note that this is different from the case of the toroidal compactification, where the factors of the boundary $\frac{5}{6}\overline{H}^{\tor} + T$ arise naturally from the ball quotient picture and the role of the toroidal boundary 
   in the theory of Mumford et al. concerning toroidal compactifications. 
   The resulting pair is quasi-dlt and can be treated in the framework of LMMP.
   We also refer to the example discussed in \cite[Section 2]{CMGHL23b} for a related example in another setting.      
    \item   More generally, if there is a reflective automorphic form that vanishes exactly on the branch divisors with multiplicities according to the ramification order, called a \textit{special reflective modular form} in \cite{MO23}, then we can prove that the associated Baily-Borel compactification is $\Q$-Gorenstein, as in \cite[Theorem 2.4]{MO23} and the proof therein. The authors do not know of good criteria which 
    guarantee that the Baily-Borel compactification is $\Q$-factorial. We note, however, that special reflective modular forms also exist in the case of cubic surfaces and 8 points on $\P^1$, as treated in  \cite{CMGHL23b, HM25}. Moreover, in the 
    case of moduli of cubic surfaces, the (unique) cusp is even a smooth point.    
\end{enumerate}   
\end{rem}
\vspace{2cm}

\section{Two applications: another proof and failure of stacky $D$-equivalence}
\label{section:two_applications}
\subsection{Another proof of the non-liftability of the period map}
\label{subsec:Applications of the main result}
As an application of Theorem \ref{thm:not_semi_toric}, we now give another proof of our claim concerning the non-liftability of the Deligne-Mostow period map.
\begin{cor}
\label{cor:nonliftabilitysecondproof}
The inverse  $\phi^{-1}:  \overline{\B^9/\Gamma}^{\BB} \to \M^{\GIT}$ of the Deligne-Mostow isomorphism does not lift
to a morphism $g^{-1}: \overline{\B^9/\Gamma}^{\tor} \to \M^{\K}$.
\end{cor}

\begin{proof}
Assume $g^{-1}: \overline{\B^9/\Gamma}^{\tor} \to \M^{\K}$  is a morphism. Then the map $\pi: \overline{\B^9/\Gamma}^{\tor} \to \overline{\B^9/\Gamma}^{\BB}$ factors 
through $\M^{\K}$. Now this implies that $\M^{\K}$ is a semi-toroidal compactification from the proof of \cite[Theorem 3.24]{AEH24}, which in turn is based on \cite[Theorem 5.14]{AE23}.
More precisely, this says that 
a normal compactification of $\B^9/\Gamma$, is a semi-toroidal compactification if and only if the map $\pi: \overline{\B^9/\Gamma}^{\tor} \to \overline{\B^9/\Gamma}^{\BB}$ factors through it.
This, however, contradicts  Theorem  \ref{thm:not_semi_toric}.
\end{proof}

\begin{rem}
    Here we add some comments on various logical implications.
    Corollary \ref{cor:nonliftabilitysecondproof}, a partial statement of Theorem \ref{thm:lift}, is in fact equivalent to the claim
    of Theorem \ref{thm:not_semi_toric} asserting that $\M^{\K}$ is not a semi-toroidal compactification.
    This is because  \cite[Theorem 5.14]{AE23} provides a necessary and sufficient condition for a compactification to be semi-toroidal.
    Hence we can obtain Theorem \ref{thm:not_semi_toric} (the latter part of) (1), (2) by a Luna slice calculation only, without an explicit calculation of discrepancies.
\end{rem}

\subsection{Failure of stacky $D$-equivalence}
\label{subsection:derived_category}
We shall now apply our results to the derived algebraic geometry of $\M^{\K}$ and $\overline{\B^9/\Gamma}^{\tor}$.
We refer to Subsections \ref{subsec:main_results} and  \ref{subsec:concepts} for a discussion of  $DK$-equivalence and the terminology concerning derived categories.
Here we shall use that we have already shown in Theorem \ref{thm:not_k_equiv} that the two varieties $\M^{\K}$ and $\overline{\B^9/\Gamma}^{\tor}$ are not $K$-equivalent.
We start with the following lemma. 
\begin{lem}
\label{lem:kodaira_dimension_of_ball_quotient}
    The $\Q$-line bundle $-K_{\overline{\B^9/\Gamma}^{\BB}}$ is ample and $-K_{\overline{\B^9/\Gamma}^{\tor}}$ is big.
\end{lem}
\begin{proof}
    By (\ref{eq:canonical_bundle_bb in terms of L}) we know that  $-K_{\overline{\B^9/\Gamma}^{\BB}}= 210 \L$, where $\L$ is the automorphic $\Q$-line bundle, which is ample on $\overline{\B^9/\Gamma}^{\tor}$.
    Hence $-K_{\overline{\B^9/\Gamma}^{\BB}}$ is ample.
For the toroidal compactification, equation (\ref{eq:canonical_bundle_toroidal in terms of L}) leads to 
  \[-K_{\overline{\B^9/\Gamma}^{\tor}} = 210\L+16T,\]
  which implies that $-K_{\overline{\B^9/\Gamma}^{\tor}}$ is big because $210\L$ is big on $\overline{\B^9/\Gamma}^{\tor}$, and the sum of a big divisor and an effective divisor is also big.
\end{proof}
    
We proved in Theorem \ref{thm:not_k_equiv} that the two varieties $\M^{\K}$ and $\overline{\B^9/\Gamma}^{\tor}$ are not $K$-equivalent.
Hence, by the generalized $DK$-conjecture (Conjecture \ref{mainconj:generalized DK conjecture}), it is expected that $D(\mathcal{X}(\M^{\K}))\not\cong D(\mathcal{X}(\overline{\B^9/\Gamma}^{\tor}))$.
This is indeed correct, as we shall show now.

\begin{thm}
\label{thm:derived_equivalent}
The two categories $D(\mathcal{X}(\M^{\K}))$ and $D(\mathcal{X}(\overline{\B^9/\Gamma}^{\tor}))$ are not equivalent as triangulated categories.    
In other words, the two varieties $\M^{\K}$ and $\overline{\B^9/\Gamma}^{\tor}$ are not stacky $D$-equivalent.
\end{thm}
\begin{proof}
The claim follows from Kawamata's work, combined with Theorem \ref{thm:not_k_equiv} and Lemma \ref{lem:kodaira_dimension_of_ball_quotient}.
Since $\M^{\K}$ and $\overline{\B^9/\Gamma}^{\tor}$ have finite quotient singularities, we can define associated smooth stacks $\mathcal{X}(\M^{\K})$ and $\mathcal{X}(\overline{\B^9/\Gamma}^{\tor})$ as in \cite{Kaw04}.
Moreover, we know that $\M^{\K}$ and $\overline{\B^9/\Gamma}^{\tor}$ are not $K$-equivalent (Theorem \ref{thm:not_k_equiv}) and $-K_{\overline{\B^9/\Gamma}^{\tor}}$ is big (Lemma \ref{lem:kodaira_dimension_of_ball_quotient}).
It then follows from \cite[Theorem 7.1 (2)]{Kaw04}, which is an application of Orlov's type representability theorem for varieties with only finite quotient singularities \cite[Theorem 1.1]{Kaw04}, that the two varieties are not stacky $D$-equivalent
\end{proof}
    We remark that Theorem \ref{thm:derived_equivalent} also holds for the moduli spaces of 8 points on $\P^1$ and cubic surfaces, which complements the results of \cite{HM25} and  \cite{CMGHL23b}, combined with Remark \ref{rem:singularities} (1), (2).
In general,  Baily-Borel compactifications need not necessarily have finite quotient singularities at the cusps.
However, in some special situations, this is known to be the case.
In such cases, the derived categories of the associated stacks can be investigated by the same method that we already employed for Theorem \ref{thm:derived_equivalent}.
We note, however, that we do not know whether $\overline{\B^9/\Gamma}^{\BB}$ has finite quotient singularities or not (Question \ref{que:quotient_singularities}).

We illustrate this approach in the case of the moduli space of cubic surfaces, which was investigated in \cite{CMGHL23b}.
The GIT compactification 
\[\M^{\GIT}_{\mathrm{cub}}\defeq \P H^0(\P^3, \OO_{\P^3}(3))/\!/_{\OO(1)}\SL_4(\C)\]
of the moduli space of smooth cubic surfaces is isomorphic to the Baily-Borel compactification of a 4-dimensional ball quotient $\overline{\B^4/\Gamma_{\mathrm{cub}}}^{\BB}$.
By a classical result of invariant theory, $\M^{\GIT}_{\mathrm{cub}}$ is isomorphic to the weighted projective space $\P(1,2,3,4,5)$. Hence the spaces
$\overline{\B^4/\Gamma_{\mathrm{cub}}}^{\BB} \cong \M^{\GIT}_{\mathrm{cub}}$ have only finite quotient singularities and
we can define an associated stack $\mathcal{X}(\overline{\B^4/\Gamma_{\mathrm{cub}}}^{\BB})$ according to \cite{Kaw04}.
Following our standard notation, we shall denote the Kirwan blow-up of $\M^{\GIT}_{\mathrm{cub}}$ by $\M_{\mathrm{cub}}^{\K}$ and the toroidal compactification by $\overline{\B^4/\Gamma_{\mathrm{cub}}}^{\tor}$.
In this case, we can even make a stronger statement than in  Theorem \ref{thm:derived_equivalent} above, including a result on autoequivalence groups of derived categories.
\begin{thm}
\label{thm:derived_cubic}
For the moduli spaces of cubic surfaces, the following holds.
\begin{enumerate}
    \item The three categories $D(\mathcal{X}(\M_{\mathrm{cub}}^{\K}))$, $D(\mathcal{X}(\overline{\B^4/\Gamma_{\mathrm{cub}}}^{\tor}))$ and $D(\mathcal{X}(\overline{\B^4/\Gamma_{\mathrm{cub}}}^{\BB}))$ are all 
    different as triangulated categories.    
In other words, no two of the three varieties $\M_{\mathrm{cub}}^{\K}$, $\overline{\B^4/\Gamma_{\mathrm{cub}}}^{\tor}$ and $\overline{\B^4/\Gamma_{\mathrm{cub}}}^{\BB}$ are stacky $D$-equivalent.
\item The group of autoequivalences $\mathrm{Auteq}(D(\mathcal{X}(\overline{\B^4/\Gamma}_{\mathrm{cub}}^{\BB})))$ is isomorphic to the semidirect product of $\Aut(\overline{\B^4/\Gamma}_{\mathrm{cub}}^{\BB})$ and $\Pic(\mathcal{X}(\overline{\B^4/\Gamma}_{\mathrm{cub}}^{\BB}))\times\Z$.
\end{enumerate}

\end{thm}
\begin{proof}
(1) The claim concerning the first two varieties can be proven in essentially the same way as in Theorem \ref{thm:derived_equivalent}. For this we note that $-K_{\overline{\B^4/\Gamma_{\mathrm{cub}}}^{\tor}}$
is big by \cite[Proposition 5.2.(i)]{CMGHL23b}. 
Since $\overline{\B^4/\Gamma_{\mathrm{cub}}}^{\BB}\cong\P(1,2,3,4,5)$, this space has Picard group has rank 1 and $-K_{\overline{\B^4/\Gamma_{\mathrm{cub}}}^{\BB}}=15 \mathcal{O}_{\P(1,2,3,4,5)}(1)$ is ample.  
By the description of \cite[Corollary 6.8]{CMGHL23b} (resp. \cite[Proposition 5.8]{CMGHL23b}), and the fact that the top self-intersection numbers of the respective exceptional divisors are non-zero (see
\cite[Theorem 2.2]{CMGHL23b} for a precise computation of these numbers), it follows that the pairs 
($\M^{\K}_{\mathrm{cub}}, \overline{\B^4/\Gamma_{\mathrm{cub}}}^{\BB}$) and ($\overline{\B^4/\Gamma_{\mathrm{cub}}}^{\tor}, \overline{\B^4/\Gamma_{\mathrm{cub}}}^{\BB}$) are not $K$-equivalent.
Again, using  \cite[Theorem 7.1 (2)]{Kaw04} this implies that the categories $D(\mathcal{X}(\M_{\mathrm{cub}}^{\K}))\not\cong D(\mathcal{X}(\overline{\B^4/\Gamma_{\mathrm{cub}}}^{\BB}))$ and  $D(\mathcal{X}(\overline{\B^4/\Gamma_{\mathrm{cub}}}^{\tor}))\not\cong D(\mathcal{X}(\overline{\B^4/\Gamma_{\mathrm{cub}}}^{\BB}))$.

(2) We have already remarked that
$-K_{\overline{\B^4/\Gamma_{\mathrm{cub}}}^{\BB}}$ is ample. Hence this is a direct consequence of \cite[Theorem 7.2]{Kaw04}.
\end{proof}
Note that the explicit description of the derived category of the associated stack of the weighted projective spaces, including $D(\mathcal{X}(\overline{\B^4/\Gamma_{\mathrm{cub}}}^{\BB}))$, can be found in \cite[Section 5]{Kaw04}.
It is indeed natural to ask about derived categories for pairs. 
In fact, Kawamata \cite{Kaw05} defined these, generalizing the stacky construction in \cite{Kaw04}, for special classes of klt pairs.
However, our pairs $(\overline{\B^9/\Gamma}^{\tor}, \frac{5}{6}\overline{H}^{\tor} + T)$ and $(\M^{\K}, \frac{5}{6}\widetilde{\D} + \Delta)$ are not klt, and hence a priori, there is no analog of Theorem \ref{thm:derived_equivalent} for these pairs.

\begin{rem}
    \label{rem:absolute_dc}
    In view of the original $DK$-conjecture, not taking associated stacks,
    we can also ask whether the \textit{original}, which means non-stacky, derived categories $D(M^{\K})$ and $D(\overline{\B^9/\Gamma}^{\tor})$ are equivalent as triangulated categories.
    We do not know the answer but would like to comment on the question. 
    \begin{enumerate}
        \item  By \cite[Theorem 2]{Bal11}, generalizing \cite[Theorem 3.5]{BO01} for projective Gorenstein varieties $X$ with mild singularities and $K_X$ or $-K_X$ ample, 
         the derived category $D(X)$ recovers the original variety $X$.
        Unfortunately, the varieties appearing in this paper are $\Q$-Gorenstein and not a priori Gorenstein.
         Nevertheless, one can describe the canonical bundles via automorphic forms.
    The Baily-Borel compactification $\overline{\B^9/\Gamma}^{\BB}$ has ample anti-canonical bundle by (\ref{eq:canonical_bundle_bb in terms of L}).
    On the toroidal compactification, the canonical bundle is also big and can be represented as in (\ref{eq:canonical_bundle_toroidal in terms of L}).
    In this situation, the ampleness of the anti-canonical bundle fails because $\L$ is not ample (only big and nef) and there is an obstruction to ampleness due to the toroidal boundary $T$.
    \item 
    Ballard, Favero and Katzarkov proved that there is a full exceptional collection for $D(\M_{n\cdot\epsilon})$.
    Related to this, Castravet and Tevelev \cite[Question 1.1]{CT20a}  asked whether there is an $S_n$-invariant full exceptional collection in $D(\M_{0,n})$, which was indeed shown to exist in \cite[Theorem 1.1]{CT20b}.
    Moreover, they constructed an $S_{12}$-invariant full exceptional collection for $\M^{\K}_{\ord}$ in \cite[Theorem 1.2]{CT20b}, hence we know that $D(\overline{\B^9/\Gamma}^{\tor})$ has a full exceptional collection.
    If we can show that there is no full exceptional collection in $D(\M^{\K})$, then this would show that $D(\M^{\K})\not\cong D(\overline{\B^9/\Gamma}^{\tor})$.
    \item For a projective variety or Deligne-Mumford stack $X$, if $D(X)$ has a full exceptional collection, then $D(X)$ admits a non-trivial semi-orthogonal decomposition.
    Hence, in order to show that $D(\M^{\K})\not\cong D(\overline{\B^9/\Gamma}^{\tor})$, it suffices to prove that $D(\M^{\K})$ admits no such decomposition.
    There are sufficient conditions that might be helpful \cite[Theorems 3.1, 3.3]{KO15}.
    However, these criteria are currently only known for smooth schemes or stacks, and hence cannot immediately be used for $D(\M^{\K})$ but $D(\mathcal{X}(\M^{\K}))$ 
    and $D(\mathcal{X}(\overline{\B^9/\Gamma}^{\tor}))$.
      \end{enumerate}
    \end{rem}
   Here we summarize some remaining open problems concerning the derived categories in our cases.
    \begin{que}
    \label{que:derived_categories}
    The above discussion leads to the following questions: 
        \begin{enumerate}
            \item Does $D(\M^{\K})$ have a full exceptional collection?
            \item Determine whether  the non-stacky  categories $D(\M^{\K})$ and $D(\overline{\B^9/\Gamma}^{\tor})$ are equivalent or not? 
            \item Can one give another proof of Theorems  \ref{thm:derived_equivalent} and \ref{thm:derived_cubic} (1) in terms of semi-orthogonal decompositions?
            \item Does a full exceptional collection in $D(\overline{\B^9/\Gamma}^{\tor})$ induce one in $D(\mathcal{X}(\overline{\B^9/\Gamma}^{\tor}))$?
        \end{enumerate}
    \end{que}

\section{Computation of the cohomology}
\label{sec:cohomology}
We shall now determine the cohomology groups of the spaces appearing in this paper.

\subsection{The cohomology of the spaces except  $\M^{\K}$ and  $\overline{\B^9/\Gamma}^{\tor}$}
Due to the work of Kirwan \cite{Kir89} and Kirwan-Lee-Weintraub \cite{KLW87}, the topology of all spaces except $\M^{\K}$ and  $\overline{\B^9/\Gamma}^{\tor}$ is known:
\begin{thm}[{Kirwan-Lee-Weintraub \cite[Table III, Theorem 8.6]{KLW87}}, Kirwan {\cite[Table, p.40]{Kir89}}]
\label{thm:coh_previous_work}
All the odd degree cohomology of $\M_{\ord}^{\K}$,  $\M_{\ord}^{\GIT}$, $\M^{\GIT}$ and $\overline{\B^9/\Gamma}^{\BB}$ vanishes.
In even degrees, the Betti numbers are as follows:

\begin{align*}
\renewcommand*{\arraystretch}{1.2}
\begin{array}{l|cccccccccc}
\hskip2cm j&0&2&4&6&8&10&12&14&16&18\\\hline
\dim H^j(\M_{\ord}^{\K})&1&474&991&1618&2410&2410&1618&991&474&1\\
\dim IH^j(\M_{\ord}^{\GIT})&1&12&67&232&562&562&232&67&12&1\\
\dim H^j(\M^{\GIT})&1&1&2&2&3&3&2&2&1&1\\
\dim IH^j(\overline{\B^5/\Gamma}^{\BB})&1&1&2&2&3&3&2&2&1&1
\end{array}
\end{align*}
\end{thm}

\subsection{The cohomology of the Kirwan blow-up $\M^{\K}$}
The seminal work of Kirwan \cite{Kir84, Kir85, Kir89} is the fundamental tool to compute the cohomology of GIT quotients and their Kirwan blow-up.
Here we follow the set-up of \cite{CMGHL23a}, which is based on Kirwan's work.  
For $X\defeq\P^{12}$, acted on by $G=\SL_2(\C)$ with the symmetric linearization induced via $\Sym^{12}(\P^)\cong \P^{12}$, let $\widetilde{X}^{\ss}$ be the blow-up of the semi-stable locus $X^{\ss}$ whose center 
is the unique polystable orbit $G\cdot Z_R^{\ss}$.
Here, $Z_R^{\ss}$ is the fixed locus of the action of the stabilizer $R$, computed in Lemma \ref{lem:Luna_slice} (2), on $X^{\ss}$.
The main idea for computing the cohomology of the Kirwan blow-up $\M^{\K}$ is to reduce the problem to a calculation on  $\widetilde{X}^{\ss}$.
In our case, $\M^{\K}$ is a blow-up of $\M^{\GIT}$ at the unique cusp, see \cite[Lemma 3.11]{Kir85}.
The space $\M^{\K}$ is defined as the GIT quotient $\widetilde{X}^{\ss}$ by $G$ (with respect to a suitable linearization).
The Poincare series of $\M^{\K}$ is given by the $G$-equivariant cohomology of $\widetilde{X}^{\ss}$:
\[P_t(\M^{\K})=P_t^G(\widetilde{X}^{\ss}).\]

Moreover, from \cite[Section 3 Eq. 3.2]{Kir89} or \cite[Subsection 4.12, (4.22)]{CMGHL23b}, 
the Poincare series of $\widetilde{X}^{\ss}$ can be computed from that of $X^{\ss}$ as
\[P_t^G(\widetilde{X}^{\ss})=P_t^G(X^{\ss})+A_R(t),\]
where $A_R(t)$ is a correction term. This consists of a ``main term" and an  ``extra term" with respect to the unique stabilizer $R$; see \cite[Section 4.1.2]{CMGHL23b} for precise definitions.
We shall outline how to determine these terms, for full details see \cite[Chapter 3, 4]{CMGHL23a}.  
Let $\mathcal{B}$ be the set of points that are closest to the origin of the convex hull spanned by some weights in the closure of a positive Weyl chamber in the Lie algebra of a maximal torus in $\SO(2)$
(as introduced in \cite[Definition 3.13]{Kir84}).
Let $\{S_{\beta}\}_{\beta\in\mathcal{B}}$ be the stratification defined in \cite[Theorem 4.16]{Kir84} and $d(\beta)\defeq\mathrm{codim}_{X^{\ss}}S_{\beta}$.
Let $\mathcal{N}$ be the normal bundle to the orbit $G\cdot Z_R^{\ss}$.
Then, for a generic point $x\in Z_R^{\ss}$, we have a representation $\rho$ of $R$ on $\mathcal{N}_x$.
Let $\mathcal{B}(\rho)$ be the set consisting of the closest point to 0 of the convex hull of a non-empty set of weights of the representation $\rho$.
For $\beta'\in\mathcal{B}(\rho)$, let $n(\beta')$ be the number of weights less than $\beta'$.

The following proposition allows us to determine $P_t^G(\widetilde{X}^{\ss})$.
\begin{prop}
The following holds:
\begin{enumerate}
    \item The equivariant cohomology of the semi-stable locus is given by  \[P_t^G(X^{\ss})\equiv 1+t^2+2t^4+2t^6+3t^8 \bmod t^{10}.\]
    \item The main correction term in $A_R(t)$ is \[ (1-t^4)^{-1}(t^2+t^4+t^6+t^8)\equiv t^2+t^4 + 2t^6+2t^8 \bmod t^{10}.  \]
    \item The extra correction term vanishes modulo $t^{10}$, i.e.,  does not contribute to $A_R(t)$.
\end{enumerate}

\end{prop}
\begin{proof}
(1) We decompose the equivariant cohomology of the semistable locus as 
\[P_t^G(X^{\ss})=P_t(X)P_t(B\SL_2(\C))-\sum_{0\neq\beta\in\mathcal{B}}t^{2d(\beta)}P_t^G(S_{\beta})\]
by \cite[Eq. 1]{Kir89}.
We can determine the set $\mathcal{B}$ of the weights and  show $2d(\beta)\geq 10$ for $\beta\neq 0$ using the same method as in the proof of \cite[Proposition 3.5]{CMGHL23a} or \cite[Proposition 5.3]{HM25}.
Thus it follows that
    \begin{align*}
    P_t^G(X^{\ss})&\equiv P_t(X)P_t(B\SL_2(\C))\bmod t^{10} \\
        &\equiv (1-t^2)^{-1}(1-t^4)^{-1}\bmod t^{10}\\
    &\equiv 1+t^2+2t^4+2t^6+3t^8 \bmod t^{10}.
\end{align*}

(2) Let $N(R)$ be the normalizer of $R$.
Then, the main correction term (modulo $t^{10}$) is given by 
\[P_t^{N(R)}(Z_R^{\ss})(t^2+t^4+t^6+t^8) \bmod t^{10}\]
by \cite[(4.24)]{CMGHL23a}.
Next we claim that $P_t^{N(R)}(Z_R^{\ss})=(1-t^4)^{-1}$ (compare to \cite[Proposition 5.4]{HM25}).
Indeed, since the normalizer of $R$ is $N(R)\cong \mathbb{T}\rtimes \Z/2\Z$,
 it follows that
\begin{align*}
    H_{N(R)}^{\bullet}(Z_R^{\ss})&=(H_{\mathbb{T}}^{\bullet}(Z_R^{\ss}))^{\Z/2\Z}\\
    &=(H^{\bullet}(BR)\otimes H^{\bullet}_{\mathbb{T}/R}(Z_R^{\ss}))^{\Z/2\Z}\\
    &=(H^{\bullet}(BR)\otimes H^{\bullet}(*))^{\Z/2\Z}\\
        &=\Q[c^4]
\end{align*}
where $*$ denotes a set of 1 point and the degree of $c$ is 1.
The last equation follows from \cite[Proposition 4.4]{CMGHL23a}.
This implies our claim.

(3) By \cite[(4.25)]{CMGHL23a}, the extra correction term is a multiple of $t^{2d(|\beta'|)}$ for any $0\neq\beta'\in\mathcal{B}(\rho)$.
Note that $d(|\beta'|)=n(|\beta'|)$, where $n(|\beta'|)$ is the number of weights less than $|\beta'|$; see \cite[Proof of Proposition 4.19]{CMGHL23a}.
Hence, it suffices to determine the weight set $\mathcal{B}(\rho)$, which follows from the description of the normal bundle $\N_x$ for $x\in Z_R^{\ss}$.
However, in our case, since $Z_R^{\ss}=\{p_{6,6}\}$, we only need to compute the weights for an element in $R$ acting on $T_{p_{6,6}}\C^{13}$.
Lemma \ref{lem:Luna_slice} (2) implies that the weights are
\[0,\pm2, \pm4, \pm6, \pm8, \pm10, \pm12\]
and $T_{p_{6,6}}(G\cdot\{p_{6,6}\})$ is generated by the weights  $\{0,\pm2\}$.
Hence, it follows that
\[\B(\rho)=\{\pm2,\pm4,\pm6,\pm8,\pm10\}\]
 and thus we have
 \[n(|\beta'|)\geq 5.\]

\end{proof}

It follows that
\begin{align*}
    P_t(\M^{\K})&= P_t^G(\widetilde{X}^{\ss})\\
    &\equiv (1+t^2+2t^4 + 2t^6 + 3t^8)+(t^2+t^4 + 2t^6+2t^8) \bmod t^{10}\\
    &\equiv 1+2t^2+3t^4 + 4t^6+5t^8 \bmod t^{10},
\end{align*}
 and hence we obtain 
\begin{thm}
\label{thm:coh_M^K}
All the odd degree cohomology of $\M^{\K}$ vanishes.
In even degrees, its Betti numbers are given as follows:
\begin{align*}
\renewcommand*{\arraystretch}{1.2}
\begin{array}{l|cccccccccc}
\hskip2cm j&0&2&4&6&8&10&12&14&16&18\\\hline
\dim H^j(\M^{\K})&1&2&3&4&5&5&4&3&2&1\\
\end{array}
\end{align*}
\end{thm}

\subsection{The cohomology of the toroidal compactification $\overline{\B^9/\Gamma}^{\tor}$}
Next, we compute the cohomology of the toroidal compactification $\overline{\B^9/\Gamma}^{\tor}$.
Recall that $T$ is the toroidal boundary of $\overline{\B^9/\Gamma}^{\tor}$, which is isomorphic to $\P^4\times\P^4/(S_6\times S_6)\rtimes S_2$; see \cite[Remark 4.6]{Has03}.

\begin{prop}
\label{prop:coh_boundary}
All the odd degree cohomology of the boundary $T$ vanishes.
In even degrees, its Betti numbers are given as follows:
\begin{align*}
\renewcommand*{\arraystretch}{1.2}
\begin{array}{l|ccccccccc}
\hskip2cm j&0&2&4&6&8&10&12&14&16\\\hline
\dim H^j(T)&1&1&2&2&3&2&2&1&1\\
\end{array}
\end{align*}
\end{prop}
\begin{proof}
We have to compute the cohomology ring
\[H^{\bullet}(T)=
 H^{\bullet}(\P^4 \times \P^4)^{(S_6\times S_6)\rtimes S_2} = H^{\bullet}((\P^4/S_6)^2, \Q)^{S_2}.
\]
A similar calculation can be found in   \cite[Proposition 7.13]{CMGHL23a}
Since $H^{\bullet}(\P^4/S_6) \cong \Q[x]/(x^5)$, 
   it sufficies to determine $S_2$-invariant part of the tensor product $\Q[x]/(x^5)\otimes\Q[y]/(y^5)$.
   The invariant part spanned by the invariant polynomials 1 in degree 0, $x+y$ in degree 1, $x^2+y^2, xy$ in degree 2, $x^3+y^3, x^2y+xy^2$ in degree 3, $x^4+y^4, x^3y+xy^3, x^2y^2$ in degree 4.
    Hence, the invariant cohomology is given by 
    \[P_t(T)=1+t^2+2t^4+2t^6+3t^8+2t^{10}+2t^{12}+t^{14}+t^{16}.\]
\end{proof}

Combining this with the decomposition formula \cite[Lemma 9.1]{GH17}, we finally obtain  
\begin{thm}
\label{thm:coh_tor}
All the odd degree cohomology of  $\overline{\B^9/\Gamma}^{\tor}$ vanishes.
In even degrees, its Betti numbers are given as follows:
\begin{align*}
\renewcommand*{\arraystretch}{1.2}
\begin{array}{l|cccccccccc}
\hskip2cm j&0&2&4&6&8&10&12&14&16&18\\\hline
\dim H^j(\overline{\B^9/\Gamma}^{\tor})&1&2&3&4&5&5&4&3&2&1\\
\end{array}
\end{align*}
\end{thm}

\section{Automorphic forms}
\label{section:Automorphic forms}
So far, we have dealt with the canonical bundles of modular varieties mostly from a geometric point of view.
At the same time, pluricanonical forms on modular varieties are closely related to automorphic forms. In our case, we have seen   
this in Remark \ref{rem:compnormalbundle} and the paragraph before it.

In this section, we will consider the aspect of automorphic forms.
First, we shall give a new construction of an
automorphic form on another, but closely related 9-dimensional ball quotient.
For this, we will recall in  Subsection \ref{subsec:A complex ball uniformization} the complex ball uniformization of the moduli space $\M_{4}^{\mathrm{nh}}$ of non-hyperelliptic curves of genus 4 given in \cite{Kon02}, and then construct an automorphic form whose zero divisor has an interesting geometric interpretation.
Second, we review the automorphic form $\Psi_1=\Psi_{\mathrm{A}}^3$ on $\B^9/\Gamma$, where $\Psi_{\mathrm{A}}$ is the form constructed by Allcock  \cite[Theorem 7.1]{All00}, 
see also Remark \ref{rem:compnormalbundle}.
We give a new construction of this form as a quasi-pullback of Borcherd's automorphic form $\Phi_{12}$. This allows us to compare the automorphic forms we and Allcock construct on $\B^9/\Gamma_{\mathrm{nh}}$ and $\B^9/\Gamma$  respectively; see Remark \ref{rem:prob_modular_forms}.
This result is independent of, and not strictly necessary, for our previous discussion, but we believe that it contributes to our overall understanding of the situation.

Below,  for any lattice $M$ of signature $(2,n)$ let us denote by 
\[\mathbb{D}_{M}\defeq\{[w]\in\P(M\otimes\C)\mid \l w,w\r=0, \l w,\overline{w}\r>0\}^+\]
the associated hermitian domain,
where $\l \, , \,  \r$ denotes the hermitian form and the superscript + indicates that one chooses one of the two connected components.  
The manifold $\mathbb{D}_{M}$ is the $n$-dimensional Hermitian symmetric domain of type IV associated with $M$.
For positive integers $m,n$ with $m-n \equiv 0\ ({\rm mod}\ 8)$, 
we denote by $II_{m,n}$ the even unimodular lattice of signature $(m,n)$.  For an even lattice $M$, let  $A_M=M^*/M$ be the {\em discriminant group} and 
\begin{alignat*}{2}
    q_M&: A_M \to \Q/2\Z,\quad &&q_M(x + M) \defeq x^2 \mod 2\\  b_M&: A_M\times A_M \to \Q/\Z,\quad &&b_M(x+M, y+M) \defeq \langle x, y\rangle \mod \Z.
\end{alignat*}
We call $q_M$ (resp. $b_M$) the \textit{discriminant quadratic form} (resp. the \textit{discriminant bilinear form}) of $M$.  

An even lattice $M'$, containing $M$, is called an {\it overlattice} of $M$ if $M$ is of finite index in $M'$.
Let $H \subset A_M$ be an isotropic subgroup with respect to $q_M$, that is, $q_M|H=0$.  Then, 
\[M_H\defeq\{ x \in M^* \mid x + M \in H\}\]
is an even lattice with $A_{M_H}= (M_H^*/M)/H$ and $q_{M_H} = (q_M|H^{\perp})/H$.  Conversely, for any overlattice $M'$ of $M$, $M'/M$ is an isotropic subgroup.  Thus isotropic subgroups in $A_M$ 
correspond bijectively to overlattices of $M$.

Let $M$ be a {\it primitive} sublattice of an even lattice $L$, that is, 
$L/M$ is torsion free.  Let $N$ be the orthogonal complement of $M$ in $L$ which is also primitive in $L$.  
Then $L$ is an overlattice of $M\oplus N$ and hence the subgroup $L/(M\oplus N) \subset A_M\oplus A_N$ is isotropic with respect to $q_M\oplus q_N$.
The primitiveness implies that the two projections $L/(M\oplus N)\to A_M$ and $L/(M\oplus N)\to A_N$ are injective.  In particular, if $L$ is unimodular, then $|L/(M\oplus N)|=|A_M|=|A_N|$, 
and hence these projections are isomorphisms.  Thus, we obtain a canonical isomorphism
$\ : A_M \to A_N$ with $q_N \circ\varphi = - q_M$.  For more details, we refer the reader to Nikulin \cite[Section 1]{Nik80}.

\subsection{Moduli of curves of genus 4}
\label{subsec:A complex ball uniformization}

Let $C$ be a non-hyperelliptic curve of genus 4.  Then the canonical model of $C$ is the
intersection of a non-singular quadric surface or a quadric cone $Q$ and a cubic surface $S$ in 
$\P^3$: $C= Q\cap S$.
It is known that $C$ has two (resp. one) $g_3^1$ (= a pencil of degree 3) if
$Q$ is non-singular (resp. a quadric cone), induced from the ruling(s) of $Q$.  
The canonical class $K_C$ is the sum of  the two $g_3^1$ in case $Q$ is non-singular, and 
$K_C = 2g_3^1$ in case $Q$ is a cone. 
In the latter case, we say that the curve $C$ has a {\em vanishing theta 
constant}.  In case $C$ is hyperelliptic
any $g_3^1$ has a base point $q$ and $g_3^1= g_2^1 + q$ where $g_2^1$ gives the hyperelliptic involution.

For a smooth quadric $Q$, 
there exists a triple covering of
$Q$ branched along $C$ whose minimal resolution is a $K3$ surface $X$.  
Let $\sigma$ be the covering transformation of $X\to Q$.
Since $Q$ is rational, $\sigma$ is a non-symplectic automorphism of $X$ of order 3, that is, 
$\sigma^*(w_X)= \omega w_X$ where $w_X$ is a nowhere vanishing holomorphic 2-form on $X$
and $\omega$ is a primitive cube root of unity.  Note that $\Pic(Q)\cong U$ is generated by the classes of
fibers of two rulings.  
The two rulings induce two elliptic fibrations on $X$ whose general fibers $F_1, F_2$ satisfy
$F_1\cdot F_2=3$.  Thus,
$\Pic(X)$ contains the primitive sublattice $U(3)$ generated by the classes of $F_1, F_2$. The order 3 automorphism $\sigma$ acts trivially on $U(3)$. 

Recall the $K3$ lattice 
\[L_{K3}\defeq II_{3,19},\]
which is (abstractly) isomorphic to $H^2(X,\Z)$ equipped with the cup product.  
Since $\Pic(X)$ is primitive in $L_{K3}$, there exists a primitive embedding of $U(3)$ into $L_{K3}$.
It follows from Nikulin's analog of the Witt theorem \cite[Theorem 1.14.4]{Nik80}, that a primitive embedding of 
$U(3)$ into $L_{K3}$ is unique up to isomorphisms, and hence the orthogonal complement is also uniquely determined by its genus.  Now we define $L$ as the orthogonal complement of $U(3)$ in $L_{K3}$.
Since the three lattices $L$, $U\oplus U(3)\oplus E_8^{\oplus 2}$ and $U^{\oplus 2}\oplus E_8\oplus E_6\oplus A_2$ have the same discriminant form isomorphic to $q_{U(3)}= -q_{U(3)}$, it follows that 
$L_{K3}$ is an overlattice of the orthogonal direct sum of $U(3)$ and any one of these three lattices.  In other words, the three lattices
are the orthogonal complement of $U(3)$ in $L_{K3}$.  Thus
we have
\[L\cong U\oplus U(3)\oplus E_8^{\oplus 2} \cong U^{\oplus 2}\oplus E_8\oplus E_6\oplus A_2.\] 
 
Let $e, f$ (resp. $e', f'$) be a basis of $U(3)$ (resp. $U$) with $e^2=f^2=0, \langle e, f\rangle =3$ (resp.
$e'^2=f'^2=0, \langle e', f'\rangle =1$).  
Let $\theta_1$ be the isometry of $U(3)\oplus U$ of order 3 defined by
\[\theta_1(e)\defeq -2e+3e',\ \theta_1(f)\defeq f + 3f',\ \theta_1(e')\defeq -e+e',\ \theta_1(f')\defeq -f -2f'.\]
Note that $\theta_1$ has no non-zero fixed vectors and that the action of the discriminant group of 
$U(3)\oplus U$ is trivial.  
Similarly, $A_2$ has an isometry of order 3 without non-zero fixed vectors which acts trivially on $A_{A_2}$.
In fact, if we denote by $e_1, e_2$ a basis of $A_2$ with $e_1^2=e_2^2=-2$, $\langle e_1, e_2\rangle =1$,
then the isometry $g$ defined by $g(e_1)=e_2, g(e_2)= -e_1-e_2$ is such an isometry.  By taking an isotropic 
subgroup of $A_{A_2^{\oplus 4}}$ of order $3^2$,
we obtain an overlattice which is an even unimodular negative definite lattice of rank 8, that is, $E_8$.  Since the isometry $(g,g,g,g)$ of $A_2^{\oplus 4}$ preserves the isotropic subgroup, 
it can be extended to an isometry of the overlattice.  Thus we have an isometry $\theta_2$ 
of $E_8$ of order 3 without non-zero fixed vectors. 
In this way, we can define an isometry $\theta$ of $L$ by setting
$\theta := (\theta_1, \theta_2, \theta_2)$. This has no non-zero fixed vectors in $L$
and acts trivially on the discriminant group $A_L$ of $L$.  A similar argument shows that 
$\theta$ can be extended to an isometry of $L_{K3}$, which we will denote by the same symbol $\theta$, satisfying $\theta\vert_{U(3)} = 1$. 
Then there exists an isomorphism $\varphi: H^2(X,\Z) \to L_{K3}$
satisfying $\varphi \circ \sigma^*\circ \varphi^{-1}= \theta$ or $\theta^2$; see \cite[p. 386]{Kon02}
for a proof.
We now consider the eigenspace decomposition 
$L\otimes \C = L_{\omega}\oplus L_{\omega^2}$
of $\theta$ and the subdomain of $\mathbb{D}_L$ defined by
\[\B_{L_{\omega}} \defeq \{[v] \in \P(L_{\omega}) \cap \mathbb{D}_L\mid\langle v, \bar{v}\rangle > 0\} \subset \mathbb{D}_L.\]
Note that $\langle v, v\rangle =0$ for any $v \in L_{\omega}$ because 
$\langle v, v\rangle = \langle \theta(v), 
\theta(v)\rangle = \omega^2 \langle v, v\rangle$.
This is equivalent to taking a Hermitian lattice $\Lambda_{\mathrm{nh}}$ over $\Z[\omega]$ of signature $(1,9)$ such that the composition of the Hermitian form and the trace map $\Q(\omega)\to\Q$ is equal to $L$.
Since the Hermitian form $\langle v, \bar{v}\rangle$ has signature $(1,9)$, the space $\B_{L_{\omega}}$ is nothing but the complex ball $\B^9$ of dimension 9, which we also considered earlier.
Let 
\begin{align}
\label{def:kondo_unitary}
    \Gamma_{\mathrm{nh}}&\defeq \{ g \in {\rm O}(L)^+ \mid g\circ \theta= \theta\circ g\}\\
    &=\U(\Lambda_{\mathrm{nh}}).\notag
    \end{align}
The Torelli  theorem for $K3$ surfaces then implies that the above correspondence gives an injection
\[\M_4^{\mathrm{nh}} \to \B^9/\Gamma_{\mathrm{nh}}.\]

Next, we consider the discriminant locus, that is, the complement of $\M_4^{\mathrm{nh}}$ in $\B^9/\Gamma_{\mathrm{nh}}$.
It is easy to see that $A_L \cong (\Z/3\Z)^2$ consists of the following 9 vectors:
\smallskip
\begin{align*}
    \mathrm{type}\ (00)&:\  q_{L}(\alpha) = 0,\ \alpha = 0;\ \# \alpha = 1,\\
    \smallskip
    \mathrm{type}\ (0)&:\ q_{L}(\alpha) = 0,\ \alpha\ne 0;\ \# \alpha = 4,\\
    \smallskip	
    \mathrm{type}\ (4/3)&:\ q_{L}(\alpha) = -\frac{2}{3};\ \# \alpha = 2,\\
    \smallskip
    \mathrm{type}\ (2/3)&:\ q_{L}(\alpha) = -\frac{4}{3};\ \# \alpha = 2.
\end{align*}
\smallskip
For each $r\in L\otimes \Q$ with $r^2<0$ and for $\alpha \in A_L, n \in \Q$, we define \textit{Heegner divisors (of discriminant $n$)} on $\mathbb{D}_L$ by
\[\D_r \defeq   r^\perp \defeq\{ [\omega] \in \D_L \mid \langle \omega, r\rangle = 0\},\quad 
\D_{\alpha,n} \defeq\bigcup_{\substack{r+L=\alpha\\ r^2=n}} \D_r;\]
\[\D_{-2} \defeq \D_{0,-2}, \quad \D_{-2/3} \defeq \D_{\alpha,-2/3}, \quad \D_{-4/3} \defeq \D_{\alpha,-4/3}.\]
Note that $\D_{\alpha,n} = \D_{-\alpha,n}$ and there are only two vectors of type $(4/3)$ and of type $(2/3)$, 
which differ by a sign, and hence $\D_{-2/3}$, $\D_{-4/3}$ are independent of the choice of $\alpha$.
We recall that we also used the notation $\D$ for the discriminant locus in the Kirwan blow-up $\M^{\K}$, but as these spaces will not play a role in what follows, we trust that this will not 
cause any problems.   
Similarly, we obtain and define
\begin{align}
\label{def:divisors_kondo}
    \H = \B^9\cap \D_{-2},\quad \H_{\mathrm{h}} \defeq \B^9\cap \D_{-2/3},\quad \H_{\mathrm{vt}} \defeq \B^9\cap \D_{-4/3}.
\end{align}
We consider $\D_{\alpha,n}$ as a divisor on $\D_L$ by attaching multiplicity 1 to all components $\D_r$ and also call it a Heegner divisor. 
By \cite[Lemmas 3, 4]{Hof14} the intersection $\H_*$ of a Heegner divisor $\D_*$ and the ball $\B^9$ is also a Heegner divisor. 
It is known that $\H/\Gamma_{\mathrm{nh}}$ is the complement of $\M_4^{\mathrm{nh}}$ in $\B^9/\Gamma_{\mathrm{nh}}$, and hence
we call $\H$ the \textit{discriminant locus}.  We shall see below that $\H$ decomposes into two components, one of which is $\H_{\mathrm{h}}$.

In the following, we recall some details concerning the discriminant locus.
For $r\in L$ with $r^2=-2$, let $\Lambda_r$ be the lattice generated by $r$ and $\theta(r)$.  
It follows from the equation $\theta^2 + \theta + 1_L=0$ that
$\langle r, \theta(r)\rangle =1$, and hence $\Lambda_r$ is isomorphic to a root lattice of type $A_2$.
Let $\Lambda_r^{\perp}$ be the orthogonal complement of $\Lambda_r$ in $L$ and 
$M$ the orthogonal complement of $\Lambda_r^{\perp}$ in $L_{K3}$.  
Recall that the first projection
\[L/(\Lambda_r\oplus \Lambda_r^{\perp}) \to A_{\Lambda_r} (\cong A_{A_2} \cong \Z/3\Z)\]
is injective.  Thus we have 
\[ [L : \Lambda_r\oplus \Lambda_r^{\perp}]=3 \ \ {\rm or}\ \ L=\Lambda_r\oplus \Lambda_r^{\perp}.\]
For example, if $\Lambda_r$ is a primitive sublattice of $E_8$ in a decomposition 
$L= U\oplus U(3) \oplus E_8^{\oplus 2}$, then the first case occurs.  If $\Lambda_r$ is the last component of
the decomposition $L=U^{\oplus 2}\oplus E_8\oplus E_6\oplus A_2$, then the second case occurs.
Altogether, we get the following two cases:

\smallskip

(i)\ $M\cong U(3)\oplus A_2$ and $\Lambda_r^{\perp}\cong U(3)\oplus U \oplus E_8\oplus E_6$;
\smallskip

(ii)\ $M\cong U\oplus A_2$ and $\Lambda_r^{\perp}\cong U\oplus U\oplus E_8\oplus E_6$.

\smallskip
\noindent
Accordingly, the divisor $\H$ decomposes into two components $\H = \H_{\mathrm{n}} + \H_{\mathrm{h}}$.   

Recall that the sublattice $U(3)\subset \Pic(X)$ is generated by two elliptic curves.  In other words, there are two elliptic fibrations on $X$.  
Now assume $\Pic(X)=U(3)$.
If these fibrations have a reducible fiber, the Picard number is at least three, which we have just excluded. Hence all fibers must be irreducible in this case.  
Since the covering transformation $\sigma$ over $Q$ preserves the fibers, a non-singular fiber
is an elliptic curve with an automorphism of order 3.  This implies that the fibration is isotrivial and
the irreducible singular fiber is of type ${\rm II}$ in the sense of Kodaira.  
Since the Euler number of a $K3$ surface is 24 and that of a singular fiber of type ${\rm II}$ is 2, 
there are exactly twelve singular fibers of type ${\rm II}$ on each fibration. We also note that this property is an open condition as we vary the $K3$ surface in its 9-dimensional moduli.

We shall now consider the special cases, starting with (i). Hence we consider $K3$ surfaces $X$ whose Picard group $\Pic(X)$ contains $U(3)\oplus A_2$ as a primitive 
sublattice and which further admit a non-sympletic automorphism $\sigma$ acting trivially on the sublattice $U(3)\oplus A_2$. This corresponds to the period domain $\H_{\mathrm{n}}$. We further claim that this 
corresponds to a nodal curve $C$ of arithmetic genus 4.
Indeed, the triple covering of $Q$ branched along $C$ has a rational double point of type $A_2$ over the node. 
The corresponding smooth $K3$ surface $X$ is its minimal resolution.
It is obtained as follows.  We first blow up the node of $Q$ and denote it by $E_0$.  
Next, we blow up the two points which are
the intersection of $E_0$ and the proper transform of $C$.  Finally, we take the triple covering of $Q$ 
branched along the proper transform of $C$ and that of $E_0$.
The pre-image of $E_0$ has self-intersection number $-1$ and hence we can contract it to a point $q_0$,  
giving us the desired surface $X$.  Let $\ell_1, \ell_2$ be
two lines on $Q$ passing through the node.  The pre-images $L_i$ of $\ell_i$ $(i=1,2)$ are the triple 
coverings of $\ell_i$ branched at two points, namely the intersection points with the proper transforms of $E_0$ and $C$, and are $(-2)$-curves on $X$. 
Let $E_1, E_2$ be the preimages of the exceptional curves of the second blow-ups, which are also 
$(-2)$-curves on $X$.  
The four $(-2)$-curves meet exactly at one point $q_0$ transversally 
and generate a sublattice $M= U(3)\oplus A_2$, where
$U(3)$ is generated by the classes of $L_1+E_1+E_2$, $L_2+E_1+E_2$ and $A_2$ is generated by $E_1, E_2$.
The linear systems $|L_1+E_1+E_2|$ and
$|L_2+E_1+E_2|$ define two elliptic fibrations with one reducible singular fiber of type IV.  
The fixed sublattice of the covering transformation is $M=U(3)\oplus A_2$.
In this case, $\Lambda_r$ is the component $A_2$ of $M$.  

We now turn to case (ii), where we again want to establish a relationship with $K3$ surfaces. For this, we consider $K3$ surfaces $X$ whose Picard group $\Pic(X)$ contains $U\oplus A_2$ as a primitive 
sublattice and which further admit a non-sympletic automorphism $\sigma$ acting trivially on the sublattice $U\oplus A_2$. 
The period domain of these surfaces is $\H_{\mathrm{h}}$.
In the case of non-hyperelliptic curves of genus $4$, we used its canonical model to obtain a correspondence between genus 4 curves and $K3$ surfaces by taking a triple covering of a rational surface. 
This is not available to us in the hyperelliptic case, but we shall give another geometric construction, which leads to 
pairs $(C, (q_+,q_-))$ with $C$ a hyperelliptic curve of genus 4 and two points $q_+, q_-$ on $C$ 
conjugated under the hyperelliptic involution.
The assumption on the Picard lattice implies that  $X$ has an elliptic fibration $\pi : X \to \P^1$ with a reducible singular fiber of type 
${\rm I}_3$ or type ${\rm IV}$ and a section (the classes of a fiber and the section generate the sublattice $U$ and the irreducible components of reducible fibers not meeting the section generate a sublattice $A_2$).  
Thus, there exist four $(-2)$-curves $E_0, E_1, E_2, E_3$ on $X$ such that $E_1+E_2+E_3$ is 
the reducible fiber and $E_0$ is the section.  We assume that $E_0$ meets $E_1$.   
By assumption, the order 3 automorphism $\sigma$ preserves the elliptic fibration. 

It follows from the topological Lefschetz fixed point formula that the Euler number of the set of fixed 
points of $\sigma$ is $2 + 4 + 9(\omega +\omega^2) = -3$.  If $\sigma$ were to act non-trivially on the base of the fibration,
then the fixed point would be contained in the two invariant fibers, implying that its Euler number is non-negative.  Thus
$\sigma$ acts trivially on the base and hence acts on each fiber as an automorphism.
This implies that the reducible fiber is of type ${\rm IV}$ and all irreducible singular fibers are of type
${\rm II}$.  Thus, the fact that the Euler number $e(X)=24$, shows that $\pi$ has one singular fiber of type ${\rm IV}$ (contributing Euler number 4) and ten singular fibers of type ${\rm II}$ (each contributing Euler number 2).
By using the topological Lefschetz fixed point formula, one can further see that 
the fixed point set of $\sigma$ is the disjoint union of a non-singular curve $C$ of genus 4, $E_0$ and 
the singular point $p_0$ of the fiber of type ${\rm IV}$. The restriction of $\pi$ to $C$ gives a double covering $C\to \P^1$ and hence $C$ is hyperelliptic.  The ramification locus of $\pi|C$ consists of 
10 singular points of singular fibers of type ${\rm II}$, and, in addition, $C$ meets the fiber of type ${\rm IV}$ at
a point $q_+$ on $E_2$ and $q_-$ on $E_3$.  In this way, we obtain a pair $(C, q=(q_+,q_-))$.
To realize $X$ as the minimal model of a cyclic covering of degree 3 over a rational surface,  
we first blow up $\widetilde{X}\to X$ at $p_0$ and denote the exceptional curve by $E$.  Then
$\sigma$ induces an automorphism $\widetilde{\sigma}$ of $\widetilde{X}$ whose fixed point set consists of
$E$ and the proper transforms  $\widetilde{E}_0$, $\widetilde{C}$ of $E_0$, $C$, respectively.  
The images of $E_1, E_2, E_3$ on $\widetilde{X}/\langle \widetilde{\sigma} \rangle$ are $(-1)$-curves and we can 
blow down them.  In this way, we obtain a Hirzebruch surface $H_5$ on which the $(-5)$-curve is the image of $\widetilde{E}_0$.  For more details, see  \cite[Example 2]{Kon02}.

Finally, there is another important divisor, namely the one parametrizing non-singular, non-hyperelliptic curves $C$ of genus 4 
with a vanishing theta constant. This is also sometimes called the theta null locus. The reason is that such a curve has a vanishing theta constant. Alternatively, these curves can be characterized by having an effective 
theta characteristic. In this case the two $g^1_3$'s coincide, $K_C=2g^1_3$ and      
the canonical model of $C$ is the intersection of a quadric cone $Q_0$ 
and a cubic surface $S$ in $\P^3$.  
The minimal resolution of the triple covering of $Q_0$
branched along $C$ is a $K3$ surface $X$ such that $\Pic(X)$ primitively contains a lattice isomorphic to 
$U\oplus A_2(2)$.  The $K3$ surface $X$ contains three disjoint $(-2)$-curves $E_1, E_2, E_3$ 
which are exceptional curves over the vertex.
The pencil of lines passing through the node of $Q_0$ induces an elliptic fibration on $X$.  In the generic case, the fibration has twelve singular fibers of type II and three sections $E_1, E_2, E_3$. 
Here $U$ is generated by the classes of a fiber $F$ and $E_1$, and $A_2(2)$ is generated by $E_2-E_3$ and $2F+E_1-E_2$.
The invariant lattice under the action of $\sigma^*$ is generated by $F$ and $E_1+E_2+E_3$ and isomorphic to
$U(3)$. The orthogonal complement of $U(3)$ in $\Pic(X)$ is isomorphic to $A_2(2)$ and is generated by
$E_1-E_2, E_3-E_1$.  Note that the discriminant group $A_{U(3)}$ of $U(3)$ is isomorphic to 
$(\Z/3\Z)^2$ generated by $(E_1 +E_2+E_3)/3, F/3$, and $A_{A_2(2)}$ is isomorphic to 
$(\Z/2\Z)^2\oplus \Z/3\Z$ generated by $(E_1-E_2)/2, (E_3-E_1)/2, (2E_1 -E_2-E_3)/3$.  Further note that 
$q_{U(3)}((E_1+E_2+E_3)/3)= -2/3$, $q_{A_2(2)}((2E_1 -E_2-E_3)/3)= -4/3$, and that the class
$E_1\mod U(3)\oplus A_2(2)$ is a generator of $U\oplus A_2(2)/U(3)\oplus A_2(2)$ using that 
$E_1= (E_1+E_2+E_3)/3+(2E_1 -E_2-E_3)/3$.  This implies that the period of $X$ is contained in 
$\B^9\cap \D_{-4/3}$.  We denote the corresponding Heegner divisor 
$\B^9\cap \D_{-4/3}$ by $\H_{\mathrm{vt}}$. We can summarize the above discussion as follows.
\begin{prop}[{}]
\label{prop:moduli_of_nonhyperelliptic_genus_4}
The moduli space of non-singular and non-hyperelliptic curves of genus $4$ is isomorphic to $(\B^9\setminus \H)/\Gamma_{\mathrm{nh}}$.  The discriminant locus $\H$ 
decomposes into $\H_{\mathrm{n}}$ and $\H_{\mathrm{h}}$, where a general point of
$\H_{\mathrm{n}}$ corresponds to a nodal curve of arithmetic genus $4$ and 
that of $\H_{\mathrm{h}}$  to a pair consisting of a hyperelliptic curve of genus $4$ and two points on it conjugated under the hyperelliptic involution. A general point of the Heegner divisor
$\H_{\mathrm{vt}}$ corresponds to a non-singular curve of genus $4$ with a vanishing theta constant.
\end{prop}

\begin{rem}\label{rem:DM_relation}
The moduli space of genus $4$ curves is closely related to the ancestral Eisenstein Deligne-Mostow
 variety parametrizing $12$ (unordered) points. For this  
let $C$ be a non-singular, non-hyperelliptic curve of genus 4.  
Then the corresponding $K3$ surface $X$ has one or two
elliptic fibration corresponding to the number of $g_3^1$'s.
Let $\M_4^{\mathrm{nh}}(g_3^1)$ be the moduli space of non-singular, non-hyperelliptic curves of genus 4 endowed with a $g_3^1$. This is the double covering of $\M_4^{\mathrm{nh}}$, branched along $\M_4^{\mathrm{vt}}$, the moduli space of non-singular, non-hyperelliptic curves with a vanishing theta constant.  Recall that in the Picard group of the associated $K3$ surface,  $U(3)$ is generated by
two fibrations. Note that, if we denote the automorphism group of $A_L$ preserving $q_L$ by $\O(q_L)$, then
$\O(q_L)$ is isomorphic to $\Z/2\Z$, and this is generated by the isometry interchanging the two fibrations.  Let
\[\widetilde{\Gamma}_{\mathrm{nh}} \defeq \Ker(\Gamma_{\mathrm{nh}} \to \O(q_L)).\]
The above implies that $\M_4^{\mathrm{nh}}(g_3^1)$ is isomorphic to 
$(\B^9\setminus \H)/\widetilde{\Gamma}_{\mathrm{nh}}$. 
Each point in $\M_4^{\mathrm{nh}}(g_3^1)$ gives an elliptic fibration structure on the associated $K3$ surface.  
In the general
case, the fibration has twelve singular fibers of type ${\rm II}$ and hence we obtain 
twelve points of $\P^1$.
Note that the sum of the Euler numbers of fibers coincides with
the Euler number of $X$: $12\cdot e({\rm II}) = 24$.  Dividing by 12, we obtain twelve points with weight
$12$
\[\left(\frac{1}{6}\right)^{12}=\left(\frac{1}{6}, \frac{1}{6}, \frac{1}{6}, \frac{1}{6}, \frac{1}{6}, \frac{1}{6}, \frac{1}{6}, \frac{1}{6}, \frac{1}{6}, \frac{1}{6}, \frac{1}{6}, \frac{1}{6}\right)\]
which is nothing but the notation of Deligne-Mostow.  

Thus we get a rational map 
\[\pi_{\mathrm{geom}} : \M_4^{\mathrm{nh}}(g_3^1) \dashrightarrow \B^9/\Gamma.\]
In the case $C$ is a nodal curve or hyperelliptic curve, the corresponding $K3$ surface has an elliptic fibration with a singular fiber of type IV and 
ten singular fibers of type II.  Note that the Euler number of a fiber of type IV is 4 and thus we have
\[\left(\frac{2}{6}, \frac{1}{6}, \frac{1}{6}, \frac{1}{6}, \frac{1}{6}, \frac{1}{6}, \frac{1}{6}, \frac{1}{6}, \frac{1}{6}, \frac{1}{6}, \frac{1}{6}\right),\]
which is a degenerate case.  In the case $C$ is a general non-singular curve of genus 4 with a vanishing theta
constant, as mentioned above, $X$ has an elliptic fibration with twelve singular fibers of type II.  
Therefore this case is included in the general case.  
\end{rem}

\subsection{Automorphic forms associated to $L=U\oplus U(3)\oplus E_8^{\oplus 2}$}
\label{AutoForm}
We keep the notation of the previous subsection.  Let $M$ be an even lattice of signature $(2,n)$ and let
\[\widetilde{\mathbb{D}}_{M} = \{ \omega \in L\otimes_{\Z} \C\mid \langle \omega, \omega\rangle =0, \ \langle 
\omega\cdot \bar{\omega}\rangle > 0\}^+\]
be the affine cone over $\mathbb{D}_{M}$.
A holomorphic (resp. meromorphic) function $F: \widetilde{\mathbb{D}}_M \to \C$ is called a \textit{holomorphic (resp. meromorphic) 
automorphic form of weight $k$ with respect to a group $\Gamma$ and with character $\chi$},  if $F$ satisfies the following two conditions:
\begin{enumerate}
    \item $F(\gamma(\omega)) = \chi(\gamma) F(\omega)$ for any $\omega \in \widetilde{\mathbb{D}}_M$ and $\gamma\in \Gamma$,
    \item $F(c\cdot \omega) = c^{-k}F(\omega)$ for any $\omega \in \widetilde{\mathbb{D}}_M$ and $c\in \C^{\times}$,
\end{enumerate}
where $\Gamma$ is a finite index subgroup of $\O^+(M)$ and $\chi$ is a character of $\Gamma$. Here $\O^+(M)$ is the subgroup of $\O(M)$ which fixes the connected component $\mathbb{D}_{M}$ (or alternatively the subgroup of elements with real spinor norm 1). 
In the present paper, we call a holomorphic automorphic form simply an \textit{automorphic form}.
Automorphic forms on balls $\B^n$ are defined in an analogous way.

We recall that the group $\widetilde{\O}^+(L)$  is the {\em stable} orthogonal group of $L$ which fixes the chosen component $\mathbb{D}_{L}$. 
Here stable means that it acts trivially on the discriminant $A_L=L^*/L$.
We will show the following theorem.
\begin{thm}\label{AutoMain}
There exists a holomorphic automorphic form $\Psi$ on $\mathbb{D}_L$ of weight $51$ (with respect to $\widetilde{\O}^+(L)$ and with character $\det$) whose zero divisor is given by
\[\D_{-2} + 27\D_{-2/3} + 3\D_{-4/3}.\]
In addition, $\Psi$ is the Borcherds lift of a weakly holomorphic vector-valued modular form.
\end{thm}
We recall that in genus 1 the Koecher principle does not hold. This means that a modular form is not automatically holomorphic at the cusps. This leads to the notion of a  
{\em weakly holomorphic} modular form which is defined as a holomorphic automorphic form on the complex upper half-plane which has poles of finite order at infinity.
Note that the latter part of Theorem \ref{AutoMain} follows from the converse theorem of Bruinier \cite[Theorem 5.11]{Bru02}, \cite[Theorem 1.2]{Bru14}.
In other words, there exists a weakly holomorphic form $f$ with $\Psi=\mathrm{B}(f)$ 
in the notation of \cite{Bru02}.  This sheds new light on $\Psi$.
Note that $L$ satisfies the assumption of \cite[Theorem 5.11]{Bru02} and by the proof therein, we can apply the converse theorem to automorphic forms with any character. 
This modular form $f$ can be explicitly computed by a result of \cite[Theorem 1.1]{Ma19}.
In Subsection \ref{subsec:Ma's result} we give the explicit form of this vector-valued modular form by using theta functions.

Using Theorem \ref{AutoMain}, we can now construct a specific automorphic form on $\B^9$.
\begin{cor}\label{AutoCor}
There exists a holomorphic automorphic form $\Psi_{\mathrm{B}}$ on $\B^9$ of weight $51$ whose zero divisor is given by
\[3(\H_{\mathrm{n}} + 28\H_{\mathrm{h}} + 3\H_{\mathrm{vt}}).\]
\end{cor}
\begin{proof}
If we restrict $\Psi$ to the complex ball $\B^9$, its weight is the same, but the multiplicity of 
zeroes is multiplied by 3.  This follows from the fact that if $r\in L$ with $r^2 <0$, 
then $r, \phi (r)$ and $\phi^2 (r)$ define the same hyperplane in $\B^9$. 
For $\H_{\mathrm{h}}$, the function $\Phi\vert_{\B^9}$ has
multiplicity $(1+27)\times 3$ coming from $\D_{-2}$ and $\D_{-2/3}$.    
\end{proof}
We will give two proofs of Theorem \ref{AutoMain}.  
One of the two proofs uses quasi-pullbacks of automorphic forms due to  \cite{BKPS98}.
This was already used in Casalaina-Martin, Jensen, Laza \cite{CMJL12}.
Their idea, based on a communication with the second named author (see \cite[Acknowledgements]{CMJL12}), is to take the quasi-pullback under the embedding $\B^9\hookrightarrow \D_{II_{2,26}}$. 
In this paper, we pass through $\B^9\hookrightarrow \D_L \hookrightarrow \D_{II_{2,26}}$.
The other method is due to \cite{AF02, Fre03} based on the results of Borcherds \cite{Bor98,Bor99}. 
This
is useful as it is independent of the embedding $\D_L \hookrightarrow \D_{II_{2,26}}$. 
In other words, this strategy can be applied to any lattice.
In this sense, our second proof gives a new proof and interpretation of \cite[Theorem 5.11]{CMJL12}.

\begin{rem}\label{Radu-Laza}
Let $H_*\defeq \H_*/\Gamma_{\mathrm{nh}}$ (in this remark only).    
Casalaina-Martin, Jensen and Laza showed that $\overline{H_{\mathrm{n}}}^{\BB}+14\overline{H_{\mathrm{h}}}^{\BB}+\frac{9}{2}\overline{H_{\mathrm{vt}}}^{\BB}$ is an ample divisor on $\overline{\B^9/\Gamma_{\mathrm{nh}}}^{\BB}$.  
Corollary \ref{AutoCor} recovers this result as follows.
The complex uniformization map $\B^9\to\B^9/\Gamma_{\mathrm{nh}}$ ramifies along the divisors 
$\H_{\mathrm{n}}$, $\H_{\mathrm{h}}$, $\H_{\mathrm{vt}}$ with indices 3, 6, 2.
In combination with Corollary \ref{AutoCor} it follows that 
\[51\L=\overline{H_{\mathrm{n}}}^{\BB}+14\overline{H_{\mathrm{h}}}^{\BB}+\frac{9}{2}\overline{H_{\mathrm{vt}}}^{\BB}.\]
Since $\L$ is ample, this shows the claim.  
\end{rem}

\subsection{First proof of Theorem \ref{AutoMain}}
\label{subsection:first_proof}
Recall that
\[L\cong U^{\oplus 2}\oplus E_8\oplus A_2\oplus E_6\]
and define
\[R\defeq E_6\oplus A_2.\]
It is known that there exist primitive embeddings of $E_6$ and $A_2$ into $E_8$ and
such an embedding is unique up to isometry.
Here we fix a primitive embedding of $R$ into $ II_{2,26}$ such that $E_6$ is a primitive sublattice in one component $E_8$ of $II_{2,26}=U^{\oplus 2}\oplus E_8^{\oplus 3}$ and $A_2$ is a primitive sublattice of another component $E_8$ of 
$II_{2,26}$.  
Since the orthogonal complement of $A_2$ (resp. $E_6$) in $E_8$ is isomorphic to
$E_6$ (resp. $A_2$), $L$ is isomorphic to the orthogonal complement of $R$ in $II_{2,26}$.  
Since $II_{2,26}$ is unimodular, there exists a canonical isomorphism $\varphi : A_L \to A_R$ with
$q_R \circ \varphi = -q_L$ as mentioned at the beginning of this section.

Borcherds \cite{Bor95} constructed a holomorphic automorphic form $\Phi_{12}$ of weight 12 on the domain
$\mathbb{D}_{II_{2,26}}$ such that 
$\Phi_{12}$ vanishes exactly along $(-2)$-hyperplanes (= the hyperplane perpendicular to a $(-2)$-vector) 
with multiplicity 1. Under the inclusion $L \subset II_{2,26}$ the domain
$\mathbb{D}_L$ is naturally embedded in $\mathbb{D}_{II_{2,26}}$.  
The restriction of $\Phi_{12}$ to $\mathbb{D}_L$ vanishes identically because $R$ contains $(-2)$-vectors.  
However, by first dividing by linear functions corresponding to $(-2)$-vectors in $R$, and then 
by restricting $\Phi_{12}$ to $\mathbb{D}_L$, we obtain an automorphic form $\Psi$ on $\mathbb{D}_L$ 
which is called a \textit{quasi-pullback} of $\Phi_{12}$; see \cite[p. 200]{Bor95}, \cite[Section 2]{BKPS98}.  
Since the weight is increased by 1 each time we divide $\Phi_{12}$ by a linear function, 
the weight of $\Psi$ is equal to the weight of $\Phi_{12}$ plus the number of positive roots 
in $R= A_2\oplus E_6$.  Since the number of positive roots of $A_2$ (resp.  $E_6$) is 3 (resp. 36),  
the weight of $\Psi$ is $12 + 3 + 36 = 51$.

Next, let us consider the zeroes of $\Psi$.  Let
$r\in II_{2,26}$ with $r^2=-2$.  Then $\Phi_{12}$ vanishes along $r^{\perp}$, 
the orthogonal complement of $r$, with multiplicity $1$.  
If $r$ is in $L$, the quasi-pullback $\Psi$ vanishes along $\D_r = r^{\perp}\cap \mathbb{D}_L$ with multiplicity 1, and, in this case, 
$\D_r \subset \D_{-2}$. 
Next we consider the case $r = r_1 + r_2$ where $r_1\in L^*$ with $r_1^2 < 0$ and $r_2\ne 0 \in R^*$.  
Since $R$ is negative definite and there are only two types of non-isotropic
vectors in $A_L \cong A_R$, that is, vectors with norm $-2/3$ or  $-4/3$, we have 
$(r_1)^2 = -2/3, -4/3$ and $(r_2)^2= -4/3, -2/3$, respectively.  
Fix one such $r_1$.  Then take any $r_2 \in R^*$ with $\varphi(r_1+L)= r_2+R$ and $r_1^2 + r_2^2 = -2$.  
Then $r=r_1+r_2 \in II_{2,26}$, and hence $\Phi_{12}$ 
vanishes along $r^{\perp}$.  This implies that
$\Psi$ vanishes along $r_1^{\perp}$ with multiplicity $k$, where $k$ is the number of such $r_2$.  
Since there are exactly two vectors
in $R^*/R$ with the same norm $(-2/3)$ or $(-4/3)$, which are the same up to a sign, 
$k$ is equal to half of the number of vectors in $R^*$ with norm $(-2/3)$ or $(-4/3)$.
If $r_2^2 = -2/3$, then $r_2 \in A_2^*$ and one can calculate the number of such $r_2$ directly,
that is 6.  Hence in this case $k =3$.  If $r_2^2 =-4/3$, then $r_2 \in E_6^*$.  The number of $(-4/3)$-vectors in $E_6^*$ is 54; see \cite[Chapter 4, Subsection 8.3 (122)]{CS99}.
Hence in this case $k =27$.  
Note that 
if $r_1^2 = -2/3$ (resp. $-4/3$), then $r^{\perp}\cap \mathbb{D}_L = \D_{r_1} \subset \D_{-2/3}$ 
(resp. $\subset \D_{-4/3}$).
Thus we have finished a proof of Theorem \ref{AutoMain}.

\subsection{Second proof of Theorem \ref{AutoMain}}
\label{subsec:second_proof}
In the following, for simplicity, we assume that $M$ has even rank, that is an even lattice of signature $(2,2n)$.  
We denote by $S, T$ the standard generators of $\SL(2, \Z)$:
\[T =
\begin{pmatrix}1&1
\\0&1
\end{pmatrix},\ 
S =
\begin{pmatrix}0&-1
\\1&0
\end{pmatrix}.\]
Let $\rho$ be the Weil representation of $\SL(2, \Z)$ on the group ring $\C[A_{M}]\defeq \oplus_{\alpha \in A_M}\C\cdot e_{\alpha}$ defined by
\[\rho(T)(e_{\alpha}) = e^{\pi\sqrt{-1}\ q_M(\alpha)} e_{\alpha}, \quad
\rho(S)(e_{\alpha}) = \frac{\sqrt{-1}^{n-1}}{\sqrt{|A_M|}} \sum_{\beta\in A_M} 
e^{-2\pi\sqrt{-1}\ b_M( \beta, \alpha )} e_{\beta}.\]
A holomorphic function $f:\mathbb{H}\to \C[A_M]$ from the upper-half plane $\mathbb{H}$ to $\C[A_M]$ is called 
a \textit{vector-valued modular form of type $\rho$} if, when we write $f(\tau)=\sum_{\alpha\in A_M}f_{\alpha}(\tau)e_{\alpha}$, each component $f_{\alpha}$ satisfies
\[f_{\alpha}(\tau+1) = e^{\pi\sqrt{-1}\ q_M(\alpha)} f_{\alpha}(\tau), \quad
f_{\alpha}\left(-\frac{1}{\tau}\right) = \tau^{k+1-n}\frac{\sqrt{-1}^{n-1}}{\sqrt{|A_M|}} \sum_{\beta\in A_M} 
e^{-2\pi\sqrt{-1}\ b_M( \beta, \alpha )} f_{\beta}(\tau)\]
and each $f_{\alpha}$ is meromorphic at the cusp.  The weight of $f$ is $k+1-n$.
For simplicity, we denote this by $f=\{f_{\alpha}\}_{\alpha\in A_L}$ and denote by
$f_{\alpha}(\tau) = \sum_m c_{\alpha}(m)q^m$ the Fourier expansion of $f_{\alpha}$.
Note that $c_{\alpha}(m) = c_{-\alpha}(m)$ by the invariance under $S^2$.

In this situation, Borcherds \cite{Bor98} constructed a
meromorphic automorphic form associated with $f$, called the  \textit{Borcherds product}.
\begin{thm}[{(Borcherds \cite[Theorem 13.3]{Bor98})}]\label{BorcherdsProductMarch}
Let $f=\{f_{\alpha}(\tau) =\sum_m c_{\alpha}(m)q^m\}$ be a weakly holomorphic vector-valued modular form of weight $1-n$ and of type $\rho$.  Assume that $c_{\alpha}(m) \in \Z$ for $m < 0$ and $c_{00}(0)\in 2\Z$.
Then there exists a meromorphic automorphic form $\Phi$ of weight $c_{00}(0)/2$ whose divisor is given by
\[\frac{1}{2}\sum_{\alpha\in A_M}\left(\sum_{m<0,\ m\in q_M(\alpha)+\Z}c_{\alpha}(m)\D_{\alpha, m}\right).\]
\end{thm}

In our case, namely $L= U\oplus U(3)\oplus E_8^{\oplus 2}$, we have $(2,n)=(2,18)$ and $|A_L|=3^2$, and hence
\[f_{\alpha}(\tau+1) = e^{\pi\sqrt{-1}\ q_L(\alpha)} f_{\alpha}(\tau), \quad
f_{\alpha}\left(-\frac{1}{\tau}\right) = \frac{\tau^{k-8}}{3} \sum_{\beta\in A_L} 
e^{-2\pi\sqrt{-1}\ b_L( \beta, \alpha )} f_{\beta}(\tau).\]
\begin{ex}\label{ExampleBorcherdsProd}
Here, we recall two famous examples of the Borcherds products which play a major role for us.
\begin{enumerate}
    \item (Borcherds \cite[Section 10, Example 2]{Bor95})  Let $M= II_{2,26}$ .
Then $A_M = 0$ and a vector-valued modular form is a usual modular form.
If we take the modular form 
\[1/\Delta(\tau) = q^{-1} + 24 + \cdots \] 
of weight $-12$, where $\Delta(\tau)=\eta(\tau)^{24}$ and $\eta(\tau)$ is Dedekind eta function, we get a holomorphic automorphic form $\Phi_{12}$ with respect to $\O^+(II_{2,26})$ and the character $determinant$ whose weight is $12$ (= half of the constant term $24$) and which vanishes along the hyperplane
perpendicular to a vector with the norm $-2$ (= 2 times the exponent of a negative power of $q$)  with multiplicity $1$ (= the Fourier coefficient of the negative power of $q$). 
This is called the \textit{Borcherds form}, used in Proposition \ref{prop:modular_forms_Allcock_quasi_pullback}.
\item  (Allcock \cite[Proof of Theorem 7.1]{All00})
  If we consider $M=II_{2,18}$ and 
take a modular form of weight $-8$: 
\[E_4(\tau)/\Delta(\tau) = q^{-1} + 264 + \cdots, \]
we get a holomorphic automorphic form on $\D_{II_{2,18}}$ of weight $132$ with respect to $\O^+(II_{2,18})$ whose zero divisor is $\D_{-2}$.
By the construction of Allcock \cite[Theorem 7.1]{All00} this automorphic form has trivial character
\end{enumerate}
\end{ex}
In general, there are as good as no methods to construct vector-valued modular forms.
Instead of constructing such an $f$ directly, 
we employ another approach, which goes back to \cite{Bor99, AF02, Fre03}, 
to show the existence of an automorphic form.

For this, we first introduce the obstruction space, which is a complex vector space, consisting of all vector-valued modular forms $\{ f_{\alpha} \}_{\alpha \in A_{L}}$ of weight $(2+18)/2=10$ and with respect to
the dual representation $\rho^*$ of $\rho$:
\[f_{\alpha}(\tau + 1) = e^{-\pi \sqrt{-1}\ q_L(\alpha)}f_{\alpha}(\tau), \quad
f_{\alpha}\left(-\frac{1}{\tau}\right) = \frac{\tau^{10}}{3} 
\sum_{\beta\in A_L} e^{2\pi \sqrt{-1} \ b_L(\alpha, \beta)} f_{\beta}(\tau).\]
We shall apply the next theorem to show the existence of 
some Borcherds products. 
\begin{thm}[{Borcherds \cite{Bor99}, Allcock-Freitag \cite{AF02}, Freitag \cite[Theorem 5.2]{Fre03}}]
\label{freitag}
 A linear combination
\[\sum_{\substack{\alpha \in A_{L}/\pm 1\\ n<0}} c_{\alpha , n} \D_{\alpha, n}\quad   (c_{\alpha, n} \in \Z)\]
is the divisor of a meromorphic automorphic form on $\mathbb{D}_L$ 
of weight $k$ if for every cusp form 
\[f = \{f_{\alpha}\}_{\alpha\in A_{L}}\]
where
\[ f_{\alpha}(\tau) = \sum_{n \in \Q} a_{\alpha, n} e^{2\pi \sqrt{-1} n \tau}\quad (a_{\alpha, n} \in \C)\] 
is in the obstruction space, the relation
\[\sum_{\substack{\alpha \in A_{L}\\ n<0}} a_{\alpha, -n/2}c_{\alpha, n} = 0\]
holds.  In this case, the weight $k$ is given by
\[k = \sum_{\substack{\alpha \in A_{L}\\ n\in \Z}} b_{\alpha, n/2}c_{\alpha, -n}\]
where $b_{\alpha, n}$ are the Fourier coefficients of the Eisenstein series in the obstruction 
space with the constant term
$b_{0, 0} = -1/2$ and $b_{\alpha, 0} = 0$ for $\alpha \ne 0$.  If all
$c_{\alpha , n}$ are non-negative, then it is the divisor of a holomorphic automorphic form.
\end{thm}
For each $u \in A_{L}$, we denote by $m_{0}$ (resp. $m_{1}$, $m_2$) the number of vectors 
$v \in A_{L}$ of given type with $b_{L}(u, v) \equiv 0$ 
(resp. $1/3, 2/3$).   Then $m_{0}, m_{1}, m_2$ are given in  Table 1.
\begin{table}[h]
\[
\begin{array}{r|lllllllllllllllllllllll}
u& 00&00&00&00&0&0& 0&0&4/3&4/3&4/3&4/3&2/3&2/3&2/3&2/3\\
v&00&0&4/3&2/3&00&0&4/3&2/3&00&0&4/3&2/3&00&0&4/3&2/3\\ \hline
m_0&1&4&2&2&1&2&0&0&1&0&0&2&1&0&2&0\\
m_1&0&0&0&0&0&1&1&1&0&2&1&0&0&2&0&1\\
m_2&0&0&0&0&0&1&1&1&0&2&1&0&0&2&0&1\\
\end{array}
\]
\caption{}
\label{vectors}
\end{table}
In the following, we shall study the divisors 
\[\sum_{\substack{\alpha \in A_{L}/\pm 1\\ n<0}} c_{\alpha , n} \D_{\alpha, n}\]
where $c_{\alpha, n}$ depends only on the type of $\alpha$.
Thus we consider the 4-dimensional representation 
\[V\defeq\{e_{\alpha}+e_{-\alpha}\mid \alpha\in A_L\}\]
induced by $\rho^*$. 
It follows from 
Table \ref{vectors} that the following is the matrix representation of $V$:
\begin{equation}\label{dualRep}
\rho^*(T) =
\begin{pmatrix}1&0&0&0
\\0&1&0&0
\\0&0&\omega&0
\\0&0&0&\omega^2
\end{pmatrix},\quad
\rho^*(S) =  \frac{1}{3}
\begin{pmatrix}1&4&2&2
\\1&1&-1&-1
\\1&-2&-1&2
\\1&-2&2&-1
\end{pmatrix}.
\end{equation}
We denote by 
\[h_{00},\ h_{0},\ h_{4/3},\ h_{2/3}\]
the sum of $f_{\alpha}$'s according to their types.  
\begin{lem}\label{dim}
The dimension of the space of modular forms of weight $10$ and of 
type $\rho^{*}$ is $4$. The dimension of the space of the Eisenstein series of weight $10$ and of type $\rho^{*}$ is $2$.
\end{lem}
\begin{proof}
In general, the dimension of the space of modular forms of weight $k > 2$ and of type $\rho^*$
is given by
\[d + \frac{dk}{12} - \alpha (e^{k\pi\sqrt{-1}/2} \rho^*(S)) - \alpha ((e^{k\pi\sqrt{-1}/3}\rho^*(ST))^{-1}) - \alpha (\rho^*(T))\]
by \cite[Secion 4]{Bor98}, \cite{AF02} and \cite[Prop. 2.1]{Fre03}.
Here 
\[d\defeq\dim \{ x \in V \mid \rho^*(-S)x \defeq (-1)^kx\}\]
and 
\[\alpha (A) \defeq \sum_{\lambda = e^{2\pi\sqrt{-1}t}} t\]
where $\lambda$ runs through all eigenvalues of the unitary matrix $A$ with $0\leq t <1.$

In our situation, $k=10$ and $d= \dim (V) = 4$.  An elementary calculation shows that
\[\alpha (e^{10\pi\sqrt{-1}/2} \rho^*(S)) = 1, \ \alpha ((e^{10\pi\sqrt{-1}/3}\rho^*(ST))^{-1}) = 4/3, \ \alpha (\rho^*(T)) = 1.\]
On the other hand, the space of Eisenstein series is isomorphic to the subspace of $V$ given by
\[\rho^*(T)(x) =x, \ \rho^*(-S)(x) = (-1)^kx;\]
see \cite[Remark 2.2]{Fre03}.
This proves the assertion. 
\end{proof}
Next, we shall determine 
the Eisenstein series $\{ h_{\alpha}\}_{\alpha\in A_L}$ 
of weight 10 and of type $\rho^*$.
By (\ref{dualRep}), the functions $\{ h_{\alpha}\}_{\alpha\in A_L}$ should satisfy the following:
\begin{align*}\begin{cases}
h_{00}(\tau + 1) &= h_{00}(\tau), \ h_{0}(\tau + 1) = h_{0}(\tau), 
\ h_{4/3}(\tau + 1) = \omega h_{4/3}(\tau),\ 
h_{2/3}(\tau + 1) = \omega^2 h_{2/3}(\tau),\smallskip\\
h_{00}(-1/\tau) &= \frac{\tau^{10}}{3} (h_{00} + h_{0} + h_{4/3} +
h_{2/3}),\smallskip\\
h_{0}(-1/\tau) &= \frac{\tau^{10}}{3} (4h_{00} + h_{0} - 2h_{4/3} 
-2 h_{2/3}),\smallskip\\
h_{4/3}(-1/\tau) &= \frac{\tau^{10}}{3} (2 h_{00} - h_{0} - h_{4/3} +
2 h_{2/3}),\smallskip\\
h_{2/3}(-1/\tau) &= \frac{\tau^{10}}{3} (2h_{00} - h_{0} + 2h_{4/3} -
h_{2/3}).
\end{cases}
\end{align*}
Let 
\[E_{1} = G_{10}^{(0,1)}(\tau),\ E_{2} = G_{10}^{(1,0)}(\tau),\ E_{3} = G_{10}^{(1,1)}(\tau),\ 
E_{4} = G_{10}^{(1,2)}(\tau)\]
be the Eisenstein series of weight 10 and level 3; Koblitz \cite[p.131]{Kob93}.  
Here $(0,1)$, $(1,0)$, $(1,1)$, $(1,2) \in (\Z/3\Z)^2$.
The actions of $S, T$ on $\{E_i\}_{\alpha\in A_L}$ are as follows:
\[S: E_1\to \tau^{10}E_2\to \tau^{20}E_1,\ E_3\to \tau^{10}E_4\to \tau^{20}E_3; \quad T(E_1)=E_1,\ T: E_2\to E_3\to E_4\to E_2.\]
Thus the Eisenstein series $\{h_{\alpha}\}_{\alpha\in A_L}$ are given by
\[\begin{cases}
h_{00} = a E_{1} + \frac{a+b}{3}(E_{2} + E_{3} + E_{4}),\cr
h_{0} = b E_1 + \frac{4a + b}{3}(E_{2} + E_{3} + E_{4}),\cr
h_{4/3} = \frac{2a - b}{3}(E_{2} +\omega^2 E_{3} + \omega E_{4}),\cr
h_{2/3} =  \frac{2a - b}{3}(E_{2} +\omega E_{3} +\omega^2 E_{4}),\cr
\end{cases}\]
where $a, b\in\C$ are parameters.
We assume that the constant terms of 
$h_{00}$, $h_0$, $h_{4/3}$, $h_{2/3}$ are $-1/2, 0,0,0$, respectively, by choosing $a, b$ appropriately.
On the other hand, the Fourier expansions of $E_i$ are as
follows, see \cite[Chapter III, Section 3, Proposition 22]{Kob93}:
\begin{align*}\begin{cases}
E_{1} &= \zeta^{1}(10) + \zeta^{-1}(10) + c_{10}((\omega + \omega^2)q + \cdots),\smallskip\\
E_{2} &= c_{10}(q^{1/3} + (2^{9}+1) q^{2/3} + 3^9 q + \cdots ),\smallskip\\
E_{3} &= c_{10}(\omega q^{1/3} + (2^{9}+1)\omega^2 q^{2/3} + 3^9 q + \cdots ),\smallskip\\
E_{4} &= c_{10}(\omega^2 q^{1/3} + (2^{9}+1)\omega q^{2/3} + 3^9 q+ \cdots ), 
\end{cases}
\end{align*}
where 
\[c_{10} = \frac{-20 \zeta(10)}{3^{10}B_{10}} = -\frac{10 \cdot (2\pi)^{10}}{3^{10}\cdot 10!},\  
\zeta^{1}(10) + \zeta^{-1}(10) = \left(\frac{2\pi}{3}\right)^{10}\cdot \frac{1}{10!}\cdot \frac{2\cdot 5 \cdot 11\cdot 61}{3}.\]
Therefore the desired Eisenstein series is given by
\begin{align}\label{eisen1}
\begin{cases}
h_{00} &= a E_{1} + \frac{a}{3} (E_{2} + E_{3} + E_{4}) 
= -\frac{1}{2} + \frac{3^{10}-3}{2\cdot 11\cdot 61}q + \cdots ,\smallskip\\
h_{0} &= \frac{4a}{3} (E_{2} + E_{3} + E_{4})
= \frac{2\cdot 3^{10}}{11\cdot 61}q + \cdots ,\smallskip\\
h_{4/3} &= \frac{2a}{3} (E_{2} +\omega^2 E_{3} + \omega E_{4})
= \frac{3}{11\cdot 61} q^{1/3} + \cdots ,\smallskip\\
h_{2/3} &= \frac{2a}{3} (E_{2} +\omega E_{3} + \omega^2 E_{4})
= \frac{3(2^9+1)}{11\cdot 61} q^{3/2} + \cdots, 
\end{cases}
\end{align}
where
\[a = -\frac{1}{2}\left(\frac{3}{2\pi}\right)^{10} \frac{3\cdot 10 !}{2\cdot 5 \cdot 11 \cdot 61}.\] 
The obstruction space has dimension 4 and it
contains a 2-dimensional subspace of cusp forms; see Lemma \ref{dim}.
In order to calculate the cusp forms in the obstruction space we consider the following two types:

\smallskip

(A) $\{ \eta^{8}(\tau ) g_{\alpha} \}_{\alpha\in A_L}$,
\smallskip

(B) $\{ \eta^{16}(\tau) g_{\alpha} \}_{\alpha\in A_L}$,
\smallskip

\noindent
where $\eta(\tau)$ is the Dedekind eta function.
\medskip

{Case (A)}:
we denote by
\[F_{1} = G_{6}^{(0,1)}(\tau),\ F_{2} = G_{6}^{(1,0)}(\tau),\
F_{3} = G_{6}^{(1,1)}(\tau),\ F_{4} = G_{6}^{(1,2)}(\tau)\]
the Eisenstein series of weight 6 and level 3; Koblitz \cite[p.131]{Kob93}.  
Their Fourier expansion is as follows, see (\cite[Chapter III, Section 3, Proposition 22]{Kob93}:
\begin{align*}
\begin{cases}
F_{1} &= \zeta^1(6)+\zeta^{-1}(6) - c_{6}q + \cdots , \\
F_{2} &= c_{6}(q^{1/3} + (2^{5} + 1) q^{2/3} + 3^{5}q +
\cdots ),\\
F_{3} &= c_{6}(\omega q^{1/3} + \omega^2(2^{5}+1)q^{2/3} + 3^5q + \cdots ),\\
F_{4} &= c_{6}(\omega^2q^{1/3} + \omega (2^5+1) q^{2/3} + 3^5q +
\cdots ),
\end{cases}
\end{align*}
where 
\[c_{6} = \frac{-6 (2\pi)^{6}}{3^{6}\cdot 6!},\ \zeta^1(6) + \zeta^{-1}(6) = \frac{26 (2\pi)^6}{3^7\cdot 6!}.\]
The action of $S, T$ on the functions $\{F_i\}$ is as follows: 
\[S: F_1 \to \tau^6F_2\to \tau^{12}F_1,\ F_3\to \tau^6F_4\to \tau^{12}F_3; \quad T(F_1)=F_1, \ T: F_2\to F_3\to F_4\to F_2.\]
Recall that $\eta^8(\tau+1)=\omega\eta^8(\tau)$ and $\eta^8(-1/\tau)=\tau^4\eta^8(\tau)$.
If we write $h_{\alpha} = \eta^{8}(\tau ) g_{\alpha}$, 
we need to find $\{ g_{\alpha}\}_{\alpha\in A_L}$ satisfying:
\begin{align*}
\begin{cases}
    g_{00}(\tau + 1) &= \omega^2 g_{00}(\tau), \ g_{0}(\tau + 1) = \omega^2 g_{0}(\tau),\
g_{4/3}(\tau + 1) = h_{4/3}(\tau),\ 
g_{2/3}(\tau + 1) = \omega g_{2/3}(\tau),\smallskip\\
g_{00}(-1/\tau) &= \frac{\tau^{6}}{3} (g_{00} + g_{0} + g_{4/3} + g_{2/3}),\smallskip\\
g_{0}(-1/\tau) &= \frac{\tau^{6}}{3} (4g_{00} + g_{0} - 2g_{4/3} - 2g_{2/3}),\smallskip\\
g_{4/3}(-1/\tau) &= \frac{\tau^{6}}{3} (2g_{00} - g_{0} - g_{4/3} + 2g_{2/3}),\smallskip\\
g_{2/3}(-1/\tau) &= \frac{\tau^{6}}{3} (2g_{00} - g_{0} + 2g_{4/3} - g_{2/3}).
\end{cases}
\end{align*}

By solving linear equations, we obtain a one-dimensional subspace of cusp forms and their expansions as follows: 
\begin{equation}\label{cuspA}
\begin{cases}
h_{00} = a \eta(\tau)^{8}(F_{2} + \omega F_{3} + \omega^2 F_{4}) = 
3(2^{5} + 1)a c_6 q + \cdots ,\cr
h_{0} = -2 a \eta(\tau)^{8}(F_{2} + \omega F_{3} + \omega^2 F_{4}) = 
-6a c_{6}(2^{5}+1)q + \cdots ,\cr
h_{4/3} = a \eta(\tau)^{8} (3F_1 - F_{2} - F_{3} - F_{4}) = 
3a (\zeta^1(6) + \zeta^{-1}(6)) q^{1/3} + \cdots ,\cr
h_{2/3} = 2a \eta(\tau)^{8}(F_{2} + \omega^2 F_{3} + \omega F_{4}) = 
6a c_{6}q^{2/3} + \cdots , \cr
\end{cases}
\end{equation}
where $a\in\C$ is a parameter.

\medskip
{Case (B):}
Let 
\[G_{1} = G_{2}^{(0,1)}(\tau), \ G_{2} = G_{2}^{(1,0)}(\tau),\
G_{3} = G_{2}^{(1,1)}(\tau),\ G_{4} = G_{2}^{(1,2)}(\tau)\]
be the Eisenstein series of weight 2 and level 3, see \cite[Chapter VII, Section 2]{Sch74}. We remark that the
$G_i$ are not holomorphic, but also that their non-analytic parts are all equal.
The actions of $S, T$ on $\{G_{i}\}$ are:
\[S : G_{1} \rightarrow \tau^2G_{2} \rightarrow \tau^4G_{1}, 
G_{3} \rightarrow \tau^2G_{4} \rightarrow \tau^4G_{3}; \quad T(G_{1}) = G_{1},\ G_{2} \rightarrow G_{3} \rightarrow G_{4}
\rightarrow G_{2}.\]
Write $h_{\alpha} = \eta^{16}(\tau)g_{\alpha}$.
Since $\eta^{16}(\tau + 1) = \omega^2\eta^{16}(\tau), \eta^{16}(-1/\tau)
= \tau^{8} \eta^{16}(\tau)$, the functions $\{ g_{\alpha}\}_{\alpha\in A_L}$ satisfiy:
\begin{align*}
\begin{cases}
    g_{00}(\tau + 1) &= \omega g_{00}(\tau),\ g_{0}(\tau + 1) = \omega g_{0}(\tau),\ 
g_{4/3}(\tau + 1) = \omega^2 g_{4/3}(\tau),\ g_{2/3}(\tau + 1) = g_{2/3}(\tau),\smallskip\\
g_{00}(-1/\tau) &= \frac{ \tau^2}{3} (g_{00} + g_{0} + g_{4/3} + g_{2/3}),\smallskip\\
g_{0}(-1/\tau) &= \frac{\tau^{2}}{3} (4g_{00} + g_{0} - 2g_{4/3} - 2g_{2/3}),\smallskip\\
g_{4/3}(-1/\tau) &= \frac{\tau^{2}}{3} (2g_{00} - g_{0} - g_{4/3} + 2g_{2/3}),\smallskip\\
g_{2/3}(-1/\tau) &= \frac{\tau^{2}}{3} (2g_{00} - g_{0} + 2g_{4/3} - g_{2/3}).
\end{cases}
\end{align*}
An easy calculation shows that the $g_{\alpha}$ can be written as a linear combination of the $\{ G_{i} \}$ as follows:
\begin{align*}
\begin{cases}
g_{00} &= a (G_{2} + \omega^2 G_{3} +\omega G_{4}),\\
g_{0} &= -2a (G_{2} +\omega^2 G_{3} +\omega G_{4}),\\
g_{4/3} &= 2a (G_{2} +\omega G_{3} +\omega^2 G_{4}),\\
g_{2/3} &= 3a G_{1} - a(G_{2} + G_{3} + G_{4}),
\end{cases}
\end{align*}
where $a\in\C$ is a parameter.
Note that the non-analytic parts are canceled.
The Fourier coefficients of the $G_{i}$, up to their non-analytic parts, are as follows
(\cite[Chapter VII, Section 2]{Sch74}):
\begin{align*}\begin{cases}
G_{1} &= \frac{4\pi^2}{27}  + 0\cdot q^{1/3} + \cdots  , \smallskip\\
G_{2} &= -\frac{4\pi^2}{9}q^{1/3} + \cdots , \smallskip\\
G_{3} &=  - \frac{4\pi^2 \omega}{9}q^{1/3} + \cdots , \smallskip\\
G_{4} &= - \frac{4\pi^2 \omega^2}{9}q^{1/3} + \cdots . 
\end{cases}
\end{align*}
Thus we have
\begin{align}\label{cuspB}
\begin{cases}
h_{00}&=\eta(\tau)^{16} g_{00} = -\frac{4\pi^2 a}{3} q + \cdots  ,\smallskip\\
h_0&=\eta(\tau)^{16} g_{0} = \frac{8 \pi^2 a}{3} q + \cdots  ,\smallskip\\
h_{4/3}&=\eta(\tau)^{16} g_{4/3} = 0\cdot q + \cdots ,\smallskip\\
h_{2/3}&=\eta(\tau)^{16} g_{2/3} = \frac{4\pi^2 a}{9}q^{2/3} + \cdots .
\end{cases}
\end{align}
\begin{thm}\label{multiplicative}
A divisor
\[ m_{00}\D_{-2} + m_{4/3}\D_{-2/3} + m_{2/3}\D_{-4/3}\]
is a divisor of a meromorphic automorphic form on $\D_L$ if $m_{2/3} = 3 m_{00}$, $m_{4/3} = 27 m_{00}$.
In this case, the weight is given by $51 m_{00}$.
The automorphic from is holomorphic if and only if  $m_{00}>0$.
\end{thm}
\begin{proof}
The first assertion follows from Theorem \ref{freitag} and the equations (\ref{cuspA}), (\ref{cuspB}):
\[-\frac{4\pi^2}{3}m_{00} + \frac{4\pi^2}{9}m_{2/3} = 0,\quad
3(2^5 + 1)c_6m_{00} + 3(\zeta^1(6) + \zeta^{-1}(6))m_{4/3} + 6c_6m_{2/3} =0.\]
By using Theorem \ref{freitag} and the equation (\ref{eisen1}), we can see that 
the weight of $F$ is given by
\[\frac{3^{10} -3}{2\cdot 11\cdot 61}m_{00} + \frac{3}{11\cdot 61}m_{4/3} +\frac {3(2^9 +1)}{11\cdot 61}m_{2/3} = 51m_{00}.\]
\end{proof}
By putting $m_{00}=1$ in Theorem \ref{multiplicative}, 
we have now finished the second proof of Theorem \ref{AutoMain}.

\subsection{Explicit construction of the input of the Borcherds lift}
\label{subsec:Ma's result}
Here we provide an explicit calculation of the vector-valued modular form $f$ which leads to the Borcherds lift $\Psi=B(f)$ which we discussed in Subsection \ref{AutoForm}. 
Our calculations follow  \cite{Ma19}.
Let $L= U\oplus U\oplus E_8\oplus E_6\oplus A_2$ as before.  Then $L$ is the orthogonal complement of
a primitive sublattice $K=E_6\oplus A_2$ in $II_{2,26}$.
By the construction of $\Psi$ as in Subsection \ref{subsection:first_proof}, this form is obtained by the quasi-pull back of the Borcherds form $\Phi_{12}$ to $\D_L$.
Here, we recall that $\Delta(\tau)$ is the Ramanujan delta function: $\Delta(\tau)=q\prod_{n>0}(1-q^n)^{24}$ as in Example \ref{ExampleBorcherdsProd}(1).
The Borcherds form $\Phi_{12}$ is the lift of $1/\Delta$.
As we remark just after Theorem \ref{AutoMain}, \cite[Theorem 1.1]{Ma19} implies that $\Psi$ is the Borcherds lift of a weakly vector-valued modular form $f =\{ f_{\alpha}(\tau)\}_{\alpha\in A_L}$ where $A_L=L^*/L$.
Below, we construct $f$ explicitly.

Let $M$ be a negative definite lattice $M$ and $M^*$ its dual.
Define the theta series associated with $M$:
\[\theta_M(\tau) = \sum_{x\in M}q^{-\langle x, x\rangle},\quad \theta_{M^*}(\tau) = \sum_{x\in M^*}q^{-\langle x, x\rangle},\quad 
\theta_{M+\alpha}(\tau) = \sum_{x\in M^*,\ x +M =\alpha}q^{-\langle x, x\rangle}.\]
Note that 
\[\theta_{M^*}(\tau)= \theta_M(\tau)+ \sum_{\alpha \in A_M,\ \alpha\ne 0}\theta_{M+\alpha}(\tau),\quad \theta_{M+\alpha}=\theta_{M+(-\alpha)}.\]

In case $M= A_2, E_6$, the finite group $M^*/M$ consists of three elements $0, \pm 1$.  Thus we have four theta series
\[\theta_{A_2}(\tau),\ \theta_{A_2+[1]}(\tau)= \theta_{A_2+[-1]}(\tau),\ 
\theta_{E_6}(\tau),\ \theta_{E_6+[1]}(\tau)= \theta_{E_6+[-1]}(\tau).\]
Their Fourier expansions are given by
\begin{align*}
    \theta_{A_2}(\tau) &= 1 + 6q + 6q^3+\cdots,\\
    \theta_{A_2+[1]} &=\theta_{A_2+[2]}=3q^{1/3}+3q^{4/3}+\cdots,\\
    \theta_{E_6}(\tau) &= 1 + 72q + 270q^2+\cdots,\\
    \theta_{E_6+[1]} &=\theta_{E_6+[2]}=27q^{2/3}+216q^{5/3}+\cdots.
\end{align*}
(see \cite[p.111, Equations (61), (63); p.127, (121),(123)]{CS99}.)

Ma \cite[Theorem 1.1]{Ma19} gave a recipe for how to get a weakly holomorphic modular form from the data
of lattices (in our case, $L=U^2\oplus E_8\oplus E_6\oplus A_2$, $K=E_6\oplus A_2$ and $II_{2,26}$).  
Applying this,
the explicit form of $f$ is given by
\begin{align*}
    f_{\alpha}(\tau)&= \frac{\theta_{A_2}(\tau)\theta_{E_6}(\tau)}{\Delta(\tau)} = q^{-1} + 102 +\cdots  &(\alpha:\ \mathrm{type}\ (00)),\\
    f_\alpha(\tau) &= \frac{\theta_{E_6+[1]}(\tau)\theta_{A_2+[1]}}{\Delta(\tau)}= 81+729q+\cdots  &(\alpha: \ \mathrm{type}\ (0)),\\
    f_\alpha(\tau) &=\frac{\theta_{E_6+[1]}(\tau)}{\Delta(\tau)}= 27q^{-1/3}+648q^{2/3}+\cdots  &( \alpha: \ \mathrm{type}\ (4/3)),\\
    f_\alpha(\tau)&=\frac{\theta_{A_2+[1]}(\tau)}{\Delta(\tau)}= 3q^{-2/3}+75q^{1/3}+\cdots  & (\alpha: \ \mathrm{type}\ (2/3)).
\end{align*}
Now, it follows from Theorem \ref{BorcherdsProductMarch} that the corresponding Borcherds product $\Psi'$ of $f$ has weight $51$ (=  half of the constant term of $f_{00}$),
has divisors along the Heegner divisors $\D_{-2}$ with multiplicity 1 (= the coefficient of $q^{-1}$ of $f_{00}$),
$\D_{-2/3}$ with multiplicity 27 (= the coefficient of $q^{-1/3}$ of $f_{4/3}$) and
$\D_{-4/3}$ with multiplicity 3 (= the coefficient of $q^{-2/3}$ of $f_{2/3}$).
It thus follows from the Koecher principle that $\Psi=\Psi'$ up to constant.
\subsection{Allcock's automorphic form}
\label{subsection:allcock}
 Allcock constructed an automorphic form directly on the ball quotient $\B^9/\Gamma$  using Borcherds infinite products. 
This is similar to Example \ref{ExampleBorcherdsProd} (2).
By pulling-back the automorphic form to $\B^9$ constructed there, one obtains a cusp form $\Psi_1$ on $\B^9$, which has weight 132 vanishing on $\H$ with multiplicity 3. 
This can also be interpreted as a Borcherds product on unitary groups \cite{Hof14}.
Here we interpret this as a quasi-pullback, which can be considered as a special case of \cite[Theorem 1.1]{Ma19}.
In our case we work with the Hermitian lattice $\Lambda$ of signature $(1,9)$  whose associate integral lattice is $II_{2,18}=U^{\oplus 2}\oplus E_8^{\oplus 2}$ (see Subsection \ref{AutoForm}). 
We have a natural embedding $\B^9  \hookrightarrow \mathbb{D}_{II_{2,18}}$.
 
\begin{prop}
\label{prop:modular_forms_Allcock_quasi_pullback}
Up to constant multiple, the automorphic form $\Psi_1$ coincides with the pullback to $\B^9$ of the quasi-pullback of the Borcherds form $\Phi_{12}$ on the $26$-dimensional type IV domain $\D_{II_{2,26}}$.
\end{prop}
\begin{proof}
 We use  the embedding $\B^9  \hookrightarrow \mathbb{D}_{II_{2,18}}$ from above.
Then, since the number of roots in $E_8$ is 240, the quasi-pullback of the Borcherds form $\Phi_{12}$, the Borcherds lift of the inverse of the Ramanujan delta function on $\mathbb{D}_{II_{2,18}}$, has 
weight $12 + 120 = 132$.
Here we note that the quasi-pullback increases the weight by $1/2$ the number of roots orthogonal to $U^{\oplus 2}\oplus E_8^{\oplus 2}$. As the orthogonal complement is $E_8$, this is  $1/2$ of the 240 roots of $E_8$. 
We also remark that the level of this automorphic form is $\Gamma$ because $\Lambda$ is unimodular.
Then, pulling back to $\B^9$ via the embedding $\B^9\hookrightarrow \mathbb{D}_{II_{2,18}}$, we obtain
 that the resulting cusp form on $\B^9$ has weight $132$.
 Its vanishing order along the Heegner divisors on $\B^9$ has multiplicity 3 because $|\OO_{\Q({\omega})}^{\times}/\pm|=3$.

Since $\Psi_1$ and the automorphic form constructed above have the same weight and vanishing loci it follows they are sections of the same line bundle defining the same divisor
differing only by a non-zero multiplicative scalar. 
This gives a description to $\Psi_1$ as (the pullback of) the quasi-pullback of the Borcherds form.
Note that by applying Ma's method, we can recover the weakly holomorphic modular form $f= \theta_{E_8}/\Delta = E_4/\Delta$ (Example \ref{ExampleBorcherdsProd}).
\end{proof}

To conclude the present paper, we compare the automorphic form we constructed with that of Allcock.
\begin{rem}
\label{rem:prob_modular_forms}
The two arithmetic subgroups $\Gamma$ and $\Gamma_{\mathrm{nh}}$ are commensurable by the second named author's result \cite[Theorem 3]{Kon02}, and hence $\Gamma$ and $\Gamma_{\mathrm{nh}}\cap\widetilde{\O}^+(L)$ are also commensurable.
Thus the two automorphic forms  $\Psi_1$ and $\Psi_2\defeq\Psi_{\mathrm{B}}$ share the following properties:
\begin{enumerate}
    \item They are automorphic forms on $\B^9$
     with respect to the same arithmetic subgroup.
    \item Both forms are obtained via a quasi-pullback of the Borcherds form $\Phi_{12}$ (Subsection \ref{subsection:first_proof} and Proposition \ref{prop:modular_forms_Allcock_quasi_pullback}).
    \item Moreover, both are obtained by the Borcherds lift of weakly holomorphic modular forms (Theorem \ref{AutoMain} and Example \ref{ExampleBorcherdsProd}).
\end{enumerate}
It would be interesting to investigate the relationship between $\M_4$ and the ancestral Eisenstein
Deligne-Mostow variety $\M^{\GIT}$ parametrizing 12 (unordered) points further, especially between our and Allcock's automorphic forms.  More concretely, we present the following question.  
Consider $II_{2,18}$ as an overlattice of $L=U(3)\oplus U\oplus E_8^{\oplus 2}$.  Then $\widetilde{\O}^+(L)\subset \O^+(II_{2,18})$ and hence $\widetilde{\Gamma}_{\mathrm{nh}} \subset \Gamma$.  
Thus, on the one hand, we have a map 
\[\pi_{\mathrm{arith}} : \B^9/\widetilde{\Gamma}_{\mathrm{nh}} \to \B^9/\Gamma\]
which gives us 
a lattice theoretic relation between $\M_4^{\mathrm{nh}}(g_3^1)$ and the Deligne--Mostow variety $\B^9/\Gamma$.
On the other hand, we gave a rational, geometrically defined map $\pi_{\mathrm{geom}}$ from $\M_4^{\mathrm{nh}}(g_3^1)$ to $\B^9/\Gamma$
in Remark \ref{rem:DM_relation}.
\end{rem}

\begin{que}
\label{que:comparisonmaps}
What can one say about the relation between the maps $\pi_{\mathrm{arith}}$ and $\pi_{\mathrm{geom}}$ ?
\end{que}


\begin{thebibliography}{99}

\bibitem[Ale96]{Ale96}
V. Alexeev,
\textit{Log canonical singularities and complete moduli of stable pairs}, arXiv:alg-geom/9608013, 1996.

\bibitem[AE23]{AE23}
V. Alexeev, P. Engel,
\textit{Compact moduli of K3 surfaces}, Ann. of Math. (2) {\bf 198} (2023), no. 2, 727–789.

\bibitem[AEH24]{AEH24}
V. Alexeev, P. Engel, C. Han,
\textit{Compact moduli of K3 surfaces with a nonsymplectic automorphism},
Trans. Amer. Math. Soc. Ser. {\bf B 11} (2024), 144–163.

\bibitem[All00]{All00}
D. Allcock,
\textit{The Leech lattice and complex hyperbolic reflections}, 
Invent. Math., {\bf 140} (2000), no. 2, 283--301.

\bibitem[ACT11]{ACT11} D. Allcock, J. A. Carlson, and D. Toledo,
\textit{The moduli space of cubic threefolds as a ball quotient},
Mem. Amer. Math. Soc., {\bf 209} (2011), no. 985, xii+70.

\bibitem[AF02]{AF02} D. Allcock, E. Freitag, 
\textit{Cubic surfaces and Borcherds products}, 
Comm. Math. Helv., {\bf 77} (2002), 270--296.

\bibitem[AMRT10]{AMRT10}
A. Ash, D. Mumford, M. Rapoport, Y.S. Tai,
\textit{Smooth compactifications of locally symmetric varieties},
second ed., Cambridge Mathematical Library, Cambridge University Press, Cambridge, 2010, With the collaboration of Peter Scholze.

\bibitem[AL02]{AL02}
D. Avritzer, H. Lange,
\textit{The moduli spaces of hyperelliptic curves and binary forms},
Math. Z., {\bf 242} (2002), 615--632.

\bibitem[Bal11]{Bal11}
M. Ballard, 
\textit{Derived categories of sheaves on singular schemes with an application to reconstruction}, 
Adv. Math., {\bf 227} (2011), No. 2, 895--919 (2011).

\bibitem[BFK19]{BFK19}
M. Ballard, D. Favero, L. Katzarkov, 
\textit{Variation of geometric invariant theory quotients and derived categories}, 
J. Reine Angew. Math., {\bf 746} (2019), 235--303.

\bibitem[Beh12]{Beh12}
N. Behrens,
\textit{Singularities of ball quotients},
Geom. Dedicata {\bf 159} (2012), 389--407.

\bibitem[BO01]{BO01}
A. Bondal, D. Orlov,
\textit{Reconstruction of a variety from the derived category and groups of autoequivalences}, 
Compos. Math., {\bf 125} (2001), No. 3, 327--344.

\bibitem[BO02]{BO02}
 A. Bondal,  D. Orlov,
 \textit{Derived categories of coherent sheaves},
International Congress of Mathematicians, page 47, 2002

\bibitem[Bor95]{Bor95}  R.\ Borcherds, \textit{Automorphic forms on ${\rm O}_{s+2,s}({\R})$ and infinite products}, Invent. math., {\bf 120} (1995), 161--213.  

\bibitem[Bor98]{Bor98} R.\ Borcherds,  \textit{Automorphic forms with singularities on
Grassmannians}, Invent. Math., {\bf 132} (1998), 491--562.

\bibitem[Bor99]{Bor99} R.\ Borcherds,  \textit{The Gross-Kohnen-Zagier theorem in higher dimensions}, Duke Math. J., {\bf 97} (1999), 219--233.

\bibitem[BKPS98]{BKPS98}  R.\ Borcherds, L. Katzarkov, T. Pantev, N.I. Shepherd-Barron, \textit{Families of $K3$ surfaces}, J. Algebraic Geometry {\bf 7} (1998), 183--193.

\bibitem[Bru02]{Bru02}
J. Bruinier,
\textit{Borcherds products on $\O(2,l)$
 and Chern classes of Heegner divisors},
Lecture Notes in Mathematics. {\bf 1780}. Springer (2002).

\bibitem[Bru14]{Bru14}
J. Bruinier,
\textit{On the converse theorem for Borcherds products},
J. Algebra {\bf 397} (2014), 315--342.

\bibitem[CMGHL23]{CMGHL23a}
S. Casalaina-Martin, S. Grushevsky, K. Hulek, R. Laza,
\textit{Cohomology of the Moduli Space of Cubic Threefolds and Its Smooth Models}, 
Mem. Amer. Math. Soc., {\bf 1395} (2023).

\bibitem[CMGHL24]{CMGHL23b}
S. Casalaina-Martin, S. Grushevsky, K. Hulek, R. Laza,
\textit{Non-isomorphic smooth compactifications of the moduli space of cubic surfaces}, 
Nagoya Math. J., {\bf 254} (2024), 315--365.


\bibitem[CMJL12]{CMJL12}
S. Casalaina-Martin, D. Jensen, R. Laza,
\textit{The geometry of the ball quotient model of the moduli space of genus four curves}, 
Compact moduli spaces and vector bundles, 107--136,
Contemp. Math., {\bf 564}, Amer. Math. Soc., Providence, RI, 2012.

\bibitem[CT20a]{CT20a}
A. M. Castravet, J. Tevelev, 
\textit{Derived category of moduli of pointed curves. I},
Algebr. Geom., {\bf 7} (2020), No. 6, 722--757.

\bibitem[CT20b]{CT20b}
A. M. Castravet, J. Tevelev, 
\textit{Derived category of moduli of pointed curves. II}, arXiv:2002.02889.

\bibitem[Con85]{Con85}
J. H. Conway et al., 
\textit{Atlas of finite groups}, Oxford 1985.

\bibitem[CS99]{CS99} J.\ H.\ Conway, N.\ J.\ A.\ Sloane,
\textit{Sphere packings, lattices and groups}, 3rd ed., Springer-Verlag, Berlin, Heidelberg, New York 1999.

\bibitem[DM86]{DM86}
P. Deligne, G. D. Mostow,
\textit{Monodromy of hypergeometric functions and nonlattice integral monodromy}, Inst. 
Hautes Etudes Sci. Publ. Math., {\bf 63} (1986), 5--89.

\bibitem[DvGK05]{DvGK05}
I. Dolgachev, B. van Geemen, S. Kond\={o}, 
\textit{A complex ball uniformization of the moduli space of cubic surfaces via periods of K3 surfaces}, 
J. Reine Angew. Math. {\bf 588}, 99-148 (2005).


\bibitem[Fre03]{Fre03}  E.\ Freitag, \textit{Some modular forms related to cubic surfaces}, Kyungpook Math. J., {\bf 43} (2003), 433--462.

\bibitem[Fuj17]{Fuj17}
O. Fujino,
\textit{Foundations of the minimal model program},
MSJ Memoirs, {\bf 35}, Mathematical Society of Japan, Tokyo, 2017.

\bibitem[GKS21]{GKS21}
P. Gallardo, M. Kerr, L. Schaffler,
\textit{Geometric interpretation of toroidal compactifications of moduli
of points in the line and cubic surfaces},
Adv. Math., {\bf 381} (2021), Paper No. 107632.

\bibitem[GHS08]{GHS08}
V. Gritsenko, K. Hulek, G.K., Sankaran
\textit{Hirzebruch-{Mumford} proportionality and locally symmetric varieties of orthogonal type}, 
Doc. Math., {\bf  13} (2008), 1--19.

\bibitem[GH17]{GH17}
S. Grushevsky, K. Hulek,
\textit{The intersection cohomology of the Satake compactification of $\A_g$ for $g\leq 4$}, 
Math. Ann., {\bf  369} (2017), 1353--1381.

\bibitem[Has03]{Has03}
B. Hassett,
\textit{Moduli spaces of weighted pointed stable curves},
Adv. Math., {\bf 173} (2003), no. 2, 316--352.

\bibitem[HL02]{HL02}
G. Heckman, E. Looijenga,
\textit{The moduli space of rational elliptic surfaces}, 
Algebraic geometry 2000, Azumino (Hotaka), 185--248,
Adv. Stud. Pure Math., {\bf 36}, Math. Soc. Japan, Tokyo, 2002.

\bibitem[Hof14]{Hof14}
E. Hofmann,
\textit{Borcherds products on unitary groups},
Math. Ann., {\bf 358} (2014), No. 3--4, 799--832.

\bibitem[HM25]{HM25}
K. Hulek, Y. Maeda,
\textit{Revisiting the moduli space of 8 points on $\P^1$}, 
Adv. Math. {\bf 463} (2025): 110126..

\bibitem[Kat84]{Kat84}
P. I. Katsylo,
\textit{The rationality of moduli spaces of hyperelliptic curves}, 
Math. USSR, Izv. {\bf 25}, 45-50 (1985).

\bibitem[Kaw02]{Kaw02}
Y. Kawamata,
\textit{D-equivalence and K-equivalence},
J. Differ. Geom., {\bf 61} (2002), No. 1, 147-171.

\bibitem[Kaw04]{Kaw04}
Y. Kawamata, 
\textit{Equivalences of derived categories of sheaves on smooth stacks.},
Am. J. Math., {\bf 126} (2004), No. 5, 1057--1083.

\bibitem[Kaw05]{Kaw05}
Y. Kawamata,
\textit{Log crepant birational maps and derived categories},
J. Math. Sci., Tokyo {\bf 12} (2005), No. 2, 211--231.

\bibitem[Kaw06]{Kaw06}
Y. Kawamata, 
\textit{Derived categories of toric varieties},
Michigan Math. J., {\bf 54} (2006) (3), 517--535.

\bibitem[Kaw13]{Kaw13}
Y. Kawamata, 
\textit{Derived categories of toric varieties II},
Michigan Math. J., {\bf 62} (2013) (2), 353--363.

\bibitem[Kaw16]{Kaw16}
Y. Kawamata, 
\textit{Derived categories of toric varieties III},
Eur. J. Math., {\bf 2} (2016), No. 1, 196--207.

\bibitem[Kaw18]{Kaw18}
Y. Kawamata,
\textit{Birational geometry and derived categories},
Cao, Huai-Dong (ed.) et al., Celebrating the 50th anniversary of the Journal of Differential Geometry: International Press. Surv. Differ. Geom., {\bf 22} (2018), 291--317.

\bibitem[KO15]{KO15}
K. Kawatani, S. Okawa,
\textit{Nonexistence of semiorthogonal decompositions and sections of the canonical bundle},
arXiv:1508.00682.

\bibitem[KM11]{KM11}
Y.-H. Kiem, H-.B. Moon,
\textit{Moduli spaces of weighted pointed stable rational curves via GIT},
Osaka J. Math., {\bf 48} (2011), no. 4, 1115--1140.

\bibitem[KM13]{KM13}
S. Keel, J. McKernan,
\textit{Contractible extremal rays on $\overline{M}_{0,n}$},
Handbook of moduli. Vol. II, 115--130,
Adv. Lect. Math. (ALM), 25, Int. Press, Somerville, MA, 2013.

\bibitem[Kir84]{Kir84}
F. C. Kirwan,
\textit{Cohomology of quotients in symplectic and algebraic geometry},
Mathematical Notes, vol. 31, Princeton University Press, Princeton, NJ, 1984.

\bibitem[Kir85]{Kir85}
F. C. Kirwan,
\textit{Partial desingularisations of quotients of nonsingular varieties and their
Betti numbers},
Ann. of Math. (2) {\bf 122} (1985), no. 1, 41--85.

\bibitem[Kir89]{Kir89}
F. C. Kirwan,
\textit{Moduli spaces of degree d hypersurfaces in $\P^n$},
Duke Math. J., {\bf 58} (1989), no. 1, 39--78.

\bibitem[KLW87]{KLW87}
F. C. Kirwan, R. Lee, S. H. Weintraub,
\textit{Quotients of the complex ball by discrete groups},
Pacific J. Math., {\bf 130} (1987), no. 1, 115--141.

\bibitem[Kob93]{Kob93}  N.\ Koblitz,  \textit{Introduction to elliptic curves and modular forms}, 2nd ed.,
Springer-Verlag, Berlin, Heidelberg, New York 1993.

\bibitem[Kol13]{Kol13}
J. Koll\'ar,
\textit{Singularities of the minimal model program},
With the collaboration of S\'andor Kov\'acs.,
Cambridge Tracts in Mathematics 200 (2013).

\bibitem[KM98]{KM98}
J. Koll\'ar, S. Mori,
\textit{Birational geometry of algebraic varieties},
With the collaboration of C. H. Clemens and A. Corti., 
Cambridge Tracts in Mathematics 134. (2008).


\bibitem[Kon02]{Kon02}
S. Kond\={o},
\textit{The moduli space of curves of genus $4$ and Deligne-Mostow's complex reflection groups}, 
Algebraic geometry 2000, Azumino (Hotaka), 383--400,
Adv. Stud. Pure Math., {\bf 36}, Math. Soc. Japan, Tokyo, 2002.

\bibitem[Kon07]{Kon07}
S. Kond\={o},
\textit{The moduli space of 5 points on $\P^1$ and K3 surfaces},
Arithmetic and geometry around hypergeometric functions, 189--206, Progr. Math., {\bf 260}, 
Birkh\"auser, Basel, 2007.

\bibitem[Kon13]{Kon13}
S. Kond\={o},
\textit{The Segre cubic and Borcherds products, Arithmetic and geometry of K3 surfaces and
Calabi-Yau threefolds},
549--565, Fields Inst. Commun., {\bf 67}, Springer, New York, 2013.

\bibitem[Loo85]{Loo85}
E. Looijenga,
\textit{Semi-toric partial compactifications I},
Report 8520 (1985), 72 pp., Catholic University Nijmegen.

\bibitem[Loo86]{Loo86}
E. Looijenga,
\textit{New compactifications of locally symmetric varieties},
Proceedings of the 1984 Vancouver conference in algebraic geometry, 341--364. 
CMS Conf. Proc., 6, AMS Province (1986).

\bibitem[Loo03a]{Loo03a}
E. Looijenga,
\textit{Compactifications defined by arrangements. I. The ball quotient case}, 
Duke Math. J., {\bf 118} (2003), no. 1, 151--187. 

\bibitem[Loo03b]{Loo03b}
E. Looijenga,
\textit{Compactifications defined by arrangements. II. Locally symmetric varieties of type IV}, 
Duke Math. J., {\bf 119} (2003), no. 3, 527--588. 

\bibitem[Loo23]{Loo23}
E. Looijenga,
\textit{A ball quotient parametrizing trigonal genus $4$ curves}, 
Nagoya Math. J., {\bf 254} (2024), 366--378.

\bibitem[LS07]{LS07}
E. Looijenga and R. Swierstra,
\textit{The period map for cubic threefolds},
Compos. Math., {\bf 143} (2007), no. 4, 1037--1049.

\bibitem[Ma19]{Ma19}
S. Ma, 
\textit{Quasi-pullback of Borcherds products},
Bull. Lond. Math. Soc. {\bf 51}, No. 6, 1061-1078 (2019).


\bibitem[MO23]{MO23}
Y. Maeda, Y. Odaka,
\textit{Fano Shimura varieties with mostly branched cusps},
Springer Proceedings in Mathematics \& Statistics(PROMS, volume 409), 2023, 633-664.

\bibitem[Mir81]{Mir81}
R. Miranda,
\textit{The moduli of Weierstrass fibrations over  $\P^1$}, 
Math. Ann., {\bf 255} (1981), no. 3, 379--394.

\bibitem[Mos86]{Mos86}
G.D. Mostow,
\textit{Generalized Picard lattices arising from half-integral conditions}, 
Publ. Math. IH\'ES {\bf 63} (1986) 91--106.

\bibitem[Nik80]{Nik80}  V. V.\ Nikulin, \textit{Integral symmetric bilinear forms and some of
their applications}, Math. USSR Izv., {\bf 14} (1980), 103--167.

\bibitem[Oda22]{Oda22}
Y. Odaka,
\textit{Semi-toric and toroidal compactifications as log minimal models, and applications to weak $K$-moduli},
arXiv:2203.09120.

\bibitem[Orl97]{Orl97}
D. Orlov,
\textit{Equivalences of derived categories and K3
 surfaces},
J. Math. Sci., {\bf 84} (1997), No. 5, 1361--1381 (1997).

\bibitem[Sch74]{Sch74}  B.\ Schoenberg,  \textit{Elliptic modular functions}, Springer-Verlag, Berlin, Heidelberg, New York 1974.

\bibitem[Thu98]{Thu98}
W. P. Thurston,
\textit{Shapes of polyhedra and triangulations of the sphere},
Geom. Topol. Monogr., {\bf 1} (1998), 511--549.
\end{thebibliography}
\end{document}